%% file: FSYr.tex
\documentclass[english,12pt]{article}

\usepackage[T1]{fontenc}
\usepackage{geometry}
\usepackage{babel}
\usepackage{float}
\usepackage{mathtools}
\usepackage{amsmath}
\usepackage{amsthm}
\usepackage{amssymb}
\usepackage{bbm}
\usepackage{color,colordvi}
\usepackage[mathscr]{eucal}

\numberwithin{equation}{section}

\makeatletter
\@ifundefined{date}{}{\date{}}
\makeatother

\usepackage{tikz}
\usepackage{tikz-cd}
\usepackage{tikzit}
\input{sample.tikzstyles}
\everymath{\displaystyle}
\usepackage{lmodern}
\usepackage{fancyhdr}
\usetikzlibrary{decorations,arrows}
\usetikzlibrary{decorations.markings}
\usetikzlibrary{decorations.pathmorphing}
\usetikzlibrary{matrix,arrows}

\makeatletter
\newcommand*{\relrelbarsep}{.386ex}
\newcommand*{\relrelbar}{%
  \mathrel{%
      \mathpalette\@relrelbar\relrelbarsep
        }%
        }
\newcommand*{\@relrelbar}[2]{%
  \raise#2\hbox to 0pt{$\m@th#1\relbar$\hss}%
    \lower#2\hbox{$\m@th#1\relbar$}%
    }
\providecommand*{\rightrightarrowsfill@}{%
      \arrowfill@\relrelbar\relrelbar\rightrightarrows
      }
\providecommand*{\leftleftarrowsfill@}{%
       \arrowfill@\leftleftarrows\relrelbar\relrelbar
        }
\providecommand*{\xrightrightarrows}[2][]{%
          \ext@arrow 0359\rightrightarrowsfill@{#1}{#2}%
          }
\providecommand*{\xleftleftarrows}[2][]{%
    \ext@arrow 3095\leftleftarrowsfill@{#1}{#2}%
    }
\makeatother

\theoremstyle{plain}
\newtheorem{thm}{Theorem}[section]
\newtheorem{cor}[thm]{Corollary}
\newtheorem{lem}[thm]{Lemma}
\newtheorem{prop}[thm]{Proposition}

\theoremstyle{definition}
\newtheorem{defi}[thm]{Definition}
\newtheorem{rem}[thm]{Remark}
\newtheorem{exa}[thm]{Example}

\def\act           {\,{\triangleright}\,}
\def\Act           {{\triangleright}}

\def\alph          {\alpha}
\def\Alpha         {\mathrm A}
\def\BB            {{\mathbb B}}
\def\bare          {bare}   
\def\bb            {b_{\scriptscriptstyle+}}  
\def\BC            {\Bic\C} 
\def\bd            {\mathfrak B}
\def\be            {\begin{equation}}
\def\bearl         {\begin{array}{l}}
\def\bearll        {\begin{array}{ll}}
\def\bet           {\beta}  
\def\Bic           {{\mathscr B}}
\def\bigboxtimes   {\mbox{\large$\boxtimes$}}
\def\Bl            {{\mathrm{Bl}}}
\def\bordoc        {\mathcal B\hspace*{-0.3pt}or\hspace*{-0.7pt}d^{\,\rm or}_{2,\rm o/c}}
\def\boti          {\,{\boxtimes}\,}
\def\BS            {\FF_\S}
\def\bv            {\mathbf B}
\def\bvb           {\bv_{\Bic}}
\def\bvcan         {\bv^{\rm can}}
\def\bvf           {\bv_{\Fr}}
\def\bvm           {\bv}
\def\bvo           {\bv^\circ}

\def\bvz           {\bv} 
\def\C             {{\ensuremath\calc}}
\def\cala          {{\mathcal A}}

\def\calb          {{\mathcal B}}
\def\calc          {{\mathcal C}}
\def\cald          {{\mathcal D}}
\def\calm          {{\mathcal M}}

\def\caln          {{\mathcal N}}
\def\calz          {{\mathcal Z}}
\def\cb            {\beta}      
\def\cB            {\gamma}     
\def\cir           {\,{\circ}\,}
\newcommand\Cite[2] {\cite[#1]{#2}}
\newcommand\coevl[1]{\mathrm{coev}^{\rm l}_{\!#1}}
\newcommand\coevr[1]{\mathrm{coev}^{\rm r}_{\!#1}}
\def\colB          {\mathrm{col}_\partial} 
\def\colE          {\mathrm{col}_E} 
\def\cancolor      {\varphi_{\rm can}}
\def\Colon         {:\quad}
\def\colorDefect   {purple!70}
\def\colV          {\mathrm{col}_V} 
\def\complex       {{\mathbbm C}}
\def\Cor           {\mathrm{Cor}}
\def\CorSN         {\mathrm{Cor}_{\rm SN}}
\def\CorSNb        {\overline{\mathrm{Cor}}_{\rm SN}}
\def\CorSNt        {\widetilde{\mathrm{Cor}}_{\rm SN}}
\def\Cup           {\mbox{\large$\cup$}}
\def\cws           {complemented world sheet} 
\def\cyl           {S^1{\times}I}
\def\cylS          {S{\times}I}
\newcommand\Cyl[1] {\mathcal C\hspace*{-1.1pt}y\hspace*{-0.2pt}l(\C,#1)}
\def\Cyle          {\Cyl{S^1}}
\newcommand\Cylo[1]{\mathcal C\hspace*{-1.1pt}y\hspace*{-0.2pt}l^{\circ\!}(\C,#1)}
\newcommand\Cylob[1]{\mathcal C\hspace*{-1.1pt}y\hspace*{-0.2pt}l^{\circ\!}(\Bic,#1)}
\def\Cyloe         {\Cylo{S^1}}
\newcommand\Cylof[1]{\mathcal C\hspace*{-1.1pt}y\hspace*{-0.2pt}l^{\circ\!}(\Fr,#1)}
\def\DC            {D_\C}
\def\dd            {d}
\def\DD            {{\mathbb D}}
\def\disk          {{\mathrm D}}
\def\dsty          {\displaystyle }
\def\ee            {\end{equation}}
\def\EE            {\widehat{E}}
\def\eear          {\end{array}}
\def\End           {\mathrm{End}}
\def\Enumerate     {\def\leftmargini{1.34em}~\\[-1.42em]\begin{enumerate}}
\def\eps           {\varepsilon}
\def\eq            {\,{=}\,}
\newcommand\evD[1] {\langle#1\rangle_\disk^{}}
\newcommand\evl[1] {\mathrm{ev}^{\rm l}_{\!#1}}
\newcommand\evr[1] {\mathrm{ev}^{\rm r}_{\!#1}}
\def\ffg           {full Frobenius graph}
\def\FF            {{\mathbb F}}
\def\FFt           {\widetilde{\mathbb F}}
\def\fg            {Frobenius graph}
\def\Fr            {{\mathscr F\hspace*{-1pt}r(\calc)}}
\def\Funre         {{\mathcal Rex}}
\def\Gama          {\varGamma}
\def\GAMA          {[\Gama]}
\def\GamaA         {\varGamma_{\!\tcl}}
\def\GamaAp        {\varGamma_{\!\tcl'}}
\def\GamaApp       {\varGamma_{\!\tcl''}}

\def\GamaS         {\varGamma_{\!\S}}
\def\GamaSt        {\varGamma_{\!\widetilde\S}}
\def\GC            {\Xi_\C}
\def\Gr            {{\mathrm G}}
\def\Hom           {\mathrm{Hom}}
\def\HomB          {1\mbox-\Hom_{\Bic}}
\def\HomC          {\ensuremath{\Hom_\calc}}
\def\HomCyle       {\Hom_{\mathcal C\hspace*{-0.5pt}yl(\C,S^1)}}
\def\HomCyloS      {\Hom_{\mathcal C\hspace*{-0.5pt}yl^{\circ\!}(\C,S)}}
\def\HomCylod      {\Hom_{\mathcal C\hspace*{-0.5pt}yl^\circ(\C,\partial\surf)}}
\def\HomCyloe      {\Hom_{\mathcal C\hspace*{-0.5pt}yl^\circ(\C,S^1)}}
\def\HomM          {\ensuremath{\Hom_\calm}}
\def\HomZ          {\ensuremath{\Hom_{\calz(\calc)}}}
\newcommand\hsp[1] {\mbox{\hspace{#1 em}}}
\def\I             {{\mathcal{I}}}
\def\id            {{\mathrm{id}}}
\def\Id            {{\mathrm{Id}}}
\def\Im            {{\mathrm{Im}}}
\def\iG            {\underline G{}}
\def\iHom          {\underline{\Hom}}
\def\iHomM         {\underline{\Hom}_\calm}
\def\imu           {\underline\mu}
\def\iN            {\,{\in}\,}
\def\iNat          {\underline{\Nat}}
\def\Itemize       {\def\leftmargini{1.05em}~\\[-1.65em]
                   \begin{itemize} \addtolength\itemsep{-7pt}}
\def\jj            {\jmath}
\def\kGr           {\ko\Gr}
\def\KK            {K}
\def\ko            {{\ensuremath{\Bbbk}}}
\def\la            {^{\rm l.a.}}
\def\LL            {L}
\def\M             {{\ensuremath\calm}}
\def\Map           {{\mathrm{Map}}}
\def\Mapr          {\widehat{\Map(\S)}}
\def\Maps          {\ensuremath{\Map(\S)}}
\def\Mod           {{\rm mod}\text-}
\def\moD           {\text-{\rm mod}}
\def\MoD           {\text-{\rm mod}\text-}
\def\mq            {\reflectbox{$?$}}
\def\muhor         {\underline\mu{}_{\rm hor}}
\def\muhore        {\muhor^{\rm l}}
\def\muhorz        {\muhor^{\rm r}}
\def\muver         {\underline\mu{}_{\rm ver}^{}}
\def\N             {{\ensuremath\caln}}
\def\Nat           {\mathrm{Nat}}
\def\Null          {{\mathrm N}}

\newcommand\Nxl[1] {\\[-1.3em]\\[#1mm]}
\def\ohr           {\reflectbox{$\rho$}}
\def\ohrp          {\reflectbox{$\rho$}\hspace*{-1pt}'}
\def\ohrs          {\reflectbox{$\scriptstyle\rho$}}
\def\OHZ           {{\mathscr O\HomZ}}
\def\one           {{\bf1}}

\def\opp           {^{\rm opp}}
\def\Opp           {{{}^{\phantom{t}\!\!\!}}\opp}
\def\OSN           {{\mathscr O\SN}}
\def\ota           {\otimes_{\!A}}
\def\otA           {\,{\ota}\,}
\def\otap          {\otimes_{\!A'}}
\def\otAp          {\,{\otap}\,}
\def\otb           {\otimes_{\!B}}
\def\otB           {\,{\otb}\,}
\def\oti           {\,{\otimes}\,}
\def\Oti           {{\otimes}}
\def\otik          {\,{\otimes_\ko}\,}
\newcommand\otx[1] {\,{\otimes_{\!#1}}\,}
\def\OWS           {{\mathscr O\mathrm{WS}_\C}}
\def\pcan          {p^{\rm can}}
\def\phiH          {\varphi^{}_{\mathrm{Hom}}}
\newcommand\pn[1]  {^{\epsilon_{#1}}}
\def\pop           {\surf_{\rm p.o.p.}}
\def\Prof          {\mathcal P\hspace*{-0.7pt}r\hspace*{-0.7pt}o\hspace*{-1.4pt}f_\ko}
\def\pz            {p} 
\newcommand\rarr[1]{\xrightarrow{~#1~}}
\newcommand\Rarr[1]{\,{\xrightarrow{\,#1\,}}\,}
\def\ra            {^{\rm r.a.}}
\def\rev           {^{\rm rev}}
\def\S             {\mathscr S}
\def\sew           {{\mathrm s}}

\def\Set           {\mathcal S\hspace*{-0.5pt}e\hspace*{-0.4pt}t}
\def\Shore         {\S_{\rm hor}^{\rm l}}
\def\Shorz         {\S_{\rm hor}^{\rm r}}
\def\SN            {\mathrm{SN}_\C^{}}
\def\SNo           {\mathrm{SN}_\C^\circ}
\def\SNob          {\mathrm{SN}_\Bic^\circ}
\def\SNof          {\mathrm{SN}_{\Fr}^\circ}
\def\SO            {\breve{\S}}
\def\surf          {{\varSigma}}
\def\Times         {\,{\times}\,}
\def\To            {\,{\to}\,}
\def\tcl           {\vartheta}  
\def\tr            {\mathrm{tr}}
\def\UC            {\mathrm U_\C}  
\def\UCorSN        {\mathrm{UCor}_{\rm SN}}
\def\vect          {\mathcal V\hspace*{-0.7pt}e\hspace*{-0.3pt}c\hspace*{-0.2pt}t}
\def\Vd            {V_\partial}
\def\Vee           {{}^{\vee\!}}
\newcommand\void[1] {}
\def\vphi          {\phi}    

\def\Zc            {\reflectbox{$\Zm$}}
\def\Zcs           {{\reflectbox{$\scriptstyle \Zm$}}}
\def\ZC            {{\ensuremath{\calz(\calc)}}}
\def\Zm            {{\mathrm Z}}

\begin{document}

\thispagestyle{empty}

~ \vskip 2.5em

\begin{center}
{\bf \Large String-net Construction of RCFT Correlators} 

\vskip 2.5em

{\large \  \ J\"urgen Fuchs\,$^{\,a},~\quad$ Christoph Schweigert\,$^{\,b}~\quad$
and $\quad$ Yang Yang\,$^{\,b}$ }

\vskip 15mm

 \it$^a$
 Teoretisk fysik, \ Karlstads Universitet\\
 Universitetsgatan 21, \ S\,--\,651\,88\, Karlstad
 \\[9pt]
 \it$^b$
 Fachbereich Mathematik, \ Universit\"at Hamburg\\
 Bereich Algebra und Zahlentheorie\\
 Bundesstra\ss e 55, \ D\,--\,20\,146\, Hamburg

\end{center}

\newpage
\tableofcontents{}
\newpage

\section*{Preface}

It is a common and often repeated saying that
for any quantum field theory the ultimate goal is to \emph{solve the theory}.
Alas, already specifying in concrete terms the very meaning of this goal is a
formidable task. After all, even the definition of the notion of a quantum 
field theory is in itself under debate. In the present book we take a pragmatic 
point of view: we understand `solving the theory' as the task to determine 
a `consistent' set of correlation functions.

For the class of two-dimensional rational conformal field theories this task has
indeed been formulated as a precise mathematical problem: to determine a solution
to a collection of consistency conditions for correlators as particular elements
of spaces of conformal blocks. These conditions, which impose invariance under
appropriate mapping class groups and compatibility with factorization, are known
explicitly. Moreover, already two decades ago a constructive solution to the problem 
has been achieved. It crucially utilizes tools from the three-dimensional surgery 
topological field theory associated with the modular fusion category that underlies a
chiral rational conformal field theory.

The purpose of the present book is to exhibit an alternative solution to the same
problem. Instead of relying on topological field theory methods, the novel 
approach is based on string-net models. This has the advantage of 
providing a purely two-dimensional construction of correlators. In particular we obtain 
concise geometric expressions for the objects that describe bulk and boundary fields,
comprising in fact more general classes of fields than what has been
standard in the literature (including earlier work of two of the present authors).
The expressions for the fields are phrased in terms of idempotents in the cylinder 
category of the underlying modular fusion category. Combining these idempotents with 
Frobenius graphs on the world sheet yields string nets that satisfy all the required 
consistency conditions and thus form a consistent system of correlators.

In our approach, correlators are the basic entities to be expounded in a quantum
field theory. In the physics literature on conformal field theory, operator products 
are abundant, though. An obvious question is thus to ask what results for operator
products can be extracted from our setup. When answering this question one must
account for the fact that to extract an operator product from the correlator for a
specific world sheet, it is necessary to endow the world sheet with additional 
geometric structure, specifically with a marking, a structure familiar from 
pair-of-pants decompositions. We show that with the help of this extra geometric
structure operator products can indeed be determined. In fact, we obtain natural 
algebraic expressions which still make sense beyond semisimplicity, i.e.\ beyond the
realm of rational conformal field theories. 

As an application of independent interest we derive an Eck\-mann-Hil\-ton relation
internal to a braided category, thereby demonstrating the utility of string nets for
understanding algebra in braided tensor categories. Finally, as an outlook we discuss
the notion of a universal correlator. This systematizes the treatment of situations
in which different world sheets have the same correlator. It also allows us to
introduce a more comprehensive mapping class group.

\medskip

To follow our arguments, actually no particular prerequisites from quantum field theory 
are required, provided the reader accepts our point of view that the relevant issue
is the mathematically well-formulated problem of finding consistent systems of 
correlators in the sense of Definition \ref{def:cosyCor}. In our opinion this 
attitude is fully justified by the beauty of the (algebraic) solution of this
problem. An exposition to a standard class in two-dimensional conformal field theory
(such as \cite{Sche2a}) could be helpful to understand the physics interpretation of
our results, but the more mathematically minded reader might rather wish to look up
such literature only after a first read of Chapter \ref{sec:corCFT}.
On the mathematical side, it is of avail to be acquainted with the basics of 
the theory of braided fusion categories. However, only a fraction of the 
material presented in e.g.\ \cite{EGno} will amply suffice, and in the Appendix
we supply a brief survey of all pertinent input. (Also, these days physicists working
on topological phases of matter will typically have the required background at their
disposal.) Moreover, it definitely helps the
intuition to know that the representation categories of a class of
vertex algebras have the structure of a modular tensor category \cite{huan21} 
and that they give rise to systems of conformal blocks.

\medskip

We thank Gregor Schaumann and Matthias Traube for helpful comments on an early version
of the text, and the anonymous referee for suggesting valuable improvements.
JF is supported by VR under project no.\ 2017-03836. CS and YY are supported by 
the Deutsche Forschungsgemeinschaft (DFG, German Research Foundation) under
SCHW1162/6-1; CS is also supported by DFG under Germany's Excellence Strategy - 
EXC 2121 ``Quantum Universe'' - 390833306.

\newpage

\section{Introduction}

The main topic of this book is a novel construction of the correlators of two-dimensional
rational conformal field theories on arbitrary world sheets, allowing in particular 
for physical boundaries as well as for defects and general defect junctions. A 
crucial tool is a modular functor defined in 
terms of string nets. We treat boundary insertions in such a way that the
open and closed sectors of the CFT can be dealt with uniformly, without the need to 
introduce string nets for surfaces with physical boundaries.
A different construction of RCFT correlators, relying instead on the modular functor 
that is furnished by three-dimensional topological field theories of 
Reshetikhin-Turaev type, has already been known for about two decades.
In this introductory chapter of the book we expose virtues and drawbacks of the 
two constructions.

\subsection{A New Construction of RCFT Correlators} %

``Solving a quantum field theory'' -- for any type of quantum field theory, in 
any number of dimensions -- is an important goal, with deep ramifications
in mathematics and physics. This goal has been pursued in a variety of approaches. In
the specific case of two-dimensional rational conformal field theories (RCFTs, in
short), the task can be reduced to a tractable mathematical
problem: correlators can be described as specific elements in spaces of conformal 
blocks, and these vector spaces are under control in terms of modular functors. The 
consistency conditions for these correlators have been discussed extensively in the
literature, starting with \cite{frsh2,sono2}, see e.g.\ \cite{kolR,fuSc27}. 
They amount to requiring that the correlator for a world sheet is invariant under 
the action of the mapping class group of the world sheet, that upon sewing of 
world sheets correlators are mapped to correlators, and that the basic two-point
correlators are non-vanishing. 

Making use of three-dimensional topological field 
theories of Reshetikhin-Turaev type, it has been shown two decades ago that solutions to
these constraints are in bijection to (Morita classes of) special symmetric Frobenius
algebras \cite{fuRs4,fjfrs2} in the modular fusion category \C\ that governs the 
monodromies of the conformal blocks.
Let us briefly recapitulate this ``TFT approach to RCFT correlators'', highlighting both
successes and drawbacks. In this approach the correlator for an oriented world sheet 
$\S$ is realized as an element of the vector space of conformal blocks for a double
cover $\widehat\surf_\S$ of the surface $\surf_\S$ that underlies $\S$. The crucial
idea is to use the Reshetikhin-Turaev topological field theory that is based on the 
modular fusion category \C\ together with the fact that any three-manifold with 
embedded Wilson lines whose boundary is $\widehat\surf_\S$ determines a vector in 
the space of conformal blocks for $\widehat\surf_\S$. Accordingly, a vital step in
the TFT construction is to determine, for given $\S$, 
a three-manifold $M_\S$ (called the connecting 
manifold) with boundary $\partial M_\S \eq \widehat\surf_\S$. In this
framework extensive existence and uniqueness results for correlators have been
established. Three-dimensional geometry is a central ingredient in the construction,
and occasionally it requires a certain amount of mental
gymnastics, such as in the proof of factorization \cite{fjfrs}.

Avoiding such an effort is a motivation to develop an approach to CFT correlators that
is inherently two-dimensional. But it is not the most relevant motivation. A more
significant one is the desire to transcend the restriction to semisimplicity
of the modular tensor category. This
restriction is inherent to the TFT approach: an extended three-dimensional topological
field theory with target the bicategory of finite tensor categories assigns to the 
circle a modular category which is necessarily semisimple \cite{bdsv3}. In contrast,
many conformal field theories of direct physical relevance are based on chiral data 
that are captured by a non-semisimple modular category. To understand such 
theories is a long term goal of our research. While the present paper does restrict 
to the case of semisimple modular categories, i.e.\ modular fusion categories,
the novel construction we present promises to admit an extension to non-semisimple
situations, in particular to modular finite tensor categories. 
This is e.g.\ supported by the existence of a state-sum construction for a
modular functor \cite{fScS4} that does not require semisimplicity. (The construction
in \cite{fScS4} is, however, inconvenient for the construction of correlators.)

The new construction produces the same correlators as the one in \cite{fuRs4} in all
cases in which a direct comparison is possible, in particular in the absence of
defect lines. In that case, by the uniqueness result of \cite{fjfrs2}, the fact
that the two-point correlators on the sphere and on the disk
coincide implies coincidence of \emph{all} correlators. Crucially, besides being
considerably simpler than the TFT construction, our new approach is entirely
two-dimensional and thereby avoids the main obstacle to an extension beyond
semisimplicity. In fact, basic ingredients of the construction are formulated
in such a manner that they still make sense in the non-semisimple case.


\subsection{Advantages of the String-Net Construction} %

Any construction of correlators requires a thorough understanding of the field
content. Accordingly, one starting point of our construction is the systematic
mathematical formalization of defect fields and of their operator products
in (not necessarily semisimple) conformal field theories that was developed 
in \cite{fuSc27}. As an essential tool this formalization uses the theory of internal 
natural transformations, as introduced in \cite{fuSc25} within the framework of 
modular finite tensor categories and pivotal module categories over them. This 
framework is convenient even in the semisimple case that applies to rational CFT. In
the present paper we actually transcend the setting of \cite{fuSc27} 
(as well as the one of \cite{fuRs4})
in that we allow for more general defect fields, at which 
any finite number (rather than only two) of
defect lines can meet, 
and similarly for more general boundary fields. (See the illustrations in Examples
\ref{exa:X.X.X} and \ref{exa:MXXN} for a first impression, and subsections
\ref{sec:bulkfields} and \ref{sec:bdyfields} for details.)

To get control over correlators for such general defect and boundary fields in the
purely algebraic framework of internal natural transformations \cite{fuSc27} is hard, 
if not impossible, in particular when it comes to their invariance under the relevant 
mapping class group. More generally, beyond semisimplicity full control over correlators
has been achieved so far only in limited cases, e.g.\ -- by combining the results of
\cite{fuSs3} and \cite{fuSc22} -- for bulk fields on oriented world sheets in the case 
that the modular finite tensor category admits a description as a category of modules
over a factorizable ribbon Hopf algebra.

The method we develop in the present paper is based on the \emph{string-net 
construction} of modular functors. It is known \cite{kirI24,goos,baGo} that this 
construction, which takes as an input a spherical fusion category, realizes a modular
functor of Turaev-Viro type. The string-net approach exhibits several crucial 
advantages. First of all, in contrast to Reshetikhin-Turaev type constructions (see
\cite{lyub11} for the non-semisimple case), the representations of mapping class
groups are geometric. Also, string nets provide a geometric realization of specific 
vectors in the spaces of conformal blocks. Moreover, in contrast to approaches
based on a Le\-go-Teich\-m\"ul\-ler game (see e.g.\ \cite{fuSc22}), the string-net
construction does not introduce any auxiliary structure like a pair-of-pants 
decomposition whose combinatorics quickly becomes unwieldy, in particular in the
presence of boundaries and of additional stratifications of the surface that result
from defects. Accordingly we expect that the string-net construction can also be
generalized beyond oriented surfaces, e.g.\ to
unoriented surfaces or to surfaces with spin structure.

The advantages of string nets
have already been exploited in \cite{scYa} to construct correlators of 
bulk fields in the special case -- called the \emph{Cardy case} in the CFT literature 
-- that the boundary conditions of the full conformal field theory are just the 
objects of the modular fusion category \C\ of chiral data. The present paper extends 
the treatment of \cite{scYa} in two directions: beyond the Cardy case, allowing for
general local field theories with chiral data given by \C\ (in which case the boundary 
conditions are objects of a module category over \C), and to include (generalized)
boundary and defect fields. It is worth stressing that we \emph{construct} these 
fields, i.e.\ do not use (as has been done in \cite{traub}) the existence of a 
consistent system of bulk and boundary fields as an input. We construct the fields 
as objects in the Drinfeld center \ZC\ of \C. More specifically, we realize them in 
terms of idempotents in the category $\Cyloe$ of boundary values on the circle, whose 
Karoubification -- the \emph{cylinder category} $\Cyle$, see Definition \ref{def:Cyl}
-- is known, for \C\ a pivotal fusion category, to be equivalent to the Drinfeld center. 
(This constitutes in fact a more categorical variant of the tube algebra, whose 
irreducible representations give the simple objects of the Drinfeld center.)

It is worth stressing that our results cover all rational conformal field theories --
that we construct fields in terms of objects 
in a Drinfeld center by no means implies that we deal only with chiral data that are
described by a modular tensor category which happens to be a Drinfeld center. Rather, 
we use the fact that for any modular tensor category \C, the Drinfeld center \ZC\
is braided equivalent to the Deligne product $C\boti\C\rev$, which is the category
that describes, in physics parlance, the ``combinations of left and right movers.''

We treat field insertions on the boundary as 
objects in the Drinfeld center as well. To this end we use the adjoint of the
forgetful functor $\ZC\,{\to}\,\C$. (As explained in Remark \ref{rem:U.L}, this can
be motivated from ideas in factorization homology, which suggest that the embedding 
of an interval into a circle gives rise to this functor.) At a technical level, this
is possible because of the cyclicity properties expressed in Lemma \ref{lem:Lxy=Lyx}.
Owing to this novel treatment of boundary insertions, we can treat the open and
closed sectors of the CFT uniformly in terms of string nets,
even though we do not define string nets for surfaces with physical boundaries
on which boundary conditions need to be prescribed.


\subsection{Organization of the Book} %

The rest of this 
                  book  %
is organized as follows. In 
                  Chapter  %
\ref{sec:corCFT} we 
first explain how the problem of finding correlation functions is translated to the
problem of finding specific elements in spaces of conformal blocks. These spaces of
conformal blocks are then described in terms of an open-closed modular functor, a
notion that is introduced in Section \ref{sec:ocmf}. Afterwards we present, in 
Section \ref{sec:S}, the class of world sheets to which we want to associate 
correlators. Concretely, we wish to work on two-dimensional oriented manifolds with
stratifications, so as to be able to account for topological line defects as well as
boundaries. The significance of studying quantum field theory on stratified spaces 
is by now well established, not only in view of their role in the study \cite{ffrs5}
of symmetries and dualities. Given a modular fusion category \C, the cells of the
world sheet need representation-theoretic labels: two-cells are labeled by simple
special symmetric Frobenius algebras in \C, labels for topological defect lines 
and for boundary conditions are bimodules and modules, respectively, internal to the
category \C. This allows us to discuss, in the subsequent Section \ref{sec:fields},
the field content, that is, the different types of field insertions in the bulk and
on the boundary.  Finally, in Section \ref{sec:Bl-Cor} we discuss correlators as
vectors in spaces of conformal blocks and give a precise mathematical meaning to
factorization. This culminates in Definition \ref{def:cosyCor} of a \emph{consistent
system of correlators}. At this point we have set up the problem that we are going to solve.

                  Chapter  %
\ref{sec:the construction} is devoted to the solution of this problem.
Our main tool are string nets; they are introduced in Section \ref{sec:SN}. 
To account for the gluing boundaries (as opposed to physical boundaries)
of world sheets, we need to deal with string nets on surfaces with boundaries.
In Section \ref{sec:ipoco} we discuss categories of boundary values and string-net
spaces for given boundary values. This allows us to set up, in Section \ref{sec:SN-mofu},
a modular functor in terms of string nets. In Section \ref{sec:SNws} we then explain
how a world sheet (including the decorations of its two-, one- and zero-cells) 
determines a string net and thereby a vector in the space of conformal blocks. Readers
who want to get a quick impression of this process are invited to see how the local
situation described in the picture \eqref{eq:pic:Ftb=Fb}, which shows a world sheet
containing two defect lines and a boundary circle that supports two boundary field
insertions as well as two physical boundaries, is translated into the idempotent 
\eqref{eq:pic:p4Ftb=Fb} of the cylinder category. This determines the form of the
string net close to the boundary of a world sheet. To also take care of the two-cells
of a  world sheet we equip them with what we call a \ffg\ -- a graph whose edges are
labeled by the Frobenius algebra assigned to the two-cell and whose vertices are 
labeled by (multiple) products and coproducts (see Definition \ref{defi:ffgraph}). The
particular choice of \ffg\ is irrelevant, which is yet another incarnation of the 
well-known compatibility of two-dimensional Pachner moves and Frobenius algebras. This
feature is also at the basis of the invariance of the correlators under the mapping
class group. The resulting prescription for the string-net correlators is summarized
in Definition \ref{def:CorSN}. Section \ref{sec:the construction} culminates in Theorem
\ref{thm:correlators}, which asserts that our purely two-dimensional construction 
indeed provides a consistent system of correlators.

                  Chapter  %
\ref{sec:specos} is devoted to the aspects of a few particular correlators and
       specifically   %
to their relation with operator products. The goal is to extract
from specific string-net correlators composition morphisms or, in other words,
products, on the field objects. This
       translation    %
from string nets to elements of morphism spaces of the relevant category requires
additional choices, concretely a marking (akin to a pair-of-pants decomposition)
of the world sheet. In Section \ref{sec:vertical} we show how the composition of
     defect  %
fields on the same line defect (as portrayed in picture \eqref{eq:pic:2.24})
leads to the vertical product of relative natural transformations as introduced
in \cite{fuSc25}. (The proposal, put forward in \cite{fuSc27}, that the vertical 
product describes this operator product, applies to all finite -- not necessarily
semisimple -- conformal field theories; our string-net techniques at present only
settle the case of semisimple categories, though.) 
In Section \ref{sec:horpros} we consider horizontal operator products, which are 
related to the fusion of two parallel topological line defects (compare the picture
\eqref{eq:pic:2.25}). Here the situation is considerably more involved than for
the vertical product: as we work with braided categories, there are two different
horizontal products, corresponding to two non-isotopic string nets (with equal
markings -- or, equivalently, isotopic string nets with different markings). 
An immediate insight is that the two different horizontal products are related
by a braiding, see Remark \ref{rem:4.2}. In the two subsequent
sections we compute
the bulk algebra and the torus partition function, finding full agreement with results
in the literature. We finally show that the operator product of two boundary fields
is given by the natural composition of internal Hom objects (Section \ref{sec:boOPE}) 
and that the bulk-boundary operator products are the structure morphisms of the end
that defines a bulk or, more generally, a disorder field (Section \ref{sec:buboOPE}).

As is already apparent from the discussion of the two different horizontal products, 
string nets are a convincing tool for understanding algebraic facts in a braided 
monoidal setting. A further illustration of this point is given in 
      Chapter  %
\ref{sec:brEHr}, where we establish 
      an internalized  %
version of an Eckmann-Hilton theorem.

It is known (see e.g.\ \cite{ffrs5}) that different world sheets can have correlators
that are closely related. In particular there are identities of correlators
which pinpoint symmetries of the conformal field theory that can be extracted from 
the presence of invertible topological defects. In 
      Chapter  %
\ref{sec:uCorr} we introduce 
the new notion of a \emph{universal correlator} which captures such relations. As it
turns out, the most natural notion of mapping class group under which correlators are
invariant is the stabilizer group of the world sheet regarded as a vector in the string-net
space for a pivotal bicategory of module categories over the fusion category \C.

In nine appendices we collect mathematical background and fix notation. The background
information concerns algebra and representation theory in finite tensor categories and
features of their Drinfeld centers. The expert reader might wish to consult these 
appendices only when needed.

 \medskip

Much work remains to be done. In particular it is a challenge to extend the string-net
construction beyond spherical fusion categories to finite tensor categories
that are not semisimple. In this respect it is gratifying that the methods of
\cite{fuSc25,fuSc27} work for general finite conformal field theories. Furthermore,
in the end even rigidity should be relaxed to the more general structure of a 
Grothendieck-Verdier duality, which is the natural duality structure for large 
classes of vertex operator algebras \cite{alsW}. It is encouraging that for ribbon
Grothendieck-Verdier categories one can still construct \cite{muWo2} blocks that 
carry a representation of a handle body group.

  ~ %

\section{Correlators in Rational Conformal Field Theory} \label{sec:corCFT}  %

In this chapter we turn the
physics
task of finding correlation
functions into the mathematically precise problem of determining specific
elements in spaces of conformal blocks. These spaces are described in terms 
of an open-closed modular functor. Section \ref{sec:S} presents a
comprehensive class of world sheets for which correlators will be 
constructed, allowing for topological line defects as well as
boundaries. The cells of these stratified manifolds are labeled 
by representation-theoretic data internal to a modular fusion category \C: 
two-cells by simple special symmetric Frobenius algebras in \C, and defect 
lines and boundary conditions by bimodules and modules, respectively. Within
this setup, Section \ref{sec:fields} provides the field content, i.e.\ the 
different types of field insertions in the bulk and on the boundary. The
field insertions we allow are
       considerably 
more general than what is usually 
discussed in the literature. Finally, Section \ref{sec:Bl-Cor} expresses 
correlators as vectors in spaces of conformal blocks and culminates in 
Definition \ref{def:cosyCor} of a \emph{consistent system of correlators}.
Finding such systems for a given modular functor can be seen as a precise 
mathematical realization of `solving a conformal field theory.'

\subsection{Open-closed Modular Functor} \label{sec:ocmf}  %

A fundamental task in any quantum field theory is to determine, for each member of
a class of manifolds with suitable additional structure -- such as
boundary conditions and field insertions -- a correlation function.
In this sense, to specify a quantum field theory one needs to provide first its 
field content and then an infinite set of correlation functions.
Specifically, in two-dimensional conformal field theory (CFT) one wants to associate 
such a correlation function to a surface with additional data, which we call a 
\emph{world sheet}. 

In the approach to CFT correlation functions we take in this paper, the starting point 
is not a Lagrangian or a partial differential equation describing classical field 
equations, and we do not rely on any type of perturbation theory. Instead, we postulate
that the correlation functions of the theory exhibit certain chiral symmetries, which 
means concretely that they are solutions to a collection of linear differential 
equations, also known as chiral Ward identities. For the two-dimensional conformal 
field theories of our interest, these symmetries can be encoded in the structure of
a conformal vertex algebra, and the chiral Ward identities on any surface can be 
derived from that algebra. Locally, the solutions to these equations -- which contain the
correlation functions we are looking for as particular elements -- form vector spaces. 
For the models of our interest these vector spaces are finite-dimensional; 
when regarded as elements of these spaces, we refer to the correlation functions 
for brevity also as \emph{correlators}. Typically the solutions are multivalued,
so the fundamental groups of the relevant parameter spaces (conformal structure of the
surfaces and positions of insertion points), i.e.\ mapping class groups of surfaces,
act on the solution spaces.
The particular elements in these spaces that are correlation functions of 
\emph{bulk fields} must be single-valued, and thus transform trivially under the 
action of the mapping class group, while correlation functions involving general
\emph{defect fields} are invariant under a subgroup of the mapping class group. (For
the latter, see Definition \ref{def:Maps} for an initial description and 
Eq.\ \eqref{eq:def:Mapr} for an extended definition based on the notion of
``universal correlators'' which we introduce in 
                  Chapter  %
\ref{sec:uCorr}.)
To appreciate these statements in our context, several ingredients need to be 
explained. First and foremost, we need a tractable notion of \emph{conformal blocks}
and a precise notion of world
sheet. These will be provided in the present and the next subsection, respectively.

The spaces of solutions to the chiral Ward identities form vector bundles with 
projectively flat connection over moduli spaces of insertion points and conformal 
structures (see e.g.\ \cite{dagT}). One legitimate use of the term conformal blocks 
is to denominate these vector bundles, or also holomorphic sections in them. On 
the other hand, the conformal blocks depend functorially
on the representations of the chiral symmetry structure, the conformal vertex algebra, 
at the insertion points. This leads to the notion of a \emph{modular functor}, which
furnishes a more algebraic formalization of the idea of conformal blocks: 
finite-dimensional vector spaces endowed with a representation of the mapping class
group and depending functorially on the representation-theoretic data at the insertions.
Indeed it is widely believed that one can associate to any rational chiral
CFT a modular functor 
such that these mapping class group representations are related to the projectively 
flat connection on the vector bundles of the CFT by a Riemann-Hilbert correspondence. 
Accordingly, another use of the term conformal blocks refers to the vector spaces
furnished by the modular functor; this is the use of the term that we adopt here.

A modular functor is then a certain 2-functor from a geometric bicategory of 
two-dimensional bordisms to an algebraic bicategory. For being relevant to CFT with 
the type of world sheets considered in this paper, which are allowed to have boundary
field insertions, the bordism category has to be of
\emph{open-closed} type, meaning that its objects are finite disjoint unions of both 
circles and intervals, while on the algebraic side it turns out to be
convenient to take the 1-morphisms as profunctors. This leads us to the following
definition, where \ko\ is an algebraically closed field of characteristic zero,
chosen once and for all. 

\begin{defi} \label{def:mofu}
An \emph{open-closed modular functor} is a symmetric monoidal 2-functor
  \be
  \Bl \Colon \bordoc \rarr~ \Prof
  \ee
from the bicategory $\bordoc$ of two-dimensional open-closed bordisms to the
bicategory of \ko-linear profunctors.
\end{defi}

Let us unwrap this definition. First, the bicategory $\bordoc$ of two-dimensional 
\emph{open-closed bordisms} is the following symmetric monoidal bicategory: An object
$\alph \iN \bordoc$ is a finite disjoint union of copies of the standard closed interval
$I \eq [0,1] \,{\subset}\, \mathbb R$, oriented from 1 to 0, and of the standard circle 
$S^1 \eq \{z\iN\complex \,|\, |z|\eq 1\} \,{\subset}\, \complex$, with counter-clockwise
orientation. A 1-morphism $\alph \,{\to}\, \alph'$ is a compact oriented
smooth surface $\surf$ with boundary (and with induced orientation of $\partial\surf$), 
together with a \emph{parametrization} $\vphi \,{\equiv}\, \{\vphi_-,\vphi_+\}$ of a
subset of $\partial\surf$; the latter consists of an orientation reversing embedding 
$\vphi_-(\surf)\colon \alph \,{\hookrightarrow}\, \partial\surf$ and an orientation
preserving embedding $\vphi_+(\surf)\colon \alph' \,{\hookrightarrow}\, \partial\surf$
whose images are disjoint.
A 2-morphism $\gamma\colon \surf \,{\to}\, \surf'$ is an isotopy class of diffeomorphisms 
from $\surf$ to $\surf'$ such that $\vphi_\pm(\surf') \eq \gamma \cir \vphi_\pm(\surf)$.
The composition of two 1-morphisms $\surf\colon \alph \,{\to}\, \alph'$ and
$\surf'\colon \alph' \,{\to}\, \alph''$ is gluing along $\alpha'$ (or rather,
as usual, along a collar of $\alpha'$, such that the 
glued surface is again smooth), the vertical composition of 2-morphisms is induced
by the composition of diffeomorphisms, and the horizontal composition of 2-morphisms is
gluing. The tensor product is given by disjoint union, with the monoidal unit the
empty set.

The bicategory $\Prof$ of \ko-\emph{linear profunctors} is a finite \ko-linear version 
of the standard \Cite{Prop.\,7.8.2}{BOrc1} bicategory of profunctors: An object of 
$\Prof$ is a \ko-linear finite abelian category $\mathcal A$. A 1-morphism 
$\mathcal A \,{\to}\, \mathcal A'$ is a left exact functor
$F\colon \mathcal A\opp \boti \mathcal A' \,{\to}\, \vect_\ko$ from the 
Deligne product of the opposite of $\mathcal A$ and of $\mathcal A'$ to the
category of finite-dimensional \ko-vector spaces. A 2-morphism $F \,{\to}\, F'$
in $\Prof$ is a natural transformation from $F$ to $F'$.
The composition of 1-morphisms $F\colon \mathcal A \,{\to}\, \mathcal A'$ and
$F'\colon \mathcal A' \,{\to}\, \mathcal A''$ is given by the coend
  \be
  (F \cir F')(\overline? \boti \mq) \,= \int^{X\in \mathcal A'}\!\!\!
  F(\overline? \boti X) \otimes_\ko F'(\overline X \boti \mq) \,,
  \ee
where $\overline X$ stands for $X\iN\mathcal A'$ regarded as an object of the 
opposite category.

\begin{rem} \label{rem:collars}
For defining the horizontal composition of 1-morphisms in $\bordoc$, the gluing must be
performed with a choice of collars \Cite{Thm.\,1.3.12}{KOck}. More specifically, one can
select collars for all composable pairs of bordisms. Different 
choices give different smooth structures on the glued cobordisms which are, 
however, all diffeomorphic. These diffeomorphisms are not canonical, but they are
all in the same isotopy class and hence descend to a unique 2-morphism in $\bordoc$.
As a consequence, the choice of collars is inessential. For more details see Section
3.1.2 of \cite{schom1}.
Also other differential-topological issues are, by standard arguments, insignificant 
for our construction.
\\
In fact, the situation for the target bicategory $\Prof$ is analogous: the 
composition of 1-mor\-phisms is defined up to a contractible choice of coends.
\end{rem}

In the sequel we tacitly strictify each of the bicategories $\bordoc$ and $\Prof$
      (as bicategories),
i.e.\ consider them as strict 2-categories. We then also take the 2-functor $\Bl$ as 
strict, and thus require strict preservation of composition.

That $\Bl$ is a 2-functor means in particular that
any gluing from $\surf\colon \alph \,{\to}\, \alph'$ and
$\surf'\colon \alph' \,{\to}\, \alph''$ along the one-manifold $\alpha'$ to the 
surface $\surf \,{\cup_{\alph'}}\, \surf'$ gives rise to a natural transformation
  \be
  \Bl(\surf)(\overline? \boti Y) \otimes_\ko \Bl(\surf')(\overline Y \boti \mq)
  \rarr~ \Bl(\surf \,{\cup_{\alph'}}\, \surf')(\overline? \boti \mq) 
  \label{eq:gluingBl}
  \ee
for $Y \iN \Bl(\alpha')$, given by the dinatural structure morphisms of the
(left exact) coend
  \be
  \Bl(\surf \,{\cup_\alph} \surf')(\overline? \boti \mq)
  \,= \int^{X\in \Bl(\alpha')}\!
  \Bl(\surf)(\overline? \boti X) \otimes_\ko \Bl(\surf')(\overline X \boti \mq) \,.
  \label{eq:coendBl}
  \ee
Moreover, for any 1-morphism $\surf$ in $\bordoc$, the functor $\Bl(\surf)$ carries 
an action of the mapping class group $\Map(\surf)$, i.e.\ the group of 2-endomorphisms of 
the bordism $\surf$, by natural endotransformations.
(Later on, we will be interested in the mapping class group $\Maps$ of a \emph{world sheet}
$\S$, as opposed to the mapping class group $\Map(\surf_\S)$ of the bordism that
underlies $\S$. This group will be introduced in Definition \ref{def:Maps}.)

\medskip

To be relevant for applications to rational conformal field theory, such a modular 
functor has to satisfy constraints. One first formalizes the symmetries of the CFT in 
the structure
of a vertex operator algebra $\mathfrak V$. Then one selects a good category \C\ of
$\mathfrak V$-representations, e.g.\ by requiring that the HLZ theory of logarithmic
tensor products \cite{hulz4} should apply, which endows \C\ with the structure of a 
\ko-linear braided monoidal category with a Grothendieck-Verdier duality. For a
\emph{rational} CFT -- the case of interest in the present paper -- the vertex operator 
algebra $\mathfrak V$ is required to be rational, meaning that \C\ is even a modular
fusion category (i.e., a semisimple modular tensor category); we fix this category \C\
throughout the paper.
We thus assume that there is an open-closed modular functor that models the conformal 
blocks of a rational conformal field theory with given modular fusion category \C;
we denote this modular functor by $\Bl_\C$.

Bulk fields of the CFT should be objects in the Deligne product $\C\rev\boti\C$ 
of the modular fusion category \C\ with its reverse $\C\rev$,
i.e.\ the same spherical fusion category but endowed with the opposite braiding. This 
captures the informal statement that bulk fields are obtained by `combining left
movers and right movers'. Remarkably, the nondegeneracy of the braiding of \C\ implies
that the Deligne product $\C\rev\boti \C$ is equivalent, as a modular category, 
to the Drinfeld center \ZC\ of the spherical fusion category underlying the modular
category \C. Accordingly we treat bulk fields as objects in \ZC.
(For the precise notion of bulk fields, see Section \ref{sec:bulkfields} below.)

We are interested in this paper in a modular functor which when restricted to the
\emph{closed sector} reproduces the known description of the conformal blocks of
rational conformal field theories. This includes in particular the requirement that
the modular functor $\Bl_\C$ assigns to the circle $S^1$ a category braided equivalent
to the Drinfeld center \ZC, for \C\ a \ko-linear spherical fusion category; in this
paper, \C\ is the modular fusion category that we have already fixed above. More
generally, we have to demand that the restriction of $\Bl_\C$ to the 
sub-bicategory of $\bordoc$ whose objects are disjoint unions of circles is isomorphic 
to the modular functor obtained by the Turaev-Viro state-sum construction for \C\
-- or, equivalently, by the Reshetikhin-Turaev surgery construction for the Drinfeld
center \ZC\ of \C. (The latter reformulation
implements the fact \Cite{Ch.\,17}{TUvi} that the Turaev-Viro construction based on \C\
and the Reshetikhin-Turaev construction based on \ZC\ are
isomorphic as once-extended topological field theories.)
For pertinent information on fusion categories and their Drinfeld center we refer to
\cite{EGno,TUvi} and to Appendices \ref{app:C} and \ref{app:Z}.
Without loss of generality we take the pivotal structure on \C\ to be strictified.

 \medskip

For ease of notation, in the sequel we write
  \be
  \Bl_\C(\surf)(\overline? \boti \mq) =: \Bl_\C(\surf;?;\mq) \,.
  \ee
We also follow the common terminology to call, for a 1-morphism 
$\surf\colon \alph \,{\to}\, \alph'$, the circles and intervals in $\partial\surf$
which are images under $\vphi_-(\surf)$ of the connected components of $\alph$ 
\emph{incoming}, and those which are images under $\vphi_+(\surf)$ of 
the connected components of $\alph'$
\emph{outgoing}; in case that $\alph'$ is the empty set (regarded as a one-manifold),
we briefly write $\Bl_\C(\surf;?)$ instead of $\Bl_\C(\surf;?;\emptyset)$.
For the closed sector, $\Bl_\C$ is the modular functor based on the spherical fusion 
category \C\ that is defined in terms of \emph{string nets} colored by \C. This
functor is implicit in the constructions in \cite{kirI24}; for the pertinent details
about the string-net construction see Section \ref{sec:SN}.
 
Note that the string-net construction gives rise to a 3-2-1 topological field theory 
(TFT) that is isomorphic to the 3-2-1-extended Turaev-Viro state sum TFT for \C\
\cite{kirI24,goos}. This isomorphism of topological field theories restricts to an
isomorphism between the respective modular functors. 
(In Section \ref{sec:SN-mofu} we will demonstrate explicitly that this is indeed 
true for the string-net modular functor $\SN$.)
Note that it follows e.g.\ that when the 1-morphism $\surf$ is a connected surface of 
genus $g$ with $\partial_-\surf \,{\cong}\, \mbox{\large$\sqcup$}_{i=1}^{\,p} S^1$ and
$\partial_+\surf \,{\cong}\, \mbox{\large$\sqcup$}_{j=1}^{\,q} S^1$, then there is
an isomorphism \cite{lyub11,BAki}
  \be
  \Bl_\C\big( \surf;\bigboxtimes_{i=1}^{\,p} X_i;\bigboxtimes_{j=1}^{\,q} Y_j \big)
  \,\cong\, \HomZ\big( \mbox{\Large$\otimes$}_{i=1}^{\,p} \widetilde X_i \,,
  \mbox{\Large$\otimes$}_{j=1}^{\,q} \widetilde Y_j \oti \KK^{\otimes g} \big) 
  \label{eq:Bl=HomZ}
  \ee
for $X_i,Y_j \iN \Bl_\C(S^1)$, where $\widetilde X_i$ and $\widetilde Y_j$ are their
images under $\Bl_\C(S^1) \Rarr\simeq \ZC$ and 
$\KK \iN \ZC$ is the distinguished Hopf algebra object
  \be
  \KK := \Zm_\ZC(\one_\ZC) = \int^{X\in\ZC}\!\! X^\vee \oti X 
  \label{eq:def:KK}
  \ee
in \ZC, i.e.\ the image of the monoidal unit under the central monad $\Zm_\ZC$ of \ZC.
The isomorphism \eqref{eq:Bl=HomZ} is far from canonical, though. For a fixed choice
of equivalence $\Bl_\C(S^1) \Rarr\simeq \ZC$ it can, however, be specified uniquely
up to unique isomorphism once the surface $\surf$ is endowed with the auxiliary 
structure of a \emph{marking}, which keeps track of the way that $\surf$ can be 
assembled by sewing pairs of pants \cite{lyub11,BAki,fuSc22}.

 \medskip

Next we extend $\Bl_\C$ from the closed sector to all of $\bordoc$. We
implement this extension with the help of the functor
  \be
  \LL := U\la \Colon \C \rarr~ \ZC
  \label{eq:def:LL}
  \ee
that is left adjoint to the forgetful functor from the Drinfeld center \ZC\ to \C.
Note that every functor between finitely semisimple categories is exact and hence
has both left and right adjoints. (Also, if \C\ is a finite tensor category, then
the forgetful functor $U$ is exact so that left and right adjoints of 
$U$ still exist, and they can be explicitly constructed via the
central monad and comonad, see Appendix \ref{app:ZZ}.)

In particular we require that the category assigned to the interval $I \iN \bordoc$ is,
up to canonical equivalence, $\Bl_\C(I) \,{=}\, \C$.
However, unlike in other approaches in the literature, we complement this 
assignment by a prescription that associates, via the functor \eqref{eq:def:LL}, 
a single category with a geometric boundary circle of a 1-morphism
that contains any finite number of embedded intervals, i.e.\ to a circle of the type
  \be 
  \scalebox{0.9}{\input{OC0.tikz}} 
  \label{eq:pic:circlewithIs}
  \ee
(Note that while \eqref{eq:pic:circlewithIs} is a boundary component of the 
underlying surface of a bordism, it is \emph{not} an object of the bordism
category $\bordoc$ we are working with.)
The string-net open-closed modular functor $\SN$ that we will obtain in Section 
\ref{sec:SN-mofu} indeed realizes the extension from the closed sector to
all of $\bordoc$ in this way.

More explicitly, for a 
geometric boundary circle containing the image under the parametrization map
of the union $I \,{\sqcup}\, I \,{\sqcup}\cdots{\sqcup}\, I$ of $n$ (ordered) copies 
of the interval we deal with $n$ copies of \C, and we combine any ordered collection 
$(X_1,X_2,...\,,X_n) \iN \C^{\times n}$ of objects in \C\ to the object
  \be
  \LL(X_1 \oti X_2 \,{\otimes} \cdots {\otimes}\, X_n) \in \ZC
  \label{eq:L(X1...Xn)}
  \ee
in the Drinfeld center. An advantage 
     of
this prescription is that it allows us to 
treat the open and closed sectors uniformly: we assign to \emph{any} connected 
component of $\partial\surf$ one and the same category, namely (up to isomorphism)
the Drinfeld center \ZC. This has various technical benefits. It will in particular
make it possible to define CFT correlators with the help of the string-net
construction even for world sheets with boundary field insertions.

\begin{rem} \label{rem:U.L}
The composite functor $U \cir \LL$ is the underlying functor of the \emph{central monad}
$\Zm$ of \C, as described in Appendix \ref{app:ZZ}. Accordingly, for any $X \iN \C$ the
object underlying the object $\LL(X)$ in \ZC\ is the coend given in \eqref{eq:def:Zm,Zc}.
In particular, if \C\ is semisimple, as considered here, we have
  \be
  U \circ \LL(X) = \bigoplus_{i \iN \I(\C)} i^\vee \oti X \oti i
  \ee
with $\I(\C)$ the chosen set of representatives for the isomorphism classes of simple 
objects of \C.
\end{rem}

The rationale behind the prescription \eqref{eq:L(X1...Xn)} is as follows. When 
the assignment of categories to one-manifolds comes from factorization homology --
which is indeed the case for the modular functor that is furnished by the string-net
construction -- then the embedding $I \,{\hookrightarrow}\, S^1$ induces a functor, 
calculated via excision, from $\Bl_\C(I) \,{\simeq}\; \C$ to $\Bl_\C(S^1)$, i.e.\
to the relative Deligne product $\C \,{\boxtimes_{\C\boxtimes \C\Opp}}\, \C$.
This functor factorizes into the composite of the functor
  \be
  \begin{array}{rl}
  \C \!\! & \rarr~\, \C \boti \C \,,
  \Nxl2
  c \!\! & \xmapsto{~~~}\, c \boti \one
  \eear
  \ee
to the ordinary Deligne product and the structure functor
$\C \boti \C \Rarr{} \C \,{\boxtimes_{\C\boxtimes \C\Opp}}\, \C$. By Proposition 2.18(i)
of \cite{fScS} one can represent the relative Deligne product as the twisted center
$\calz^{\mathrm D} (\C)$, i.e.\ the Drinfeld center twisted by the double dual 
functor. Here the structure map is the left adjoint $\LL^{\mathrm D}$ of the forgetful
functor $U^{\mathrm D} \colon \ZC^{\mathrm D} \rarr{} \C \boti \C$. (Since a modular
category is in particular unimodular, $\LL^{\mathrm D}$ is also a right adjoint of 
$U^{\mathrm D}$.) Moreover, owing to the pivotality of \C\ the
category $\calz^{\mathrm D}(\C)$ is canonically equivalent to \ZC. This leads us to
consider the left adjoint $\LL$ in \eqref{eq:def:LL}.
In the same spirit, embedding two copies of $I$ into $S^1$ gives rise to the functor
$L \cir {\otimes}$, and analogously for any number $n \,{>}\,2$ of copies of $I$.

Also note that it can happen that a boundary circle of a 1-morphism $\surf$ does not
contain the image of any interval or circle. In this case we can consider the sequence
$\emptyset \,{\hookrightarrow}\, I \,{\hookrightarrow}\, S^1$ of embeddings; from
factorization homology the embedding $\emptyset \,{\hookrightarrow}\, I$ gives us
the object $1 \iN \C$, and hence $\emptyset \,{\hookrightarrow}\, S^1$ gives
$L(\one) \iN \ZC$. This is reproduced by setting $n \eq 0$ in
\eqref{eq:L(X1...Xn)} and interpreting the empty tensor product as the monoidal unit.


\subsection{World Sheets} \label{sec:S} %

We now specify the world sheets to which a conformal field theory assigns correlators.
A world sheet is an oriented surface equipped with quite a lot of additional structure.
We therefore proceed in two separate steps: First we describe purely geometric features 
of a world sheet; afterwards we complement these geometric structures with labels
which come from the spherical fusion category \C\ (recall from Section \ref{sec:ocmf}
that \C\ has been selected once and for all). The former step will in particular 
allow us to extract from a world sheet specific 0-, 1- and 2-morphisms of the 
bicategory $\bordoc$, while the latter makes it possible to take the spaces of 
conformal blocks in which the correlators of the CFT are elements as the vector spaces 
that are furnished by a modular functor $\Bl_\C$ of the type discussed in Section
\ref{sec:ocmf}. (Later on, in Section \ref{sec:SNcor} we will define correlators
concretely with the help of the string-net modular functor $\SN$. But for the
considerations from here on until the end of 
                  Chapter  %
\ref{sec:corCFT} we can just
rely on the required properties of $\Bl_\C$ and do not need to invoke the concrete
realization of $\Bl_\C$ as the particular modular functor $\SN$.)

\begin{defi} \label{def:S0}
An \emph{unlabeled world sheet} $\SO$ is a smooth oriented compact surface (possibly
with boundary) equipped with a collection of 
closed
submanifolds of dimension 0, 1 and 2 -- to be referred to as zero-cells, 
one-cells and two-cells, respectively -- satisfying the following conditions:
 \Itemize
 \item
The number of cells is finite.
 \item
The set of zero-cells is the set of boundaries of all one-cells, the union of the 
one-cells is the union of the boundaries of all two-cells, and the union of the 
two-cells is $\SO$.
 \item
The intersection of any pair of one-cells is contained in the set of
zero-cells, and the intersection of any pair of two-cells is contained
in the set of one-cells.
 \item
The interior of every one-cell is either contained in the interior of $\SO$
or contained in the boundary $\partial\SO$.
 \item
Every one-cell in the interior of $\SO$ is oriented as a one-manifold.
 \item
There are two types of \emph{boundary one-cells}, i.e.\ one-cells on $\partial\SO$:
oriented and unoriented. The 1-orientation of an oriented boundary one-cell 
is opposite to the one induced on it\,%
 \footnote{The 2-orientation of an oriented surface induces a 1-orientation on
 its boundary, see e.g.\ \Cite{Ch.\,3}{BOtu}.}
from the 2-orientation of $\SO$.
 \item
A zero-cell in $\partial\SO$ at which two unoriented boundary one-cells
meet is in addition met by precisely one oriented one-cell, whose
        interior is
contained in the interior of $\SO$.
A zero-cell in $\partial\SO$ at which an unoriented and an oriented boundary one-cell
meet, is not met by any further one-cell.
 \item
Every connected component of $\partial\SO$ contains at least one zero-cell.
 \end{itemize}
\end{defi}

\begin{exa} \label{exa:S0}
The following picture shows an example of an unlabeled world sheet of genus 1 
with three geometric boundary circles:
  \be 
  \scalebox{1.6}{\input{WS0.tikz}}
  \label{eq:pic:exampleSO}
  \ee
Here for clarity we have dyed the two two-cells of $\S$ by different colors.
\end{exa}

The following terminology is motivated by the role that world sheets play for
conformal field theory correlators. The one-cells in the interior of an unlabeled world 
sheet $\SO$ are called \emph{defect lines}.
The oriented one-cells in $\partial\SO$ are called \emph{physical boundaries}.
Zero-cells in the interior and zero-cells on the boundary at which two
physical boundaries meet are called \emph{defect junctions}.
The unoriented one-cells in $\partial\SO$ are called \emph{gluing boundaries},
or also \emph{field insertion boundaries}.
We refer to the connected components of the boundary of $\SO$ as
\emph{geometric boundary circles} and denote the set of these by $\pi_0(\partial\SO)$.
If all one-cells in a geometric boundary circle $b$ are gluing boundaries, then $b$
is called a \emph{gluing circle}. If a geometric boundary circle $b$
contains at least one physical boundary, then
each connected component of the complement of the interior of the union of all
physical boundaries in $b$ is called a \emph{gluing interval}.
A geometric boundary circle can contain any finite number of gluing intervals.

Having this terminology at hand, we can now formulate the prescription for adding 
labels to an unlabeled world sheet. These labels involve in particular Frobenius 
algebras in a fixed fusion category 
\C\ and (bi-)modules over them; for information about these structures and for 
pertinent conventions we refer to Appendices \ref{app:A} and \ref{app:M}.

\begin{defi} \label{def:S}
A \emph{world sheet} $\S$ is an unlabeled world sheet $\SO$ together with
the following assignments of labels to the strata of $\SO$:
 \Itemize
 \item
To any two-cell of $\SO$ there is assigned 
a simple special symmetric Frobenius algebra in \C.
 \item
To any defect line there is assigned an $A_{\rm l}$-$A_{\rm r}$-bimodule in \C,
where $A_{\rm l}$ and $A_{\rm r}$ are the labels for the adjacent two-cells 
to the left and right of it, respectively, where the distinction between left and
right is determined by the orientations of the one- and two-cells.
 \item
To any physical boundary there is assigned a right $A$-module, where $A$
is the label for the adjacent two-cell.
 \item
To any defect junction $v$ in the interior of $\SO$ there is assigned a bimodule 
morphism in a space $H_v$ determined by the labels of the defect lines that meet at 
$v$, and to any defect junction $w$ on the boundary of $\SO$ there is assigned a 
module morphism in a space $H_w$ determined by the labels of the defect lines and the
two physical boundaries that meet at $w$.
 \\
(The spaces $H_v$ and $H_w$ will be 
specified in \eqref{eq:def:Hv} and \eqref{eq:def:Hw} below.)
 \item
Gluing boundaries, as well as zero-cells that are not defect junctions (and thus are
located at an end point of a gluing boundary), are not labeled.
 \end{itemize}
\end{defi}

To make this definition complete, we still must specify the morphism spaces in 
which the labels for labeled zero-cells take values. Before doing so, let us first 
illustrate the already somewhat lengthy description of a world sheet by an example. 
{}From here on, unless declared otherwise, by a ``Frobenius algebra'' we mean a
simple special symmetric Frobenius algebra in \C.

\begin{exa} \label{exa:S}
The world sheet
  \be 
  \scalebox{1.6}{\input{WS1.tikz}}
  \label{eq:pic:exampleS}
  \ee
with underlying unlabeled world sheet \eqref{eq:pic:exampleSO} has:
 \Itemize
 \item
two two-cells,
in two different shadings (blue and green in the color version),
which are labeled by Frobenius algebras $A$ and $B$ in \C, respectively;
 \item
six defect lines whose labels are
an $A$-$B$-bimodule $X_6$, $B$-$A$-bimodules $X_1$, $X_2$, $X_3$ and $X_5$,
and a $B$-bimodule $X_4$, respectively;
 \item
two physical boundaries on the geometric boundary circle located to the left,
labeled by an $A$-module $M_2$ and a
       $B$-module
$M_1$, respectively; 
 \item
two gluing circles (the two geometric boundary circles on the right), and one gluing
interval (drawn as a straight line) contained in the geometric boundary circle on
the left;
 \item
three defect junctions with morphism labels $\varphi_1$, $\varphi_2$ and $\varphi_3$
(lying in morphism spaces that will be given in \eqref{eq:Hom4varphi123} below).
 \end{itemize}
\end{exa}

\begin{exa} \label{exa:MXXN0}
Let us also display the situation near a slightly more general gluing interval in 
detail: the following picture shows part of a world sheet near a gluing interval that
connects physical boundaries labeled by $M$ and $N$ and is in addition met by two 
defect lines labeled by bimodules $X_1$ (outgoing) and $X_2$ (incoming):
  \be 
  \scalebox{1.0}{\input{FM1.tikz}} 
  \ee
(The gluing interval, located on the left, 
the physical boundaries and the defect lines are drawn in different shadings
(black, red, and blue, respectively, in the color version).
\end{exa}

Concerning the morphism space $H_v$ assigned to a defect junction $v$ in the interior
of a world sheet, we first note that the orientation of the world sheet furnishes 
a cyclic order of the defect lines that meet at $v$ in such a way that this cyclic 
ordering is clockwise if the orientation of the surface is counter-clockwise.
Suppose that for some selected linear order compatible with this cyclic order, 
the defect lines meeting at $v$ are labeled clockwise by bimodules $X_i$ for 
$i \eq 1,2,...\,,n$. We write $X_i^\vee$ for the bimodule dual to $X_i$, 
as defined at the end of Appendix \ref{app:M}, and set
  \be
  X\pn{} := \left\{ \bearll X & \text{for}~ \epsilon \eq + \,, \Nxl2
  X^\vee & \text{for}~ \epsilon \eq {-} \,. \eear \right.
  \label{eq:def:Xpn}
  \ee
Then $X_i\pn i$ 
is
an $A_i$-$A_{i+1}$-bimodule for $i \,{<}\, n$ and $X_n\pn n$ is an 
$A_n$-$A_1$-bimodule, for Frobenius algebras $A_j$, where
$\epsilon_i \eq +$ if the defect line labeled by $X_i$ is oriented away from $v$,
while $\epsilon_i \eq {-}$ if that defect line is oriented towards $v$.
For the chosen linear order, the space $H_v$ is now defined to be
  \be
  H_v := \Hom_{A_1|A_1}(A_1 \,,
  X_1\pn1 \otx{A_2} X_2\pn2 \otx{A_3}\! \cdots \otx{A_n} X_n\pn n) \,.
  \label{eq:def:Hv}
  \ee
If a different choice of linear order of the defect lines meeting at $v$ is made,
the same prescription gives instead the space
  \be
  \Hom_{A_j|A_j}(A_j \,, X_j\pn j \otx{A_{j+1}} X_{j+1}\pn{j+1} \otx{A_{j+2}}\! \cdots
  \otx{A_{j-1}} X_{j-1}\pn{j-1})
  \ee
for some $j \iN \{2,3,...\,,n\}$ (with labels counted modulo $n$). Owing to pivotality, 
this space is isomorphic to $H_v$ as defined in \eqref{eq:def:Hv}, and indeed 
pivotality provides a distinguished isomorphism of order $n$. Accordingly, the 
choice of linear order is immaterial.

For a defect junction $w$ in the boundary $\partial \SO$, there is directly 
a linear order on the set of physical boundaries and defect lines that meet at 
$w$. If the physical boundary that is oriented towards $w$ is labeled by a right
$A_{n+1}$-module $M$ and the physical boundary oriented away from $w$ is labeled by a 
right $A_1$-module $N$, then the morphism space $H_w$ assigned to $w$ is defined to be
  \be
  H_w := \Hom_{A_{n+1}}(M \,, N \otx{A_1}  X_1\pn1 \otx{A_2} X_2\pn2 \otx{A_3}\! \cdots
  \otx{A_n} X_n\pn n) \,,
  \label{eq:def:Hw}
  \ee
with analogous conventions about the labels of the defect lines as in \eqref{eq:def:Hv}.
In particular, in case $n \eq 0$, so that only the two physical boundaries meet at
$w$, we deal with the space $\Hom_{A_1}(M,N)$.

\begin{exa}
For the world sheet shown in Example \ref{exa:S} the situation near the three
defect junctions looks as follows:
  \be 
  \scalebox{1.0}{\input{HV1.tikz}} \hspace*{3.7em}
  \scalebox{1.0}{\input{HV2.tikz}} \hspace*{-0.2em}
  \scalebox{1.0}{\input{HV3.tikz}} 
  \ee 
Then for a suitable choice of linear order the prescription \eqref{eq:def:Hv} amounts to
  \be
  \varphi_1 \in \Hom_{B|B}(B\,, X_5^{} \otA X_3^\vee \otB X_4^\vee)
  \quad \text{and} \quad
  \varphi_2 \in \Hom_{A|A}(A\,, X_6^{} \otB X_4^{} \otB X_6^\vee) \,,
  \label{eq:Hom4varphi123}
  \ee
while \eqref{eq:def:Hw} gives $\varphi_3 \iN \Hom_B(M_1^{} \,, M_2^{} \otA X_5^\vee)$.
\end{exa}

\begin{rem}
Morita equivalent Frobenius algebras in \C\ describe the same full conformal
field theory \cite{fuRs4}. To perform our string-net construction, for each two-cell
of the world sheet we need to pick a specific Frobenius algebra representing its
Morita class. The final result for the correlators will 
   depend on this choice only very mildly.
Our approach shares this feature with the TFT construction \cite{fuRs4,fuRs10,fjfrs} 
of CFT correlators.
More specifically, as explained in \Cite{Eq.\,(3.37)}{ffrs5},
the correlators calculated with two Morita equivalent
algebras $A$ and $A'$ differ by an over-all factor of
$\big(\sqrt{\dim(A')/\dim(A)}\big)_{}^{\chi(\SO)}$, with $\chi(\SO)$ the Euler character 
of the world sheet. (In  applications of the free boson RCFT to string theory,
this factor combines with a change of the string coupling constant to make
string amplitudes invariant under T-duality; see Section 5.4 of \cite{fGrs}.)
\end{rem}

Now we would like to assign to any world sheet a correlator, and thus, as a first step,
a space of conformal blocks in which the correlator is a specific element. To get the
conformal blocks we associate to a world sheet $\S$ a compact oriented smooth
surface $\surf_\S$ that carries the structure of a 1-morphism of the open-closed 
bordism category $\bordoc$, from which we then obtain the conformal blocks by applying
an open-closed modular functor $\Bl_\C$ that satisfies the requirements formulated in
Section \ref{sec:ocmf}. As an oriented surface, $\surf_\S$ is just the surface that 
underlies the world sheet $\S$, i.e.\ the surface obtained by forgetting all extra
structure of the unlabeled world sheet $\SO$ except for the boundary one-cells. To
make this surface into a 1-morphism we endow it with a \emph{parametrization} --
that is, we regard the union of all gluing circles and gluing intervals of
$\S$ as the image of an object of $\bordoc$ under a diffeomorphism. In particular,
for every gluing interval $r$ we specify a diffeomorphism from $I \iN \bordoc$ to $r$,
and for every gluing circle $b \iN \pi_0(\partial\SO)$ a diffeomorphism from 
$S^1 \iN \bordoc$ to $b$, imposing in addition the condition that $-1 \iN S^1$ does 
not get mapped to any zero-cell of $\S$. 
These diffeomorphisms are required to preserve respectively reverse the orientation, 
depending on whether the intervals or circles are outgoing or incoming boundaries
of $\surf_\S$. Moreover, for a geometric boundary circle containing several gluing 
intervals, the parametrization provides a linear order of those gluing intervals.
{}From now on we will suppress the parametrization in the notation
and just write $\surf_\S$ for the 1-morphism of $\bordoc$ that is obtained
from the world sheet $\S$ in the manner just described. Also, for brevity we will
say that $\surf_\S$ is obtained from a \emph{world sheet with parametrized boundary},
even though the subset of $\partial\SO$ that is given by the physical boundaries is
not parametrized.

\begin{rem}
(i) Note that the parametrization is \emph{not} part of the definition of a world 
sheet. The definition of the conformal block spaces $\Bl_\C(\S)$ and of the string-net 
correlators $\CorSN(\S)$ -- see \eqref{eq:def:Bl(S)} and \eqref{eq:CorSN} below --
makes use of the surface $\surf_\S$, but nevertheless $\Bl_\C(\S)$ and $\CorSN(\S)$
only depend on the world sheet $\S$ but not on the parametrization of $\surf_\S$.
 \\[3pt]
(ii) The parametrization also plays a role when discussing conformal blocks within the
interpretation as vector bundles formed by the spaces of solutions to chiral Ward 
identities that we mentioned in the beginning of Section \ref{sec:ocmf}. In that
context, the parametrization of a gluing circle corresponds to a local holomorphic
coordinate on the disk that is bounded by the circle, and similarly for gluing 
intervals. A coordinate independent definition of conformal blocks can be obtained
in this framework by using the action of the Virasoro algebra.
\end{rem}


\subsection{Field Contents} \label{sec:fields} %

A further ingredient needed for defining the concept of a correlator is a precise notion
of \emph{field insertions}. To achieve this, we now provide a prescription for 
associating an object in \ZC, to be called a field insertion, with each geometric 
boundary circle of a world sheet $\S$ with parametrized boundary: suppressing again 
the dependence on the parametrization in the notation, we define a map
  \be
  \FF \Colon \pi_0(\partial\SO) \rarr{\,} \ZC 
  \label{eq:theFF}
  \ee
from the set $\pi_0(\partial\SO)$ of geometric boundary circles of $\SO$ to the 
Drinfeld center of \C. We refer to $\FF$ as the \emph{field map}.
In the literature, it is customary to associate field insertions, often called
\emph{boundary fields}, also to individual gluing intervals. In our approach we
combine instead such boundary fields for all gluing intervals on a given 
geometric boundary circle $b$ into a single compound field insertion $\FF(b)$ that is 
assigned to the entire circle $b$; for the precise form of this `central lift' of
boundary fields see formula \eqref{eq:def:FF(interval)} below.
This has the advantage that the so obtained compound fields are objects in \ZC,
just like the bulk fields associated to gluing circles are, whereas the boundary 
fields for individual gluing intervals are objects in \C.
Also, when defining the map $\FF$ below, in the concrete formulas we will assume that 
the circles or intervals are \emph{incoming} boundaries of $\surf_\S$; in the case
of outgoing boundaries, the objects given in those formulas must be replaced by their
duals.
 
To define the field map $\FF$ we make use of further categorical tools. We first note 
that for \M\ and \N\ left module categories over \C, to any pair of module functors
$F$ and $F'$ from \M\ to \N\ there is associated an object in the Drinfeld center of \C,
namely the internal Hom
  \be
  \iNat_{\M,\N}(F,F') := \iHom_{\Funre_\C(\M,\N)}(F,F') ~\in \ZC
  \label{eq:def:iNat}
  \ee
for $\Funre_\C(\M,\N)$ as a left module category over the Drinfeld center \ZC. (See 
Appendix \ref{app:Z}
for the description of the category $\Funre_\C(\M,\N)$ of right exact module functors
as a \ZC-mo\-dule category.) The notation $\iNat$ is chosen here in order to indicate
that these objects arise as internal Homs for a functor category. For the same reason 
they are \cite{fuSc25} also called \emph{internal natural transformations}.\,%
 \footnote{~The \ZC-valued internal Homs \eqref{eq:def:iNat} also play a role in the
 context of the so-called full center functor from \C\ to \ZC\ \cite{dakr3}.}
This terminology is further justified by the fact that 
the standard expression of ordinary natural transformations as an end has an
internalized analogue as an end over internal Homs: we have \Cite{Thm.\,9}{fuSc25}
  \be
  \iNat_{\M,\N}(F,F') \,= \int_{\!M\in\M} \iHom_{\N}(F(M),F'(M)) \,.
  \label{eq:iNat=end}
  \ee
When the module categories \M\ and \N\ are obvious from the context, we omit them 
in the notation and just write $\iNat(F,F')$. 

In the considerations below, \M\ and \N\ will often be represented in the form
$\M \eq \Mod A$ and $\N \eq \Mod A'$ for $A$ and $A'$ (simple special symmetric)
Frobenius algebras in \C. The module functors of our interest are then of the form
  \be
  G^Y = -\, \otimes_{\!A} Y \Colon \Mod A \to \Mod A'
  \label{eq:GY=}
  \ee
for some $A$-$A'$-bimodule $Y$ (compare Eq.\ \eqref{eq:def:GB}).


\subsubsection{Fields in the Bulk} \label{sec:bulkfields} %

Internal natural transformations for functors of the form \eqref{eq:GY=} provide us 
with the field insertions for those geometric boundary circles of $\SO$ which are 
gluing circles, i.e.\ which do not contain any physical boundary. For 
giving the precise prescription, consider such a boundary circle
$b \iN \pi_0(\partial\SO)$, and denote by $O_b$ the set of zero-cells contained in $b$
(which are all unlabeled, as none of them lies 
       on
a physical boundary). The orientation 
of $\SO$ induces a cyclic order on $O_b$, such that if the induced orientation of the 
circle is clockwise, then the points are ordered clockwise as well.
The parametrization chosen for the world sheet $\S$ provides a linear order on $O_b$
compatible with that cyclic order.

We can thus proceed
analogously as in the case of the morphism spaces $H_v$ in \eqref{eq:def:Hv}: Suppose 
that for the chosen linear order the defect lines that meet the zero-cells in $O_b$
are labeled clockwise by bimodules $X_i$ for $i \eq 1,2,...\,,n$. Set
$\epsilon_i \eq +$ if the defect line labeled by $X_i$ is oriented away from the
zero-cell at which it meets $b$ and $\epsilon_i \eq {-}$ if it is oriented towards 
the zero-cell, and count labels modulo $n$. Then $X_i\pn i$, with $X\pn{}$ as defined 
in \eqref{eq:def:Xpn}, is an $A_i$-$A_{i+1}$-bi\-mo\-dule for Frobenius algebras $A_j$, 
for $1 \,{\le}\, i \,{\le}\, n$.
Hereby we are ready to define the value of the field map on the gluing circle $b$: we set
  \be
  \FF(b) := \iNat( \Id_{\Mod A_1}  \,, G^{X_n\pn n} {\circ}\, G^{X_{n-1}\pn{n-1}}
  {\circ}\, \cdots \cir G^{X_1\pn 1} ) ~\iN \ZC \,.
  \label{eq:def:FF(circle)}
  \ee

\begin{exa} \label{exa:X.X.X}
To illustrate this prescription, in the following picture we show the world sheet
near a gluing circle with three zero-cells; the position of the image of
$-1 \iN S^1 \,{\subset}\, \complex$ is marked by a
thick dot (drawn in red in the color version);
  \be 
  \scalebox{1.0}{\input{FM0.tikz}} 
  ~~ \xmapsto{~~~} ~~ \iNat(\Id,G^{X_3^{}} \cir G^{X_2^\vee} {\circ}\, G^{X_1^{}}) \,.
  \label{eq:exa:circle3}
  \ee
(Note that the ordering of the labels $X_i$ gets reversed when the functors are
realized explicitly as tensor products over the relevant algebras, as e.g.\ done
in the expression \eqref{eq:def:Hv}.)
\end{exa}

It is worth indicating how the prescription \eqref{eq:def:FF(circle)} relates to more
traditional approaches to fields in conformal field theory (see \cite{fuSc27} and
references given there).
Consider first the special case of a gluing circle $b$ with two zero-cells, to which 
there are attached an incoming and an outgoing defect line, respectively:
  \be 
  \scalebox{1.0}{\input{FM2.tikz}} 
  \label{eq:pic:2.19}
  \ee
In this case we refer to the resulting field insertion $\FF(b)$ as a (two-pronged)
\emph{defect field}. If the two two-cells adjacent to a defect line are labeled by
the same Frobenius algebra and the label of a defect line is that Frobenius algebra,
regarded as a bimodule over itself, then the defect line is called \emph{transparent}. 
A \emph{bulk field} is a special defect field $\FF(b)$ for which the two defect lines 
that meet the gluing circle $b$ are both transparent. If one of the defect lines is
transparent, but the other is not, then $\FF(b)$ is called a \emph{disorder field}.

Now recall that, for \C\ a modular fusion category, any full local conformal field theory
based on chiral data that are encoded in \C\ is characterized by an indecomposable pivotal
left module category \M\ over \C. According to the proposal in \cite{fuSc27}, the types
of defect lines separating two full CFTs based on the same chiral data and on
\C-module categories \M\ and $\M'$, respectively, are given by \C-module functors
from $\M$ to $\M'$, and the defect field for a gluing circle with one incoming defect line 
of type $G_1$ and one outgoing defect line of type $G_2$ -- regarded as
a \emph{field insertion on the defect line} that changes the type of defect from $G_1$
to $G_2$ -- is the object
  \be
  \DD^{G_1,G_2} := \iNat(G_1,G_2) ~\in \ZC
  \label{eq:deffield=iNat}
  \ee
of internal natural transformations. This is indeed in full accordance with the 
prescription \eqref{eq:def:FF(circle)} for the field map $\FF$: there is a canonical
isomorphism
  \be
  \DD^{G_1,G_2} = \iNat(G_1,G_2) \cong \iNat(\Id,G_1\ra \,{\circ}\, G_2^{}) = \FF(b)
  \label{eq:DD=FF}
  \ee
of objects in \ZC, with the gluing circle $b$ as in \eqref{eq:pic:2.19}, namely the one 
obtained by combining the description \eqref{eq:iNat=end} of internal natural 
transformations with the functorial isomorphism \Cite{Lemma\,2}{fuSc25} $\iHom_{\M'}
(G_1(M),G_2(M)) \,{\cong}\, \iHom_{\M'} (\Id,G_1\ra(M) \,{\circ}\, G_2(M))$ for $M\iN\M$.
It is also worth pointing out that taking the Poincar\'e dual of the picture
\eqref{eq:pic:2.19} gives
  \be
  \bearl ~\\[-1.8em]
  \begin{tikzpicture}
  \node[left] at (0,0) {$\M$};
  \node[right] at (1.96,0) {$\M'$};
  \draw[very thick,color=\colorDefect,->] (-0.2,0.23) .. controls (0.7,1.05) and (1.3,1.05) .. (2.2,0.21);
  \node[color=\colorDefect] at (1.02,1.09) {$G_2$};
  \draw[very thick,color=\colorDefect,->] (-0.2,-0.25) .. controls (0.7,-1.05) and (1.3,-1.05) .. (2.2,-0.26);
  \node[color=\colorDefect] at (1.02,-1.14) {$G_1$};
  \draw[thick,double,->] (1,-0.7) --
       node[sloped,xshift=-3pt,yshift=6pt] {$\scriptstyle \DD_{}^{G_1,G_2}$} (1,0.73) ;
  \end{tikzpicture}
  \\[-2.3em]~ \eear
  \label{eq:Poincare}
  \ee
The reminiscence of this picture with the standard graphical description of natural 
transformations fits well with our choice of terminology for $\iNat$.

 \medskip

In the case of a general gluing circle $b$ with any number of zero-cells, in analogy 
with the previous special case we refer to the object $\FF(b)$ a \emph{generalized
defect field}. In the same vein as in \eqref{eq:pic:2.19} we can think of the 
generalized defect field for a gluing circle with $n$ zero-cells as describing
an $n$-pronged defect junction.\,%
 \footnote{~This is in fact the way defect junctions
 are treated in e.g.\ \cite{fScS4,caRun3}.}
If all defect lines that meet a gluing circle $b$ are transparent, then the
field $\FF(b)$ is of the form $\iNat(\Id,\Id)$ and is again called a bulk field.
Since we may trade gluing circles for defect junctions, we refer to all fields $\FF(b)$,
for $b$ any type of gluing circle, also collectively as \emph{fields in the bulk}.

To make contact to the expressions for defect fields, and in particular for bulk fields, 
that are familiar from the conformal field theory literature, we just note that using the
semisimplicity of \C\ it can be shown \Cite{Eq.\,(3.18)}{fuSc27} that for 
$A$-$A'$-bimodules $X_1$ and $X_2$ one has
  \be
  \iNat(G^{X_1}_{},G^{X_2}_{}) \,\cong \bigoplus_{i,j\in\I(\C)}\! \Hom_{A|A'}
  (i \,{\otimes^-} X_1 \,{\otimes^+} j\,, X_2) \otimes_\ko \GC(i \boti j)
  \label{eq:iNat=oplus-i-j}
  \ee
as objects in \ZC.
Here $\I(\C)$ is a set of representatives for the isomorphism classes of simple objects
of \C; we fix the finite set $\I(\C)$ once and for all, in such a way that it contains
the monoidal unit $\one$ of \C. $\GC$ denotes the equivalence \eqref{eq:def:GC}
between $\C\rev \boti \C$ and \ZC, and $i \,{\otimes^-} X_1 \,{\otimes^+} j$
is an $A$-$A'$-bimodule with underlying object $i \oti X_1 \oti j$ and specific
$A$- and $A'$-actions that are defined with the help of the braiding of \C\
(for details see \Cite{p.\,21}{fuSc27}).
In the case of bulk fields for the full conformal field theory that corresponds to a
\C-module category $\M \,{\simeq}\, \Mod A$, \eqref{eq:iNat=oplus-i-j} specializes to
  \be
  \iNat(\Id_{A\MoD A},\Id_{A\MoD A}) \,\cong \bigoplus_{i,j\in\I(\C)}\! \Hom_{A|A}
  (i \,{\otimes^-} A \,{\otimes^+} j\,, A) \otimes_\ko \GC(i \boti j) \,.
  \label{eq:HomAA()}
  \ee
The dimensions of the spaces $\Hom_{A|A'}(i \,{\otimes^-} X_1 \,{\otimes^+} j\,, X_2)$,
respectively $\Hom_{A|A}(i \,{\otimes^-} A \,{\otimes^+} j\,, A)$, are known as
the coefficients of the \emph{torus partition function} for defect fields and bulk
fields, respectively (see e.g.\ Table 1 of \cite{fuRs4}).

 \medskip

When dealing with the objects \eqref{eq:def:FF(circle)} various results about internal 
natural transformations that were obtained in \cite{fuSc25} (in the setting of finite
tensor categories) will be instrumental.


\subsubsection{Fields on Physical Boundaries} \label{sec:bdyfields} %

Next we define the field map on boundary components which do contain a physical 
boundary. Given such a geometric boundary circle $b \iN \pi_0(\partial\SO)$, we proceed 
in two separate steps. First consider an individual gluing interval $r$ that is 
contained in $b$. We associate to $r$ an object $\BB_r \iN \C$ as follows. Denote the 
set of zero-cells at which a defect line meets $r$ by $O_r$. On this set there is a linear
order inherited from the prescription for the orientation of physical boundaries. Assume
that the physical boundary that is oriented towards one of the end points of $r$ is labeled
by a right $A_{n+1}$-module $M$ and the physical boundary oriented away from the other end
point of $r$ is labeled by a right $A_1$-module $N$, and that the defect lines that meet 
the zero-cells in $O_r$ are labeled clockwise by bimodules 
$X_i$ for $i \eq 1,2,...\,,n$. Define $\epsilon_i \iN \{\pm1\}$ and $X_i\pn i$ in the 
same way as in the case of gluing circles, so that $X_i\pn i$ is an 
$A_i$-$A_{i+1}$-bimodule for $1 \,{\le}\, i \,{\le}\, n$. Then we define the object 
$\BB_r$ as the internal Hom
  \be
  \BB_r := \iHom_{\Mod A_{n+1}} \big( M \,, G^{X_n\pn n} {\circ}\, G^{X_{n-1}\pn{n-1}}
  {\circ}\, \cdots \cir G^{X_1\pn 1}(N) \big) ~\in \C \,.
  \label{eq:def:BBr}
  \ee
In particular, in case $n \eq 0$ we have 
  \be
  \BB_r = \iHom_{\Mod A_1}(M,N) =: \BB^{M,N} .
  \label{BB-MN}
  \ee
Note that \Cite{Eqs.\,(2.10),(2.11)}{fuSc27}
  \be
  \iHom_{\Mod A}(M,N) \,\cong \int^{C\in\C}\! \Hom_{A}(C\act M,N) \otimes_\ko C
  \,\cong \bigoplus_{i\in\I(\C)}\! \Hom_{A}(i\act M,N) \otimes_\ko i \,,
  \ee
where the second isomorphism uses semisimplicity of \C. This reproduces the
description of boundary fields familiar from earlier approaches to full conformal
field theory (see again Table 1 of \cite{fuRs4}). Accordingly we call the object
$\BB_r$ in \C\ that is assigned to a gluing interval $r$ 
      a
\emph{boundary field} if $r$ is not met by any defect line, so that $\BB_r$ 
is of the form \eqref{BB-MN}, and a \emph{generalized boundary field} otherwise. 

\begin{exa} \label{exa:MXXN}
For the gluing interval described in Example \ref{exa:MXXN0}, our prescription gives
  \be 
  \scalebox{1.0}{\input{FM1.tikz}} 
  ~~ \xmapsto{~~} ~~ \iHom(M,G^{X_2^\vee} {\circ}\, G^{X_1^{}}(N)) \,.
  \ee
Thus the field $\iHom(M,G^{X_2^\vee} {\circ}\, G^{X_1^{}}(N))$ changes the 
type of boundary condition from $M$ to $N$ in a way such that also defect lines
of type $X_1$ and $X_2$ intervene.
\end{exa}

It should be appreciated that $\BB_r$ as defined by \eqref{eq:def:BBr}
is an object in \C, not in \ZC, unlike the field insertions for gluing circles. Still,
we can attain a uniform treatment of all geometric boundary circles -- both gluing
circles and circles containing gluing
        intervals: 
In a second step we construct from the
objects $\BB_r \iN \C$ for all intervals $r$ on a geometric boundary circle $b$
an object in the Drinfeld center, by treating boundary intervals analogously as we
did for open-closed modular functors $\Bl_\C$ in Section \ref{sec:ocmf}. Thus
we take the tensor product in \C\ of the objects $\BB_r$ for all gluing intervals on $b$
and then apply the left adjoint \eqref{eq:def:LL} of the forgetful functor from \ZC\ 
to \C. Denoting by $Q_b$ the linearly ordered set of gluing intervals in $b$, this gives
  \be
  \FF(b) := \LL\big( \mbox{$\bigotimes_{r\in Q_b}$} \BB_r \big) ~\in \ZC \,,
  \label{eq:def:FF(interval)}
  \ee
where the ordering of factors in the tensor product is according to the linear order
of the boundary intervals $r\iN Q_b$ that results from the parametrization of the
world sheet.
(If the parametrization is changed in such a way that the linear order on $Q_b$ 
is changed, while still keeping the cyclic order that results from the orientation
of the world sheet, then the prescription \eqref{eq:def:FF(interval)} yields an
isomorphic object, since according to Lemma \ref{lem:Lxy=Lyx} 
we have $\LL(C \oti C') \,{\cong}\, \LL(C' \oti C)$ for $C,C' \iN \C$.)
Analogously as for \eqref{eq:L(X1...Xn)}, in the extremal case that the set $Q_b$
is empty, i.e.\ that the boundary circle $b$ does not contain any gluing boundary, 
the expression \eqref{eq:def:FF(interval)} is interpreted as the empty tensor 
product, i.e.\ $\FF(b) \eq \LL(\one) \iN \ZC$.

 \medskip

We summarize our prescriptions in

\begin{defi}
The field map $\FF$ from the set of geometric boundary circles
of $\SO$ to the Drinfeld center of \C\ is given by
  \be
  \FF(b) := \left\{ \bearll
  \eqref{eq:def:FF(circle)} & \text{if $b$ does not contain a physical boundary} \,,
  \Nxl2 \eqref{eq:def:FF(interval)} & \text{if $b$ contains a physical boundary} \,.
  \eear \right.
  \label{eq:def:FF2}
  \ee
\end{defi}

We call the field insertion $\FF(b)$ in the Drinfeld center \ZC\ that according to 
\eqref{eq:def:FF(interval)} is assigned to an entire boundary circle $b$ with physical 
boundaries the \emph{central lift} of the collection of (generalized) boundary fields 
that are associated to the individual gluing intervals on the circle.


\subsection{Correlators as Vectors in Spaces of Conformal Blocks} \label{sec:Bl-Cor} %

As seen in Section \ref{sec:S}, our prescriptions associate to every unlabeled world
sheet $\SO$ endowed with a boundary parametrization a 1-morphism $\surf_\S$ of the 
category $\bordoc$. And as seen in Section \ref{sec:fields}, for every world sheet 
$\S$ with boundary parametrization, with labels given in terms of the input datum of 
a spherical fusion category \C, they associate to every gluing circle of $\S$ an object 
in $\ZC \,{\simeq}\, \Bl_\C(S^1)$ and to every gluing interval of $\S$ an object in 
$\C \,{\simeq}\, \Bl_\C(I)$.
(Also recall that in the formulas that were given for those objects in Section 
\ref{sec:fields}, it is assumed that the relevant boundary circle or interval is 
incoming; for an outgoing boundary one must replace the object by its dual.)
We may then rephrase the statements above as follows:
For any world sheet $\S$ with boundary parametrization we are given a functor
  \be
  \Bl_\C(\surf_\S;-) \Colon \Bl_\C(S^1)^{\boxtimes p} \boti \Bl_\C(I)^{\boxtimes q} 
  \rarr~ \vect_\ko
  \label{eq:BlC(S)}
  \ee
with $p$ and $q$ the number of gluing circles and gluing intervals of $\S$, 
respectively, as well as an object $\BS$ in the domain of this functor which comprises
the (generalized) defect and boundary fields of the world sheet $\S$.
Combining this information we get a finite-dimensional vector space
  \be
  \Bl_\C(\S) := \Bl_\C(\surf_\S;\BS) \,;
  \label{eq:def:Bl(S)}
  \ee
we call this vector space the \emph{space of conformal blocks} for the world sheet $\S$.
By definition of $\Bl_\C$, the functor $\Bl_\C(\surf_\S;-)$, and thus also the 
space $\Bl_\C(\S)$ of conformal blocks, carries an action of the mapping class group 
$\Map(\surf_\S)$.
Note that up to isomorphism the space $\Bl_\C(\S)$ only depends on the topology of the
unlabeled world sheet $\SO$ and on the objects that the field map $\FF$ assigns to the
geometric boundary circles of $\S$. In particular, world sheets which differ in the
configuration of defect lines in their interior but agree on their boundary have the
same conformal block spaces and the same representations of the mapping class group.

\begin{rem}
By suitably defining a category of world sheets, the prescription \eqref{eq:def:Bl(S)}
can be turned into a functor from that category to vector spaces. When restricting
to world sheets for which all field insertions are bulk fields, this gives the
\emph{pinned block functor} considered in \Cite{Sect.\,3.3}{fuSc22}.
\end{rem}

Recall further that the vector spaces $\Bl_\C(\S)$ are isomorphic to morphism spaces 
in \ZC.\,%
 \footnote{Also recall that these isomorphisms, like \eqref{eq:Bl=HomZ} and hence
 \eqref{eq:cong:Bl}, are not canonical, but can be specified uniquely
 by the extra datum of a marking of $\surf_\S$. The conformal block space $\Bl_\C(\S)$
 can thus be understood as a colimit of marked blocks over a category of markings,
 compare e.g.\ \Cite{Sect.\,2.3}{scWo5}.}
In particular, if the world sheet $\S$ is connected, then
according to \eqref{eq:Bl=HomZ} we have isomorphisms
  \be
  \Bl_\C(\S) \cong \HomZ\big( \bigotimes_{b \in \pi_0(\partial\SO)} \FF(b) \,,
  \KK^{\otimes g} \big) \,,
  \label{eq:cong:Bl}
  \ee
where $\FF(b)$ is the generalized defect field, respectively central lift of boundary 
fields, for the geometric boundary circle $b$, as given by \eqref{eq:def:FF2}, while
$\KK \,{=} \int^{X\in\ZC}\!\! X^\vee \oti X \eq \mbox{\Large$\oplus$}_{j\in\I(\ZC)}\,
j^\vee {\otimes} j$ is the distinguished object \eqref{eq:def:KK}.
If $\S$ is the disjoint union of connected world sheets $\S^{(\ell)}$ with
$\ell \iN \{1,2,...,N\}$, then, by monoidality of the modular functor,
$\Bl_\C(\S) \,{\cong}\, \mbox{\Large$\otimes$}_{\ell=1}^N \Bl(\S^{(\ell)})$.

 \medskip

In short, for each world sheet the modular functor provides the conformal block space
$\Bl(\S)$. These spaces come with an action of the mapping class group of the underlying 
surface as well as with sewing maps that relate the conformal block spaces for different
world sheets.
To explain the mapping class group actions and the behavior of conformal blocks 
under sewing, it is convenient to introduce a few further concepts. Recall first
that the data of a world sheet $\S$ do not include any labels for 
those zero-cells on the boundary $\partial\SO$ which are not defect junctions.
There is, however, a natural way to label these zero-cells:
\Itemize
 \item
To an unlabeled zero-cell at which a defect line starts, assign the object $X$ of \C\ 
that underlies the bimodule labeling the defect line, and to one at which a defect line 
ends, assign the dual object $X^\vee$.
 \item
To a zero-cell at which a physical boundary starts or ends, respectively, assign the
object $M$ of \C\ underlying the module label of the physical boundary, respectively
the object $M^\vee$. 
\\
(In the special case that a geometric boundary circle consists of only a single
physical boundary and a single zero-cell at which the physical boundary both starts and 
ends, we take $M$ as a label.)
 \\[-1.4em]~
\end{itemize}
   
\noindent
While these additional labels are entirely redundant, 
it will be convenient to introduce them anyway, as they turn out to facilitate the
formulation of various statements below. For concreteness, we refer to them as the
\emph{natural labels} for the zero-cells on the boundary.

\begin{defi} \label{def:bd}
(i) The \emph{boundary datum} $\bd(b)$ of a connected component $b$ of
the boundary of a world sheet with boundary parametrization
is the set consisting of the labels for the defect junctions on $b$ and 
of the natural labels of the other zero-cells on $b$, ordered according to the 
linear ordering of the zero-cells that is furnished by the boundary parametrization.
 \\[3pt]
(ii) Given a boundary datum $\bd \eq (\bet_1,\bet_2,...\,,\bet_n)$, the associated 
\emph{dual boundary datum} $\bd^\vee$ is the cyclically ordered set 
  \be
  \bd^\vee
  = (\bet_1,\bet_2,...\,,\bet_n)^\vee := (\bet_n^\vee,\bet_{n-1}^\vee,...\,,\bet_1^\vee)
  \ee
that is	obtained from $\bd$ by dualizing all elements as objects and morphisms 
in \C\ and \ZC, respectively, and inverting the order.
\end{defi}

The version of mapping class group that is relevant for the definition of correlators
is now as follows.

\begin{defi} \label{def:Maps}	
The \emph{mapping class group} \Maps\ of a world sheet $\S$ is the
group of homotopy classes of orientation preserving diffeomorphisms from $\SO$ to $\SO$ 
that satisfy the following conditions:
 \Itemize
 \item
Each non-transparent defect line and each defect junction that meets at least one
non-trans\-parent defect line is mapped to itself.
 \item
Each geometric boundary circle is mapped to a geometric boundary circle
with compatible boundary parametrization and with the same boundary datum.
\end{itemize}
\end{defi}

\begin{exa}
In the case of the world sheet \eqref{eq:pic:exampleS}, the boundary data for its
three boundary circles are
  \be
  ( \varphi_3 , M_2^\vee, X_1^{} , M_1^{} )\, , \quad ( X_2^\vee , X_3^{} ) 
  \quad \text{and}\quad ( X_1^\vee , X_2^{} ) \,.
  \ee
Thus in particular any mapping class has to map each of the three boundary circles to itself.
\end{exa}

\begin{rem} \label{rem:Maps-lyub}
(i) Certain boundary circles, subject to the second condition, can be interchanged 
by acting with an element of \Maps; thus \Maps\ can in particular contain braid groups.
Also, \Maps\ contains Dehn twists around gluing circles for bulk field insertions;
these are rigidified, by the choice of a distinguished point (the image of $-1 \iN S^1$
as given by the parametrization) on each such circle.
 \\[3pt]
(ii) Recall that by its construction through a modular functor, the conformal block 
space $\Bl(\S)$ is actually equipped with an action of the entire mapping class group 
$\Map(\surf_\S)$ of the geometric surface $\surf_\S$, of which \Maps\ generically 
is a proper subgroup. In 
                  Chapter  %
\ref{sec:uCorr} we will, however, explain that
besides the group \Maps\ indeed also a larger subgroup of $\Map(\surf_\S)$
plays a role, for which the requirements on the diffeomorphisms are relaxed.
 \\[3pt]
(iii) The isomorphisms \eqref{eq:cong:Bl} are isomorphisms of vector spaces with 
$\Map(\surf_\S)$-action. In 
                  Chapter  %
\ref{sec:the construction}, conformal blocks will
be expressed in terms of string nets. When doing so, the mapping class group action 
is induced by the action of diffeomorphisms on graphs embedded in $\surf_\S$ (see
Section \ref{sec:SN}) and is thus purely geometric. In contrast, on the morphism 
spaces $\HomZ(\mbox{\Large$\otimes$}_{b\,} \FF(b),\KK^{\otimes g})$ furnished by the
isomorphisms \eqref{eq:cong:Bl}, the action is in the first place given algebraically
by composition with braiding and twist isomorphisms and with automorphisms of the
object $\KK^{\otimes g}$ that involve the structural morphisms of $\KK$ as a Hopf 
algebra in \ZC\ \cite{lyub11}. There is, however, also a more geometric skein 
theoretical interpretation of this action \cite{dggpr2}.
\end{rem}

Next we describe how world sheets are sewn.
We say that two distinct gluing circles $b'$ and $b''$ of a world sheet $\S$ 
\emph{match} iff there exists a diffeomorphism from $b'$ to $b''$ that is compatible 
with the boundary parametrization and for which the boundary data satisfy
$\bd(b'') \eq (\bd(b'))^\vee$ as ordered sets 
(which implies in particular that $\FF(b'') \,{=}\, \FF(b')^\vee$).
Matching gluing intervals are defined analogously. The \emph{sewing}, or gluing, of a
world sheet $\S$ is performed either along matching gluing circles or along matching
gluing intervals of $\S$. In the case of circles, the sewn world sheet $\Cup_{b',b''}\S$
is obtained from $\S$ by a diffeomorphism, compatible with the boundary parametrization,
between two matching gluing circles $b'$ and $b''$ of $\S$, whereby every 
zero-cell of $b'$ is mapped to a zero-cell of $b''$ with dual label,
and afterwards omitting the gluing circle.
(More precisely, the one-cells of the gluing circle are deleted, while each zero-cell
together with the two defect lines meeting at it are replaced by a single defect line.
Further aspects, like the fact that,
as described in Remark \ref{rem:collars},
the smooth structure on the glued world sheet $\Cup_{b',b''}\S$ depends on a choice
of collars, are inessential for our purposes.)
Locally the sewing along a circle looks schematically like
  \be
  \bearl ~\\[-1.8em]
  \scalebox{1.0}{\input{SC0.tikz}}
  ~~ \xmapsto{~~~} ~~~
  \scalebox{1.0}{\input{SC1.tikz}}
  \\[-2.1em]~ \eear
  \ee
 ~\\[-0.5em]
(the grey circle in the picture on the right hand side indicates the omitted gluing 
circle, whose intersections with the defect lines are the omitted zero-cells).

The sewing of a world sheet along matching gluing intervals $r' \iN O_{b'}$ and 
$r'' \iN O_{b''}$ proceeds analogously, 
with the additional requirement that the labels of the physical boundaries that meet the
gluing interval must match. The local picture is now schematically
  \be
  \bearl
  \scalebox{1.0}{\input{SI0.tikz}}
  ~ \xmapsto{~~~} ~~~
  \scalebox{1.0}{\input{SI1.tikz}}
  \\[-0.5em]~ \eear
  \ee
 ~\\[-1.6em]
(the grey line in the picture on the right indicates the omitted gluing interval).

 \medskip

Any such sewing gives rise to a gluing transformation between the corresponding 
1-mor\-phisms $\surf_\S$ and $\surf_{{\displaystyle\cup}_{b',b''}\S}$ in $\bordoc$ 
which generalizes the natural transformation \eqref{eq:gluingBl}. As a consequence it
induces a \emph{sewing map}
  \be
  \sew_{b',b''}^{} \Colon \Bl_\C(\S) \to \Bl_\C(\Cup\S_{b',b''})
  \ee
between the respective spaces of conformal blocks. The sewing map can be defined with the 
help of the dinatural structure morphism $\imath$ of a coend that describes the sewing 
procedure, similarly as for the coend \eqref{eq:coendBl}.
This will be explained in detail in Section \ref{sec:SN-mofu}.
For now, let us just mention that upon applying the isomorphisms \eqref{eq:cong:Bl}
we obtain corresponding coends over morphism spaces in the Drinfeld center \ZC, and
present the so obtained coends for the simplest case of sewing along a gluing circle.
The precise form of the coend depends on whether $b'$ and $b''$ lie on the same 
connected component of $\SO$ or not. We present only the former case, whose description 
involves less technicalities. We can then without loss of generality assume that $\SO$ 
is connected. Then $\sew_{b',b''}$ is the morphism that makes the diagram
  \be
  \begin{tikzcd}
  \Bl_\C(\S) \ar{d}[swap]{\sew_{b',b''}} \ar{r}{\cong}
  & \HomZ( \FF(b')^\vee {\otimes}\, \FF(b') \,{\otimes}\!\!\!
  \bigotimes_{\scriptstyle b\in \pi_0(\partial\SO)
  \atop \scriptstyle b \ne b',b''}\!\!\!\! \FF(b), \KK^{\otimes g})
  \ar[start anchor={[yshift=3.9ex]},end anchor={[yshift=-1.2ex]}]{d}[swap]{\imath_{\FF(b')}^{}}
  \\
  \Bl_\C(\Cup_{b',b''}\S) \ar{r}{\cong}
  & \int^{Y \in \ZC}\! \HomZ( Y^\vee {\otimes}\, Y \,{\otimes}\!\!\! 
  \bigotimes_{\scriptstyle b\in \pi_0(\partial\SO) \atop \scriptstyle b \ne b',b''}\!\!\!\!
  \FF(b) , \KK^{\otimes g}) 
  \end{tikzcd}
  \label{coendHomZC-connected}
  \ee
commute, where $\imath_{\FF(b')}^{}$ is the $\FF(b')$-component of the structure
morphism of the coend.
Since the two horizontal isomorphisms in this diagram depend on the extra datum of a
marking of $\surf_\S$ and of $\surf_{{\displaystyle\cup}_{b',b''}\S}$, respectively,
actually some further effort is needed to make \eqref{coendHomZC-connected} into a full
definition of the sewing map $\sew_{b',b''}$. We refrain from giving any details since,
as already stated, below we will give a definition of $\sew_{b',b''}$ directly in terms 
of string nets. (Also, in the generalization to non-semisimple
\C\ the coend has to be taken in the sense of coends of left exact functors, see
Appendix B of \cite{lyub11} and Section 2.3 of \cite{fuSc22} for details.) 

It is worth noting that
when realizing a sewing as a 1-morphisms in $\bordoc$, all outgoing boundary
circles of $\S'$ and all incoming boundary circles of $\S''$ participate in the
sewing. To realize also `partial sewings' as 1-morphisms one has to invoke in
addition dualities to turn incoming to outgoing circles or vice versa. 
This may be taken as an indication that besides the structure of a bicategory,
also other paradigms, such as the structure of a modular operad \cite{joKo2,rayn},
can be convenient for formalizing bordisms.

 \medskip

We are now finally ready to introduce the objects of our primary interest, 
the correlators that are assigned to world sheets:

\begin{defi} \label{def:cosyCor}
A \emph{consistent system of correlators}, for a given spherical fusion category \C,
is a collection of vectors $\Cor(\S) \iN \Bl_\C(\S)$, one for each world sheet $\S$,
that satisfy:
\\[2pt]
(i) Invariance under the mapping class group of $\S$:
  \be
  \gamma(\Cor(\S)) = \Cor(\S)
  \ee
for every mapping class $\gamma \iN \Maps$.
\\[2pt]
(ii) Compatibility with sewing:
  \be 
  \Cor(\Cup_{b,b'}\S) = \sew_{b,b'}^{} (\Cor(\S))
  \ee
for every sewing of a world sheet along gluing circles, and
  \be 
  \Cor(\Cup_{r,r'}\S) = \sew_{r,r'} (\Cor(\S))
  \ee
for every sewing of a world sheet along gluing intervals.
\end{defi}

\begin{rem} \label{rem:Cor=nattrafo}
(i) We may think of a correlator also as the linear map from \ko\ to $\Bl_\C(\S)$ that 
takes the value $\Cor(\S)$ at $1 \iN \ko$.
\\[2pt]
(ii) World sheets can be taken to constitute the objects of a symmetric monoidal category,
with tensor product given by disjoint union and with morphisms given by combinations of 
mapping classes and sewings. The collection of conformal blocks then furnishes a
symmetric monoidal functor from the category of world sheets to the category of
finite-dimensional \ko-vector spaces. In this setting, a consistent system of
correlators is the same as a monoidal natural transformation from a trivial functor that
assigns to every world sheet the ground field $\ko$ to the functor of conformal blocks.
For the case of surfaces without physical boundaries or defect lines, this approach to
correlators is described in detail in \Cite{Sect.\,3.4}{fuSc22}.
\end{rem}

  ~ %

\section{Correlators from String Nets} \label{sec:the construction} %

This chapter contains our central construction, the one of \emph{string-net correlators}, 
as defined in Definition \ref{def:CorSN}. The chief result -- Theorem 
\ref{thm:correlators} -- states that the construction does provide a consistent
system of correlators in the sense of Definition \ref{def:cosyCor}. To prepare 
the ground for this result, Section \ref{sec:SN} introduces the
notion of string nets. These are defined on surfaces with boundaries,
which account for the presence of gluing boundaries (as opposed to 
physical boundaries) of world sheets. The next step, taken in
Section \ref{sec:ipoco}, introduces categories of \emph{boundary values} and string-net
spaces for given boundary values. We can then set up a string-net modular functor
(Section \ref{sec:SN-mofu}), and can explain (Section \ref{sec:SNws}) how a world sheet
determines a string net and thereby a vector in the space of conformal blocks. 
The prescription involves for each two-cell a \emph{\ffg}, a certain type of graph with
edges labeled by a special symmetric Frobenius algebra and vertices labeled by structural
morphisms of the Frobenius algebra. The string-net correlator does not depend on
the particular choice of \ffg; this is crucial for achieving
invariance of the correlators under the mapping class group.

\subsection{String Nets} \label{sec:SN} %

As we will see in Section \ref{sec:SNws}, the data of a world sheet $\S$ provide
us in particular with sufficient information to apply the string-net construction as
formulated in \cite{kirI24} to the surface $\surf_\S$ that is obtained from $\S$,
thereby yielding finite-dimensional vector spaces carrying an action of the mapping
class group $\Maps$.  As a preparation, we here provide some pertinent information about
this construction (for further details see \cite{kirI24,scYa}). We start with the notions
of a coloring by \C\ of an oriented graph and of embedding a graph into a surface. 
By a \emph{surface} $\surf$ we mean a compact oriented smooth surface which may have a 
non-empty boundary. We consider embeddings $\Gama \To \surf$ of 
finite graphs $\Gama$ in a surface $\surf$
for which the intersection of the image of $\Gama$ with the boundary $\partial\surf$
consists of isolated points each of which is the image of a univalent vertex of
$\Gama$. In this case we simply say that the graph $\Gama$ is \emph{embedded in the
surface $\surf$}, and we denote by $\Vd(\Gama,\surf)$ the set of the latter univalent
vertices of $\Gama$, and by $V(\Gama,\surf)$ the set of all other vertices of $\Gama$.

Together with each edge $e$ of an oriented graph $\Gama$ we consider also the same 
underlying unoriented edge endowed with the opposite orientation, which we denote by
$\overline e$.
We write $E(\Gama)$ for the set of edges of $\Gama$ and 
$\EE(\Gama) \,{:=}\, E(\Gama) \,{\cup}\, \{\overline e \mid e \iN E(\Gama)\}$.
For $v$ a vertex of $\Gama$ we denote by $E_v$ the set of edges that are incident to $v$. 
 
\begin{defi} \label{def:coloring}
Let $\Gama$ be a finite graph embedded in a surface $\surf$.
 \\[3pt]
(i) A \emph{coloring} of $\Gama$ by the category \C\ consists of a map
  \be
  \colE \Colon \EE(\Gama) \rarr~ \C
  \ee
satisfying $\colE(\overline e) \eq (\colE(e))^\vee$ for all $e \iN \EE(\Gama)$, and of
an assignment $\colV$ of a morphism
  \be
  \colV(v) \in \HomC(\one, \colE(e_1)\oti\colE(e_2) \,{\otimes}\cdots{\otimes}\, \colE(e_n))
  \label{eq:colVv}
  \ee
to any vertex $v$ of $\Gama$, with
$e_1,e_2,...\,,e_n$ the edges in $E_v$, taken in clockwise order\,%
 \footnote{~Strictly speaking, $\colV(v)$ is not an element of the morphism space
 shown in \eqref{eq:colVv}, but rather in the space that is obtained as a limit
 over all linear orders compatible with the cyclic order, analogously as in
 \eqref{eq:lim-cyclic}. By abuse of notation, we suppress this issue.}
and with orientation away from $v$.
 \\[3pt]
(ii) An \emph{isomorphism} of colorings $(\colE,\colV)$ and $(\colE',\colV')$ of $\Gama$ 
is a collection $\{ f_e \mid e \iN \EE(\Gama) \}$ of isomorphisms
$f_e \colon \colE(e) \Rarr\cong \colE'(e)$ satisfying
$f_e^\vee \cir f_{\overline e}^{} \eq \id_{\colE(\overline e)}$ and
  \be
  \colV'(v) = \big( \bigotimes_{e\in E_v} f_e\big) \circ \colV(v) \,.
  \ee
(iii) The \emph{boundary value} for a coloring of $\Gama$ is the map 
  \be
  \colB^\Gama \Colon \Vd(\Gama,\surf) \to \C
  \label{eq:def:colB}
  \ee
that is given by $\colB^\Gama(v) \,{:=}\, \colE(e_v)$ with $e_v$ the edge incident to $v$ 
in case $e_v$
is oriented towards $v$, and by $\colB^\Gama(v) \,{:=}\, (\colE(e_v))^\vee$ otherwise.
\end{defi}
 
\begin{defi} \label{def:bv-etc}
Let $\surf$ be a surface.
 \\[3pt]
(i) A \emph{boundary value} for $\surf$ is a finite collection $B$ of mutually disjoint
points on $\partial\surf$ together with a map $\bvo\colon B \To \C$ to the objects of \C.
 \\[3pt]
(ii) The set $\Gr(\surf,\bvo)$ is the set of all \C-colored finite graphs $\Gama$ embedded in
         $\surf$ 
with boundary value $\colB^\Gama \eq \bvo$.
 \\[3pt]
(iii) $\kGr(\surf,\bvo)$ is the \ko-vector space freely generated by the set
         $\Gr(\surf,\bvo)$.
\end{defi}

If $\surf \eq \disk$ is a disk with counter-clockwise oriented boundary
and $X_i \iN \C$ for $i \eq 1,2,...\,,n$, denote by $\bvo_X$ the map that associates the
objects $X_i$ in clockwise cyclic order to a collection of $n$ points on $\partial\disk$.
Via the graphical description of morphisms in a spherical fusion category as given
in \eqref{eq:roundcoupon}, a graph $\Gama \iN \Gr(\disk,\bvo_X)$ can then
in an obvious way be regarded as representing a morphism in 
$\HomC(\one,X_1 \,{\otimes}\cdots{\otimes}\,X_n)$, i.e.\ we have a canonical surjection
  \be
  \evD- \Colon \kGr(\disk,\bvo_X) \twoheadrightarrow
  \HomC(\one,X_1 \,{\otimes}\cdots{\otimes}\,X_n) 
  \label{eq:evD}
  \ee
of vector spaces. Thus the morphism space $\HomC(\one,X_1 \,{\otimes}\cdots{\otimes}\,X_n)$
is isomorphic to the quotient
$\kGr(\disk,\bvo_X)/\mathrm{ker}(\evD-)$ in which any two
elements of $\kGr(\disk,\bvo_X)$ that represent the same morphism are identified. 
By mapping graphs via \eqref{eq:evD} to morphism spaces we recover the graphical calculus 
for the monoidal category \C\ as a calculus for graphs on the disk $\disk$.
 
It is natural to extend this procedure to arbitrary surfaces $\surf$, such that we can
use the graphical calculus for \C\ locally for any disk $\disk \,{\subset}\, \surf$
on $\surf$. To this end we must provide a suitable analogue of the quotient 
$\kGr(\disk,\bvo_X)/\mathrm{ker}(\evD-)$ for general surfaces. 

\begin{defi}
Let $\surf$ be a surface and $\bvo$ a boundary value for $\surf$.
 \\[3pt]
(i) Let $\disk \,{\subseteq}\, \surf$ be an embedded disk with 
$\partial\disk \,{\cap}\, \partial\surf \eq \emptyset$. A \emph{null graph with 
respect to $\disk$} is an element $\sum\nolimits_i c_i\,\Gama_{\!i}$ of $\kGr(\surf,\bvo)$
such that, for each $j$, no vertex of $\Gama_{\!j}$ lies on $\partial\disk$ and every edge
of $\Gama_{\!j}$ intersects $\partial\disk$ transversally, that all $\Gama_{\!i}$
coincide on $\surf {\setminus} \disk$, and such that
$\sum\nolimits_i c_i\,\evD{\Gama_{\!i} \,{\cap}\, \disk} \eq 0$. 
 \\[3pt]
(ii) The \emph{null space} $\Null(\surf,\bvo)$ for $\surf$ and $\bvo$ is the subspace
of $\kGr(\surf,\bvo)$ that is spanned by all null graphs for all disks embedded in $\surf$.
 \\[3pt] (iii) The \emph{string-net space} for $\surf$ and $\bvo$ is the quotient
  \be
  \SNo(\surf,\bvo) := \kGr(\surf,\bvo) / \Null(\surf,\bvo)
  \ee
of $\kGr(\surf,\bvo)$.
\end{defi}

We call a vector in the quotient space $\SNo(\surf,\bvo)$ that is the image of an 
element of the generating set $\Gr(\surf,\bvo)$ of $\kGr(\surf,\bvo)$
a \emph{\bare\ string net}, or also just a \emph{string net}. (The reason for the 
qualification ``\bare'' and for the choice of notation $\bvo$ will become clear 
in the next subsection.) A string net that has a graph $\Gama$ as representative is 
denoted by $\GAMA$; by abuse of language, the term string net is also used for individual 
graphs that represent an element $\GAMA \iN \SNo(\surf,\bv)$.
The string-net space is linear in the color of each vertex of a graph and additive with 
respect to taking direct sums of objects labeling the edges. Isotopic graphs and graphs 
with isomorphic colorings give the same string net. Furthermore, all identities valid
in the graphical calculus for \C\ also hold inside any disk embedded in $\surf$. Thus
for string nets the graphical calculus for \C\ applies locally on $\surf$.

Diffeomorphisms of the surface $\surf$ act naturally on embedded graphs. Since in the
string-net space isotopic graphs are identified, this action descends to an
action of the mapping class group of $\surf$ on $\SNo(\surf,\bvo)$.

String nets can be \emph{concatenated}, in a manner similar to the sewing of world sheets. 
Consider a string net $\GAMA \iN \SNo(\surf,\bvo)$ represented by an embedded graph 
$\Gama$. Given the set $B$ of points on $\partial\surf$ that underlies the boundary
datum $\bvo$, denote by $B|_b$ and $B|_{b'}$ the subsets of $B$ consisting of 
those points which lie on two distinct boundary circles $b$ and $b'$ of $\surf$. Assume 
that $B|_b$ and $B|_{b'}$ have the same cardinality, and let $\varphi\colon b \To b'$ be 
a diffeomorphism of one-manifolds that reverses orientation (with respect to the 
orientations on the circles induced by the 2-orientation of $\surf$)
and maps each point in $B|_b$ to a point in $B|_{b'}$. Assume 
further that $\bvo(\varphi(\bet)) \,{=}\, (\bvo(\bet))^\vee$ for each $\bet \iN B|_b$.
In this situation we can sew the surface $\surf$ by identifying the circles $b$ and $b'$
in the same way as is done when sewing world sheets (see Section \ref{sec:Bl-Cor}). 
The sewing identifies in particular the (univalent) vertices of the graph $\Gama$ that 
lie on the boundary circle $b$ with those that lie on $b'$. This results in a finite
embedded graph $\Gama'$ that intersects the identified circles $b$ and $b'$ at points $v$
which are two-valent vertices of $\Gama'$ with an incoming and an outgoing edge carrying
the same label $X(v) \iN \C$. We denote by $\Cup_{b,b'}\Gama$ the graph obtained 
from $\Gama'$ by replacing each such two-valent vertex $v$ and the two edges that meet
at $v$ by a single edge with label $X(v)$. The \emph{concatenated string net} 
$[\Cup_{b,b'}\Gama]$ is the string net on the sewn surface $\Cup_{b,b'}\surf$
that is represented by the graph $\Cup_{b,b'}\Gama$.
Note that by construction we have $[\Cup_{b,b'}\Gama] \,{\cong}\, [\Cup_{b',b}\Gama]$.
Similarly, if $r,r' \,{\subset}\, \partial\surf$ are intervals such that $r\,{\cap}\,r' \eq
\emptyset$ and $\Gama\,{\cap}\,\partial r \eq \emptyset \eq \Gama\,{\cap}\, \partial r'$,
and if $\bvo(\varphi(\bet)) \,{=}\, (\bvo(\bet))^\vee$ for each $\bet \iN B|_r$ for 
an orientation reversing diffeomorphism $\varphi\colon r \To r'$, then 
we can concatenate the string net along the intervals, resulting in a string net
$[\Cup_{r,r'}\Gama]$ on the sewn surface $\Cup_{r,r'}\surf$.

 %
 %
 %
 %
 %


\subsection{Idempotent Completion} \label{sec:ipoco} %

A further ingredient needed for a string-net construction based on world sheets is 
a suitable abelian category constructed from boundary values. This is obtained by 
first considering a \bare\ variant, which is not abelian, but for which taking its 
Karoubian envelope, i.e.\ completing it such that any idempotent has an image,
yields an abelian category (using that the category \C\ of input data is semisimple).
Recall that an object of the Karoubian envelope
$\mathrm{Kar}(\cald)$ of a category $\cald$ is a pair $(D,p)$ consisting of an object
$D\iN\cald$ and an idempotent $p \iN \Hom_\cald(D,D)$, such that $\id_{(D,p)} \eq p$.
Following \Cite{Sect.\,6}{kirI24} we give

\begin{defi} \label{def:Cyl}
Let $S$ be an oriented one-dimensional manifold, possibly with non-empty boundary.
 \\[3pt]
(i) Let $B$ be a finite collection of points on $S$. The \emph{category} $\Cylo S$ 
\emph{of boundary values} for $S$ is the following category: An object $\bvo$ of $\Cylo S$ 
is a finite collection $B$ of points on $S$ together with a map $\bvo\colon B \To \C$.
The space of morphisms between two objects is the string-net space 
on a cylinder $\cylS$ over $S$ with boundary
value induced by the two inclusions of $S$ as the boundary of the cylinder,
  \be
  \HomCyloS(\bvo,\bvo{}') := \SNo(\cylS,\bvo{}^\vee \,{\cup}\, \bvo{}') \,.
  \ee
Composition of morphisms is given by concatenating string nets.
 \\[3pt]
(ii) The \emph{cylinder category} for $S$ associated with \C\ is the Karoubian envelope
  \be
  \Cyl S := \mathrm{Kar}(\Cylo S)
  \ee
of the category of boundary values for $S$.
\end{defi}

For any two one-manifolds $S$ and $S'$ we have
  \be
  \Cyl{S{\sqcup}S'} = \Cyl S \times \Cyl{S'} \,;
  \ee
since all functors with cylinder categories as domains that we consider in this
paper are left exact in each factor, we may, and will, replace
$\Cyl S \Times \Cyl{S'}$ with $\Cyl S \boti \Cyl{S'}$. As a consequence, the cylinder 
categories for the interval and for the circle are the most interesting ones.
Consider thus first the case that $S \eq I \eq [0,1]$ is a single interval; there is 
then a canonical equivalence
  \be
  \Cylo I \rarr\simeq \C
  \label{eq:CyloI2C}
  \ee
which sends the object in $\Cylo I$ that is given by a collection $\{\xi_j\}$ of points
on $I$ satisfying $0 \,{<}\, \xi_1 \,{<}\, \xi_2 \,{<} \cdots {<}\, \xi_n \,{<}\, 1$ 
and labeled by $\bvo(\xi_j) \eq X_j \iN \C$ to the tensor product
$X_1 \oti X_2 \,{\otimes}\cdots{\otimes}\, X_n$. Moreover, since \C\ is idempotent 
complete, after choosing a splitting for each idempotent in \C\ there is a unique
extension of the functor \eqref{eq:CyloI2C} to an equivalence
  \be
  \Cyl I = \mathrm{Kar}(\Cylo I) \rarr\simeq \C \,.
  \ee
Similarly, in the case that $S \eq S^1$ is a single circle with distinguished point 
$-1 \iN S^1$ (and analogously for any pointed connected compact one-manifold), there
is a canonical functor
  \be
  \Phi^\circ \Colon\Cylo{S^1} \rarr~ \ZC \,.    %
  \label{eq:CyloSe2C}
  \ee
The functor \eqref{eq:CyloSe2C} is fully faithful. It sends the object in $\Cylo{S^1}$ 
that is given by a collection $\{\xi_j\}$ of points on $S^1$ satisfying
$ \pi \,{>}\, \mathrm{arg}(\xi_1) \,{>}\, \mathrm{arg}(\xi_2) \,{>}\cdots {>}\,
\mathrm{arg}(\xi_n) \,{>}\, {-}\pi$ and labeled by $\bvo(\xi_j) \eq X_j \iN \C$ to 
$\LL(X_1 \oti X_2 \,{\otimes}\cdots{\otimes}\, X_n)$.
To describe the action of $\Phi^\circ$ on morphisms, consider
$\varphi \iN \HomC(X,Z^\vee {\otimes}\, Y \oti Z)$ and set
  \be
  \widetilde\varphi ~:=~~
  \scalebox{0.9}{\input{IC0.tikz}} 
  ~~ \in \HomCyloe(S_X,S_Y) \,.
  \ee
The morphism $\Phi^\circ(\widetilde\varphi)\colon \LL(X) \To \LL(Y)$ in \ZC\ is
then defined by the dinatural family
  \be
  U^\vee {\otimes}\, X \oti U \rarr{\id_{U^\vee_{}}\otimes\varphi\otimes\id_U}
  U^\vee {\otimes}\, Z^\vee {\otimes}\, Y \oti Z\oti U
  \rarr{\imath^\Zm_{Y;Z\otimes U}} \LL(Y) \,,
  \ee
where $\imath^\Zm_{Y;-}$ is the dinatural structure morphism of the coend $\LL(Y)$
(compare \eqref{eq:imathZ}).
Using that \C\ is semisimple, and adopting the summation convention \eqref{eq_sum-alpha},
we have
  \be
  \Phi^\circ(\widetilde\varphi) ~=~ \sum_{i,j\in\I(\C)} d_j ~~
  \scalebox{1.1}{\input{IC1.tikz}}
  \label{eq:Phi0(tildephi)}
  \ee

Since \ZC\ is idempotent complete, after choosing a splitting for 
each idempotent in \ZC\ there is a unique extension of $\Phi^\circ$ to a functor
  \be
  \Phi \Colon \Cyle = \mathrm{Kar}(\Cylo{S^1}) \rarr\simeq \ZC
  \label{Cyl=Z}
  \ee
which sends $(\bvo,p) \iN \Cylo{S^1}$ to $\Im(\Phi^\circ(p))$ with
$\Phi^\circ(p) \iN \HomZ(\Phi^\circ(\bvo),\Phi^\circ(\bvo))$. The so obtained functor 
$\Phi$ is an equivalence \Cite{Thm.\,6.4}{kirI24}.
Below we will largely work directly with the cylinder category $\Cyle$ rather than 
with \ZC, so that the specific form of the equivalence \eqref{Cyl=Z} will actually
not be important.

\begin{exa} \label{exa:pcan}
For $Y \eq (U(Y),\cB) \iN \ZC$ consider the string net
  \be 
  \pcan_Y ~~:=~~ 
  \scalebox{0.9}{\input{PC0.tikz}}
  \label{eq:def:PY}
  \ee
in $\HomCyloe(S_{U(Y)},S_{U(Y)})$, where the unlabeled edge, which runs along the 
non-contractible cycle of the cylinder, stands for the \emph{canonical color} 
\eqref{eq:cancolor} and the half-braiding $\cB$ is indicated by an over-crossing, 
as explained in \eqref{pic:cB}. 
           (Also, here, as well as in several pictures below,
           we abuse notation by abbreviating $U(Y)$ to $Y$.)
Using the identity \eqref{eq:dominance} it is readily 
seen that $\pcan_Y$ is an idempotent (see e.g.\ \Cite{Rem.\,2.6}{scYa}). We have
  \be
  \Phi^\circ(\pcan_Y) ~=~ \sum_{i,j,k\in\I(\C)} \frac{d_j d_k}{\DC^2} ~~
  \scalebox{1.1}{\input{IC2.tikz}}
  ~~~=~ \sum_{i,j\in\I(\C)} \frac{d_j}{\DC^2} ~~
  \scalebox{1.1}{\input{IC3.tikz}}
  \ee
with $\DC^2$ the global dimension \eqref{eq:globaldim} of \C.
It is then easily seen \Cite{Lemma\,8.3}{kirI24} that the image of the morphism
$\Phi^\circ(\pcan_Y) \iN 
	   \HomZ(\LL U(Y),\LL U(Y))
$ is the object $Y \iN \ZC$. We write
  \be
  \bvcan_Y := (S_{U(Y)},\pcan_Y)
  \label{eq:def:bvcan}
  \ee
for the object of $\Cyle$ that is given by the boundary value $S_{U(Y)} \iN \Cyloe$
and the idempotent \eqref{eq:def:PY}. We thus have 
  \be
  \Phi(\bvcan_Y) = \Im(\Phi^\circ(\pcan_Y)) = Y ~\in \ZC
  \ee
for every $Y \iN \ZC$.
\end{exa}

\begin{exa} \label{exa:pXY}
For any two Frobenius algebras $A$ and $B$ in \C\
and any two $A$-$B$-bimodules $X$ and $Y$ consider the string net
  \be
  \pz_{X,Y}^{} ~~:=~~
  \scalebox{0.9}{\input{PD0.tikz}}
  \label{eq:def:pXY}
  \ee
in $\HomCyloe(S_{Y,X^\vee},S_{Y,X^\vee})$.
Using that $A$ and $B$ are special Frobenius algebras, one sees again
directly that $\pz_{X,Y}^{}$ is an idempotent. Furthermore we find
  \be
  \Im(\Phi^\circ(\pz_{X,Y}^{}))
  \,\cong\, \iNat(G^X,G^Y) \stackrel{\eqref{eq:deffield=iNat}}= \DD^{X,Y} \,\iN \ZC \,.
  \label{eq:ImPhipXY=DXY}
  \ee
Let us explain this isomorphism in detail, using manifestly that \C\ is semisimple. (The
statement is, however, in fact valid beyond semisimplicity.) Inserting the explicit
form \eqref{eq:GY=} of the functors $G^X$ and $G^Y$, in terms of the
bimodules $X$ and $Y$ we have
  \be
  \DD^{X,Y} = \bigoplus_{m\in\I(\Mod A)} m \otx A Y \otx B X^\vee \otx A m^\vee .
  \ee
Now consider the two string nets
  \be 
  \hspace*{-0.6em}
  \begin{array}{l}
  e_{X,Y}^{} ~:=~ \sum_{m\in\I(\Mod A)} \frac{d_m}{\DC^2} ~~
  \scalebox{1.0}{\input{PD1.tikz}}
  \\~\\[-4.9em]
  \hspace*{12.5em} \text{and} \qquad
  r_{X,Y}^{} ~:=~ \sum_{m\in\I(\Mod A)} ~~
  \scalebox{1.0}{\input{PD2.tikz}}
  \label{eq:eXYrXY}
  \eear
  \ee
where $d_m \eq d_{\dot m}$ is the dimension of the object in \C\ that underlies
$m\iN\Mod A$. Using again that $A$ and $B$ are special Frobenius and invoking 
Corollary \ref{cor:shrink} we see that 
  \be
  e_{X,Y}^{} \circ r_{X,Y}^{} = \pz_{X,Y}^{} \,.
  \label{eq:eXYrXY=pXY}
  \ee
For the composition in the opposite order we get
  \be
  r_{X,Y}^{} \circ e_{X,Y}^{} = \pcan_{\DD^{X,Y}}
  \ee
with $\pcan_Z$ as in \eqref{eq:def:PY}; this is shown by the following identities
between 
string nets:
  \be
  \bearl
  r_{X,Y}^{} \circ e_{X,Y}^{} ~= \sum_{m,n\in\I(\Mod A)} \frac{\dd_m}{\DC^4} ~~ 
  \scalebox{1.1}{\input{PD3.tikz}}
  \\[-2.8em] \hspace*{15em}
  = \!\sum_{\scriptstyle i\in\I(\C) \atop \scriptstyle m,n\in\I(\Mod A)}
  \!\!\!\!\frac{\dd_{i}\,\dd_{m}}{\DC^2}~~
  \scalebox{1.1}{\input{PD4.tikz}}
  ~~=~ \pcan_{\DD^{X,Y}} .
  \\[-2.8em] ~
  \eear
  \label{eq:2halfcirc-cancol}
  \ee
~\\
(Here the cylinder is drawn in a slightly different manner -- deformed by a
diffeomorphism -- than in \eqref{eq:eXYrXY}. The $\alpha$-summation is over 
module morphisms, and the middle equality uses the identity \eqref{eq:proofshrink2}.)
We write
  \be
  \bvz_{X,Y} := (S_{X,Y},\pz_{X,Y})
  \label{eq:def:bvz}
  \ee
for the object of the category $\Cyle$ that is given by the boundary value consisting of
two points labeled by $U(X^\vee)$ and by $U(Y)$, respectively, and by the idempotent
\eqref{eq:def:pXY}.
We have $\Phi(\bvz_{X,Y}) \eq \iNat(G^X,G^Y) \eq \DD^{X,Y} \iN \ZC$. Also,
in the Karoubian envelope the morphisms $\pz_{X,Y}^{}$ and $\pcan_{\DD^{X,Y}}$ are
nothing but the identity morphisms on the objects $\bvz_{X,Y}$ and $\bvcan_{\DD^{X,Y}}$,
respectively. It follows that the string nets \eqref{eq:eXYrXY} provide inverse isomorphisms
  \be
  \begin{tikzcd}
  e_{X,Y}^{} : ~~ \bvcan_{\DD^{X,Y}}\! \ar[yshift=2pt]{r}{\cong}
  & \bvz_{X,Y} ~ :\, r_{X,Y}^{} \ar[yshift=-4pt]{l}
  \end{tikzcd}
  \label{eq:bvcan=bvz}
  \ee
for any pair of $A$-$B$-bimodules $X$ and $Y$
over (simple special symmetric) Frobenius algebras $A$ and $B$.
Since according to Example \ref{exa:pcan} we have $\Phi(\bvcan_{\DD^{X,Y}}) \eq \DD^{X,Y} 
{\in}\, \ZC$, we then get
  \be
  \Im(\Phi^\circ(\pz_{X,Y}^{}))
  \,\cong\, \Phi(\bvz_{X,Y}) \cong \Phi(\bvcan_{\DD^{X,Y}}) = \DD^{X,Y} .
  \ee
\end{exa}

\begin{defi}
Let $\surf$ be a surface, $\bvo$ be a boundary value for $\partial\surf$, and
$p \iN \HomCylod(X,X)$ be an idempotent. The map
  \be
  p_*^{} \Colon \SNo(\surf,\bvo) \rarr{} \SNo(\surf,\bvo)
  \label{eq:def:p*}
  \ee
is the linear endomorphism of the string-net space $\SNo(\surf,\bvo)$ that maps
$\GAMA \iN \SNo(\surf,\bvo)$ to the string net that is obtained by concatenating
$\GAMA$ with the string net $p$ close to the boundary.
\end{defi}

The following observation is immediate:

\begin{lem}
For any idempotent $p \iN \HomCylod(X,X)$ the map $p_*^{}$ defined by 
\eqref{eq:def:p*} is an idempotent.
\end{lem}

\medskip

With this result at hand, we can adapt the definition of string-net spaces to
the Karoubification of the boundary category.

\begin{defi}
Let $\surf$ be a surface and let $\bv \eq (\bvo,p_{\bv})$ be an object in 
the cylinder category $\Cyl{\partial\surf}$.
The \emph{string-net space} for $\surf$ and $\bv$ is the subspace
  \be
  \SN(\surf,\bv) := \Im\big( (p_{\bv})_*^{} \big)
  \ee
of $\SNo(\surf,\bvo)$.
\end{defi}

We refer to an element of $\SN(\surf,\bv)$ again as a string net, or also, if we
want to stress that it is not a \bare\ string net, as a \emph{Karoubified string net}. 
The subspace $\SN(\surf,\bv)$ is invariant under the action of the mapping
class group of $\surf$ on $\SNo(\surf,\bvo)$, and thus carries an action of that group
as well.

\begin{exa}
The picture
  \be
  \bearl ~\\[-1.6em]
  \scalebox{1.0}{\input{GKS0.tikz}}
  \eear
  \ee
shows a (Karoubified) string net on a torus $\mathrm T$ with three boundary circles 
which is an element of the space $\SN(\mathrm T\,, \bvcan_{Y_1} \,{\boxtimes}\,
\bvcan_{Y_2} \,{\boxtimes}\, \bvcan_{Y_3})$.
\end{exa}


\subsection{The String-net Modular Functor} \label{sec:SN-mofu} %

Next we note that as a consequence of the equivalence \eqref{Cyl=Z} we have
a non-canonical equivalence
  \be
  \Cyl{\partial\surf} \simeq {\ZC}^{\boxtimes\ell} \qquad\text{if}\qquad
  \partial\surf \,{\cong}\, (S^1)^{\sqcup\ell} .
  \ee

\begin{prop} \label{prop:SN=HomZ}
Let $\surf$ be a connected surface with boundary $\partial\surf$ diffeomorphic to 
$(S^1)^{\sqcup\ell}$. Let $\Phi_{\partial\surf} \colon \Cyl{\partial\surf} \,{\to}\, 
{\ZC}^{\boxtimes\ell}$ be an equivalence, 
and let $\bv \eq (\bvo,p_{\bv}) \iN \Cyl{\partial\surf}$ be such that
$\Phi_{\partial\surf}(\bv) \,{\cong}\, \bigboxtimes_{i} Y_i \iN {\ZC}^{\boxtimes\ell}$. 
Then there is an isomorphism
  \be
  \SN(\surf,\bv) \,\cong\, \HomZ(\one,\Phi_{\partial\surf}^\otimes(\bv) \oti K^{\otimes g})
  \label{eq:SN=HomZ}
  \ee
of vector spaces, where $g$ is the genus of $\surf$ and
$\Phi_{\partial\surf}^\otimes(\bv) \eq \mbox{\Large{$\otimes$}}_i Y_i \iN \ZC$.
Moreover, the isomorphism \eqref{eq:SN=HomZ} intertwines the action of the mapping
class group of $\surf$ on these vector spaces.
\end{prop}

\begin{proof}
As has been shown in \cite{kirI24,goos}, the (Karoubified) string-net construction based 
on a spherical fusion category \C\ gives rise to a 3-2-1 topological field theory (TFT)
that is isomorphic to the (3-2-1-extended) Turaev-Viro state sum TFT based on \C.
The latter Turaev-Viro TFT is, in turn, isomorphic as a 3-2-1 TFT
to the Reshetikhin-Turaev surgery TFT based on the Drinfeld center \ZC\
\cite{balKi,TUvi}. These isomorphisms of TFTs restrict to isomorphisms between the
respective induced modular functors. 
\end{proof}
	
\begin{rem}
By the conventions entering the notion of boundary value (see Definition \ref{def:bv-etc}),
the unparametrized string-net functor $\SN(\surf,-)$ is naturally a covariant
functor with domain $\Cyl{\partial\surf}$. That field insertions appear
in the contravariant entry of the Hom functor (see e.g.\ \eqref{eq:cong:Bl}) is due to
our choice of convention for the field map $\FF$.
\end{rem}

It is worth recalling that the mapping class group action on the string-net space
on the left hand side of \eqref{eq:SN=HomZ} is directly induced by the action of 
diffeomorphisms on graphs embedded in $\surf_\S$ and is thus geometric, while
the action on the morphism spaces on the right hand side of \eqref{eq:SN=HomZ}
is more intricate (see Remark \ref{rem:Maps-lyub}(iii)).

The result of Proposition \ref{prop:SN=HomZ} extends in an obvious manner to 
non-connected surfaces: If $\surf$ is the 
disjoint union of connected surfaces $\surf^{(q)}$ with $q \iN \{1,2,...,N\}$, then 
  \be
  \SN(\surf,\bv) = \bigotimes_{q=1}^N \SN(\surf^{(q)},\bv|_{\surf^{(q)}}) \,,
  \ee
and for each $\SN(\surf^{(q)},\bv|_{\surf^{(q)}})$ there is an isomorphism of the form 
\eqref{eq:SN=HomZ}.

It follows in particular that for any surface $\surf$, string nets furnish an exact functor
$\SN(\surf,-) \colon 
      $\linebreak[0]$
\Cyl{\partial\surf} \Rarr{} \vect$.
It is worth pointing out that the string-net spaces considered so far are constructed
for surfaces without a boundary parametrization, while the surfaces
obtained from world sheets are equipped with this extra datum. To account for the latter,
we observe that the assignment
  \be
  b \longmapsto \Cyl b
  \ee
extends to a symmetric monoidal 2-functor from the bicategory of one-manifolds, with
orientation preserving embeddings as 1-morphisms and isotopies as 2-morphisms, to the
2-category of categories. Moreover, denoting by $\overline b$ the one-manifold $b$ 
with opposite orientation, there is a canonical equivalence
  \be
  \Cyl{\overline b} \rarr\simeq \Cyl b\opp
  \ee
which maps a boundary value $\bd$ to its dual $\bd^\vee$, i.e.\ to the boundary value
consisting of the same points on the one-manifold $b$ as $\bd$, but each labeled with the
dual of the original object of \C. 
This implies that from any boundary parametrization
  \be
  \vphi \Colon 
  (S^1)^{\sqcup(p+q)}_{} \sqcup I^{\sqcup(r+s)}_{} \rarr{~} \partial\surf
  \ee
of a surface $\surf$ with $p$ incoming and $q$ outgoing boundary circles and with $r$
incoming and $s$ outgoing boundary intervals we naturally obtain a functor
  \be
  \vphi_*^{} \Colon \Cyl{S^1}^{\boxtimes(p+q)} \boxtimes \Cyl{I}^{\boxtimes(r+s)}
  \rarr{~\,} \Cyl{\partial\surf} \,.
  \label{eq:functorPhi*}
  \ee

\begin{rem}
To treat all surfaces that can be obtained as $\surf \eq \surf_\S$ for any arbitrary
world sheet $\S$, we must also account for the case that a geometric boundary circle of
$\surf$ is not pa\-ra\-metrized at all, which happens if the corresponding boundary circle 
of $\S$ entirely consists of a single physical boundary (together with a single 
zero-cell). Accordingly, if $\S$ has $\ell$ such 
boundary circles, then one deals with a functor $\vphi_*^{}$ with domain
$\Cyl{S^1}^{\boxtimes(p+q)} \boti \Cyl{I}^{\boxtimes(r+s)} \boti \vect^{\boxtimes\ell}$
instead of \eqref{eq:functorPhi*}. (Here we use that $\Cyl\emptyset \eq \vect$).
Since $\vect$ acts as a monoidal unit under the Deligne product, we can safely
suppress this complication.
\end{rem}

We can now define a functor $\SN$ for surfaces with boundary parametrization:

\begin{defi}
Let $\surf$ be a surface and $\vphi$ a boundary parametrization of $\surf$. The
functor $\SN(\surf_\vphi,-)$ is the composite
  \be
  \SN(\surf_\vphi,-) \,:=\, \SN(\surf,-) \circ \vphi_*^{} \Colon
  \Cyl{S^1}^{\boxtimes(p+q)} \boxtimes \Cyl{I}^{\boxtimes(r+s)} \rarr{~\,} \vect \,.
  \ee
When the boundary parametrization is clear from the context, we often abuse notation
and simply write $\SN(\surf,-)$ in place of $\SN(\surf_\vphi,-)$.
\end{defi}

We are now ready to state

\begin{thm} \label{thm:sewing}
Let $\surf$ be a $($compact oriented smooth$)$ surface with a boundary parametrization 
$\vphi(\surf) \,{\eq}\, \{\vphi_-(\surf),\vphi_+(\surf)\}$ of the form
  \be
  \vphi_-(\surf) \colon~ \alpha \,{\sqcup}\, \beta \rarr~ \partial\surf
  \qquad \text{and} \qquad 
  \vphi_+(\surf) \colon~ \beta \,{\sqcup}\, \alpha' \rarr~ \partial\surf \,,
  \label{eq:alpha-beta-alpha'}
  \ee
where $\alpha$, $\alpha'$ and $\beta$ are finite disjoint unions of the circle $S^1$ and
the interval $I$.
Denote by $\Cup_\beta\surf$ the surface obtained by gluing $\surf$ along the images
of $\beta$ under $\vphi_\pm(\surf)$ via the orientation reversing diffeomorphism 
$\vphi_+(\surf)|_\beta \cir ( \vphi_-(\surf)|_\beta )^{-1}$.
Then the corresponding gluing transformation 
  \be
  \SN(\surf,\overline ? \boti \bv \boti \overline \bv \boti \mq)
  \rarr~ \SN(\Cup_\beta\surf, \overline ? \boti \mq)
  \label{eq:coendS-0}
  \ee
exhibits the functor $\SN(\Cup_\beta\surf, -)$ as the coend 
  \be
  \SN(\Cup_\beta\surf, \overline ? \boti \mq)
  \,= \int^ {\bv'\in\Cyl{\beta}}\!
  \SN(\surf,\overline ? \boti \bv' \boti \overline{\bv'} \boti \mq)
  \label{eq:coendS}
  \ee
   $($taken in the category of left exact functors$)$, i.e.\ \eqref{eq:coendS} has
\eqref{eq:coendS-0} as its structure morphisms.
\end{thm}

\begin{proof}
Fix for the moment boundary values $\bv_\alpha \iN \Cyl\alpha$ and 
$\bv_{\alpha'} \iN \Cyl{\alpha'}$. We are going to show that the family
of gluing morphisms
  \be
  \imath_{\bvo}^{} \Colon \SN(\surf,\overline{\bv_\alpha} \boti \bvo
  \boti \overline{\bvo} \boti \bv_{\alpha'})
  \rarr~ \SN(\Cup_\beta\surf, \overline {\bv_\alpha} \boti \bv_{\alpha'}) \,,
  \label{eq:iBo}
  \ee
for boundary values $\bvo {\in}\, \Cylo\beta$, is dinatural and 
satisfies the universal property of the structure morphism of the coend 
\eqref{eq:coendS}. Note that here we restrict the family to the subcategory 
$\Cylo\beta \,{\hookrightarrow}\, \mathrm{Kar}(\Cylo\beta) \eq \Cyl\beta$; this
is sufficient, owing to Lemmas \ref{lem:5.1.7} and \ref{lem:Kar-pgen}.
Dinaturality of the family \eqref{eq:iBo} is immediate: The two expressions that 
are required to be equal can both be described as an additional gluing with a 
cylinder $\surf'\colon \beta \To \beta$, namely either gluing the images of $\beta$ 
under $\vphi_+(\surf)$ and $\vphi_-(\surf')$ via $\vphi_+(\surf)|_\beta \cir 
( \vphi_-(\surf')|_\beta )^{-1}$ first, or gluing the images under $\vphi_+(\surf')$ 
and $\vphi_-(\surf)$ via $\vphi_+(\surf')|_\beta \cir ( \vphi_-(\surf)|_\beta )^{-1}$
first; dinaturality then directly follows from the associativity of composition.
 \\[2pt]
To address universality of the family $\{\imath_{\bvo}^{}\}$ we construct, given 
an arbitrary dinatural family
  \be
  \{\, g^{}_{\bvo} \colon \SN(\surf,\overline{\bv_\alpha} \boti \bvo
  \boti \overline{\bvo} \boti \bv_{\alpha'}) \Rarr{} V \,\}_{\bvo \in \Cylo\beta}^{}
  \label{eq:fam:..toV}
  \ee
to some vector space $V$, a linear map
  \be
  g \Colon
  \SN(\Cup_\beta\surf, \overline{\bv_\alpha} \boti \bv_{\alpha'}) \rarr{} V
  \label{eq:def:g..toV}
  \ee
as follows: Any string net
$\GAMA \iN \SN(\Cup_\beta\surf, \overline{\bv_\alpha} \boti \bv_{\alpha'})$
contains a representative graph $\Gama$ that intersects the one-manifold
$\beta \,{\hookrightarrow}\, \Cup_\beta\surf$ transversally.
Cutting open $\Cup_\beta\surf$ along $\beta$ thus gives rise to an object 
$\bv_\Gama \iN \Cylo\beta$ and
to a string net on the cut surface. We denote the latter by
  \be
  [\mathrm{cut}(\Gama)] \in \SN(\surf,\overline{\bv_\alpha} \boti \bv_\Gama
  \boti \overline{\bv_\Gama} \boti \bv_{\alpha'})
  \ee
and set $g(\GAMA) \,{:=}\, g^{}_{\bv_\Gama}([\mathrm{cut}(\Gama)])$. By the dinaturality
of the family \eqref{eq:fam:..toV} and the defining invariance of string nets under local
relations, the so obtained linear map \eqref{eq:def:g..toV} is well defined. Moreover,
by construction it is the unique map through which all members $g^{}_{\bvo}$ of the
family \eqref{eq:fam:..toV} factor via the gluing map. This proves universality.
 \\[2pt]
In summary, we have shown that
  \be
  \bearll
  \SN(\Cup_\beta\surf, \overline{\bv_\alpha} \boti \bv_{\alpha'})
  \!\! & \dsty
  \cong \int^ {\bvo\in\Cylo{\beta}}\! \SN(\surf,\overline{\bv_\alpha} \boti \bvo
  \boti \overline{\bvo_{}} \boti \bv_{\alpha'})
  \Nxl1 & \dsty
  \cong \int^ {\bv'\in\Cyl{\beta}}\! \SN(\surf,\overline{\bv_\alpha} \boti \bv'
  \boti \overline{\bv'} \boti \bv_{\alpha'})
  \eear
  \ee
for all $\bv_\alpha \iN \Cyl\alpha$ and $\bv_{\alpha'} \iN \Cyl{\alpha'}$,  
where the second isomorphism holds by Lemmas \ref{lem:5.1.7} and \ref{lem:Kar-pgen}.
This directly upgrades to a coend of functors as in \eqref{eq:coendS} (where we write
an equality instead of an isomorphism because the coend is determined up to unique
isomorphism), and indeed, by exactness of $\SN(\surf,-)$, to a left exact coend.
\end{proof}
	
In particular, when $\surf$ is a disjoint union of two surfaces that are glued 
together, we have

\begin{cor}
Let $\surf \eq \surf_1 \,{\sqcup}\, \surf_2$ be a $($compact oriented smooth$)$ 
surface with a boundary parametrization of the form \eqref{eq:alpha-beta-alpha'} such 
that $\surf_1 \colon \alpha \To \beta$ and $\surf_2 \colon \beta \To \alpha'$
are composable 1-mor\-phisms in $\bordoc$. Then the gluing transformation
  \be
  \SN(\surf_1,\overline ? \boti \bv) \boxtimes \SN(\surf_2,\overline \bv \boti \mq)
  \rarr~ 
  \SN(\Cup_\beta\surf, \overline ? \boti \mq)
  \equiv \SN(\surf_1 {\cup_{\beta}} \surf_2, \overline?\boti\mq)
  \label{eq:sewingS1S2}
  \ee
exhibits the functor $\SN(\surf_1 {\cup_{\beta}} \surf_2, \overline?\boti\mq)$ as 
the $($left exact$)$ coend 
  \be
  \bearll
  \SN(\surf_1 {\cup_{\beta}} \surf_2, \overline?\boti\mq) \!\! &= \dsty
  \int^{\bv'\in\Cyl{\beta}}\!
  \SN(\surf_1 {\sqcup} \surf_2, \overline? \boti \bv' \boti \overline {\bv'} \boti \mq)
  \Nxl1 & \dsty
  = \int^{\bv'\in\Cyl{\beta}}\! 
  \SN(\surf_1,\overline? \boti \bv') \otimes_\ko \SN(\surf_2,\overline{\bv'} \boti \mq) \,.
  \label{eq:coendS1S2}
  \eear
  \ee
As a consequence, the two-functor $\SN$ obtained by the string-net construction 
is an open-closed modular functor $\SN \colon \bordoc \To \Prof$ in the sense of
Definition {\rm \ref{def:mofu}}.
\end{cor}
	
\begin{proof}
The first statement is just the specialization of Theorem \ref{thm:sewing} to the
particular type of surface considered here. The second statement holds because in
this case the gluing can be seen as the composition in $\bordoc$ of the two
1-mor\-phisms $\surf_1 \colon \alpha \To \beta$ and $\surf_2 \colon \beta \To \alpha'$.
\end{proof}

In summary, we have shown that $\SN$ is an open-closed modular functor and that it
satisfies the two requirements that we demand for $\Bl_\C$: When restricted
to the closed sector it is isomorphic to the Turaev-Viro modular functor for \C, and
thus to the Reshetikhin-Turaev modular functor for \ZC, which (conjecturally)
models the conformal blocks for a conformal field theory with chiral data
encoded in \C; and the extension to the open-closed case is implemented by the
functor $\LL$. Accordingly, from now on we take, as already announced, $\Bl_\C$ to
be realized as the string-net modular functor $\SN$.


\subsection{String Nets from World Sheets} \label{sec:SNws} %

To any world sheet $\S$ without physical boundaries we can naturally associate a 
string net in the space $\SNo(\surf_\S,\bvo_\S)$, where $\surf_\S$ is the surface 
that underlies $\S$ and where the boundary value $\bvo_\S$ is obtained from the collection
of boundary data $\bd$, in the sense of Definition \ref{def:bd}, for all geometric 
boundary circles of $\S$, seen as boundary circles of $\surf_\S$. (That is, for each
boundary circle $b\iN \pi_0(\partial\SO)$, take $B|_b$ to be the set of zero-cells of $b$
and define $\bvo_\S(v) \,{:=}\, \bd(v)$ for $v \iN B|_b$.)
We specify the string net associated to $\S$ by the following prescription that yields 
a representative graph $\GamaS$: As an embedded graph, $\GamaS$ consists of all 
one-cells in the interior of $\SO$ as oriented edges and of all zero-cells of $\SO$ 
as vertices. The coloring map $\colE$ from Definition \ref{def:coloring} assigns to
each oriented edge $e$ of $\GamaS$ the object in \C\ that underlies the bimodule label 
of $e$ when seen as a defect line, and the map $\colV$ assigns to each vertex $v$ in the
interior of $\surf_\S$ the morphism label of $v$ when seen as a defect junction on $\S$.

We refer to the graph $\GamaS$, respectively the string net 
  \be
  [\GamaS] \in \SNo(\surf_\S,\bvo_\S)
  \label{eq:[GamaS]}
  \ee
represented by it, as the \emph{partial defect network} on $\surf_\S$ that is 
associated with the world sheet $\S$. Here we include the qualification `partial'
in order to remind us of the fact that the prescription ignores all information
related to the Frobenius algebra labels of the two-cells of the world sheet $\S$.
The partial defect network constitutes an important intermediate step towards the
construction of correlators in Section \ref{sec:SNcor}. 

\begin{rem}
Recall that a defect junction is labeled by a bimodule morphism between objects that 
are of the form $Y\otB Y'$, with $Y$ an $A$-$B$-bimodule, $Y'$ a $B$-$A'$-bimodule and
$\otb$ the tensor product over the Frobenius algebra $B$. In the present context, such
a morphism is to be regarded as a morphism involving the tensor product $\otimes$ in \C\
via the identification ${}_A\Hom_{A'}(X,Y\otB Y') \eq {}_A\Hom_{A'}^{(B)}(X,Y\oti Y')$,
with ${}_A\Hom_{A'}^{(B)}(X,Y\oti Y')$ the subspace of ${}_A\Hom_{A'}(X,Y\oti Y')$
that is invariant under pre-composition with the idempotent $P_{Y \otb Y'}$ as explained
in Appendix \ref{app:M}.
\end{rem}

In the case of a general world sheet, which may possess physical boundaries, the
prescription given above is no longer applicable. Indeed, if $\S$ has at least one 
physical boundary, then the corresponding geometric boundary circle of $\S$ contains a 
labeled and unparameterized one-cell, and for such one-cells the string-net
construction does not provide a suitable counterpart. However, it is not hard to
reformulate the general situation in such a manner that it effectively reduces to the 
one for world sheets without physical boundaries. To this end we convert a world sheet
$\S$ with non-empty physical boundary to a new world sheet $\widetilde\S$, to be
referred to as the \emph{\cws}, that resembles $\S$ but whose physical boundaries are
all labeled by the monoidal unit of \C.  We then also consider an associated
\emph{complemented surface} $\widetilde\surf_\S^{} \,{\equiv}\, \surf_{\widetilde\S}$.

The \cws\ $\widetilde\S$ is constructed as follows. Given a boundary circle $b$ of 
$\S$ containing physical boundaries, but no defect junctions, we partially fill $b$ 
at each physical boundary by a two-cell that is labeled by the monoidal unit $\one$ of
\C. Hereby the circle $b$ is replaced by a new geometric boundary circle $\widetilde b$,
and the former physical boundary, which is labeled by a right $A$-module for some 
Frobenius algebra $A$, is turned into a defect line labeled by the same object, but 
regarded as a $\one$-$A$-bimodule. The boundary of the new $\one$-labeled two-cell is 
the union of this defect line and a new boundary segment $\bb$. 

In case the boundary circle $b$ contains two physical boundaries that meet, together 
with some defect line, at a defect junction on $b$, we partially fill $b$ by a single
new two-cell labeled by $\one$ that is bounded by the two former physical boundaries, 
which have both turned into defect lines, and a single new boundary segment $\bb$.
Analogously we proceed if more than two physical boundaries consecutively meet
at defect junctions. We refrain from writing out the resulting
prescription for the \cws\ $\widetilde\S$ in a formal definition, but content ourselves 
to illustrate what $\widetilde\S$ looks like locally in the case of two physical 
boundaries labeled by $M$ and $N$ that meet at a defect junction:
  \be
  \scalebox{1.0}{\input{MB0.tikz}} ~~~ \xmapsto{~~} ~~~
  \scalebox{1.0}{\input{MB1.tikz}}
  \label{eq:exa:cws}
  \ee

In order that the prescription for the \cws\ complies with the requirements on a world 
sheet in Definitions \ref{def:S0} and \ref{def:S}, the point at which a new defect line 
of $\widetilde\S$ (resulting from a physical boundary of $\S$) meets the circle 
$\widetilde b$ must lie in the interior of a gluing interval. Accordingly we slightly
enlarge the gluing intervals that are adjacent to the new boundary segment 
$\bb \,{\subset}\, \widetilde b$. To this end we regard $\bb$ as the union
of three intervals and declare the two outer ones to be new parts of the enlarged
adjacent gluing intervals. The middle part of $\bb$ must then be regarded as a new
physical boundary; as such it must be labeled by a right $\one$-module, i.e.\ by
an object of \C, and we take the distinguished object $\one \iN \C$ as that label.
Let us illustrate this part of the prescription by redrawing the new boundary segment
on the right hand side of \eqref{eq:exa:cws} in full detail:
  \be
  \scalebox{1.0}{\input{MB2.tikz}}
  \label{eq:exa:cws1}
  \ee
Since the new physical boundaries of $\widetilde\S$ that arise this way all carry the 
same distinguished label $\one$, we can, and will, safely suppress this label in all 
pictures below. Moreover, the presence of these new physical boundaries, as well as of
the adjacent $\one$-labeled two-cells, is in fact entirely irrelevant for the rest of the
construction, for the same reason as we can ignore transparent defect lines. Specifically,
we may treat the geometric boundary circles $\widetilde b$ of $\widetilde\S$ just as
if they were gluing circles. As a consequence, for all practical purposes the \cws\
can be thought of as a world sheet that does not have any physical boundaries; this 
implies in particular that after converting $\S$ to $\widetilde\S$ there is no longer 
any conflict with the string-net construction.
In the pictures below, we account for this fact by suppressing the detailed structure
of the new boundary segments of $\widetilde\S$ that is shown in \eqref{eq:exa:cws1},
and instead just draw the geometric boundary circles $\widetilde b$ in the same 
way as we do for ordinary gluing circles. 

Note that a right $A$-module $N_A$ is the same as a $\one$-$A$-bimodule ${}_\one N_A$, and
hence indeed the one-cells in $\widetilde\S$ coming from physical boundaries in $\S$ are 
properly labeled defect lines. Moreover, a module morphism in 
$\Hom_B(M_B, N_A \otx{A} {}_AX_B)$ is the same as a bimodule morphism in ${}_\one\Hom_B({}
_\one M_B, {}_\one N_A \otx{A} {}_AX_B)$.  Thus $\widetilde\S$ as obtained from $\S$ by our
construction is indeed again a proper world sheet in the sense of Definition \ref{def:S}. 

\begin{exa}
As an illustration, for the world sheet $\S$ displayed in \eqref{eq:pic:exampleS} 
the \cws\ $\widetilde\S$ is given by
  \be 
  \scalebox{1.6}{\input{WS2.tikz}}
  \label{eq:pic:exampleStilde}
  \ee
\end{exa}

Replacing $\S$ by the \cws\ $\widetilde\S$ does not change the field insertions, in the
following sense:

\begin{lem} \label{lem:FF(tildeb)=FF(b)}
Let $\S$ be a world sheet and $b$ be a geometric boundary circle of $\S$ containing a 
physical boundary. Let the boundary circle $\widetilde b$ of $\widetilde\S$ be obtained 
from $b$ by the construction described above, and regard $\widetilde b$ as a gluing 
circle. Then $\FF(\widetilde b) \,{\cong}\, \FF(b)$ as objects in \ZC.
\end{lem}

\begin{proof}
We prove the statement for the case that $b$ is a boundary circle with two (generalized)
boundary field insertions that has the specific form
  \be
  \scalebox{1.0}{\input{FB0.tikz}}
  \label{eq:pic:Ftb=Fb}
  \ee
with $X$ an $A$-$B$-bimodule, $Y$ a $B$-$A$-bimodule, $M$ a right $B$-module and $N$ a
right $A$-module for Frobenius algebras $A$ and $B$. The general case can be treated 
analogously. In the case at hand the circle $\widetilde b$ is given by
  \be
  \scalebox{1.0}{\input{FB1.tikz}}
  \ee
Consider now the following two elements $\bv \eq (\bvo,p)$ and $\widetilde{\bv} \eq 
(\widetilde{\bvo_{\phantom|}},\widetilde p)$ of the cylinder category $\Cyle$:
$\bvo$ consists of two points on $S^1$ that are mapped to 
  \be
  \iHom(M,G^X(N)) = \iHom(M, N{\ota}X) \quad \text{and} \quad
  \iHom(N,G^Y(M)) = \iHom(N, M{\otb}Y) \,,
  \ee
respectively, and $p \eq \id_{\bvo_{\phantom|}}$,
while $\widetilde{\bvo_{\phantom|}}$ consists of six points on $S^1$ that are mapped to
$N$, $X$, $M^\vee$, $M$, $Y$ and $N^\vee$, respectively, and $\widetilde p$ is the
string net represented by
  \be
  \scalebox{1.0}{\input{FB2.tikz}}
  \label{eq:pic:p4Ftb=Fb}
  \ee
(Note that here the two-cells whose boundary contains two segments labeled by $M$ and $N$,
respectively, are the $\one$-labeled two-cells that arose from complementing the world 
sheet.) By a slight variation of Example \ref{exa:pXY} one sees that 
$\Im(\Phi(\widetilde p)) \eq \FF(\widetilde b)$. Moreover, $\bvo$ is isomorphic to
$\widetilde{\bvo_{\phantom|}}$ via the string net
  \be
  \scalebox{1.0}{\input{FB3.tikz}}
  ~~ \in ~ \HomCyle(\bvo,\widetilde{\bvo_{\phantom|}}) \,,
  \ee
where we abbreviate $\iHom$ by $\underline{\mathrm H}$, and 
where the unlabeled boxes represent the canonical morphisms from the internal Homs to 
the corresponding (non-relative) tensor products.

Choosing any equivalence $\Phi_{S^1_{}} \colon \Cyle \Rarr\simeq \ZC$ we thus arrive at
  \be
  \begin{array}{rl}
  \FF(b) = \LL(\iHom(M,G^X(N)) \,{\otimes}\, \iHom(N,G^Y(M)))
  \!\!\! & = \Im(\Phi(\id_{\bvo_{\phantom|}}))
  \Nxl2
  \cong \Phi_{S^1_{}}(\bv)
  \!\!\! & \cong \Phi_{S^1_{}}(\widetilde{\bv}) 
  \cong \Im(\Phi(\widetilde p)) = \FF(\widetilde b) \,.
  \eear
  \ee
Thus $\FF(b) \,{\cong}\, \FF(\widetilde b)$, as claimed.
\end{proof}

Analogously as we did for the world sheet $\S$ at the end of Section \ref{sec:S},
we associate also to the \cws\ $\widetilde\S$ a surface, the complemented surface
$\widetilde\surf_\S^{} \,{\equiv}\, \surf_{\widetilde\S}$. To this end we must
specify a parametrization for $\widetilde\S$. Since when complementing the world
sheet we extend some gluing intervals $r$ to slightly larger intervals
$\widetilde r\,{\supset}\,r$, we cannot just take over a parametrization for
$\S$, but rather specify a diffeomorphism from $I \iN \bordoc$ to $\widetilde r$
instead of to $r$. This modification also guarantees 
that any sewing of the surface $\widetilde\surf_\S$ that is associated with 
$\widetilde\S$ and its parametrization agrees with the surface that arises when first 
sewing the surface associated with the original world sheet $\S$ and its 
parametrization and then complementing, i.e.\ that (with a slight abuse of notation) 
$\Cup_{\widetilde\beta}\widetilde\surf_\S \eq \widetilde{\Cup_\beta\surf_\S}$.

 \medskip

As a consequence of Lemma \ref{lem:FF(tildeb)=FF(b)} the spaces of conformal blocks
for the world sheets $\S$ and $\widetilde\S$ are isomorphic, albeit not canonically
isomorphic.


\subsection{String-net Correlators} \label{sec:SNcor} %

We are now ready to assign to a world sheet $\S$ a string net that realizes a 
correlator for $\S$. We obtain this string net from the partial defect network 
$[\GamaS]$, as constructed in the previous subsection, by embedding into each two-cell 
of the world sheet a suitable graph whose edges are labeled by the Frobenius
algebra label of the two-cell. (If $\S$ has physical boundaries, then we use 
instead the corresponding \cws\ $\widetilde\S$.)

We say that two graphs of the latter type are related by a \emph{Frobenius move}
if the graphs can be transformed to each other by using, locally on the world sheet,
a finite sequence of moves each of which implements a defining relation for the 
structural morphisms of the special symmetric
Frobenius algebra $A$ that labels the two-cell and of the representation morphisms
for the left or right $A$-module structures carried by the defect lines that are
contained in the boundary of the two-cell.

\begin{defi} \label{defi:ffgraph}
Let $\tcl$ be a two-cell of a world sheet and $A(\tcl) \eq (A,\mu,\eta,\Delta,\eps)$
its Frobenius algebra label. A \emph{\fg\ on $\tcl$} is a \C-colored 
oriented graph $\GamaA$ embedded in $\tcl$ having the following features:
 \Itemize
 \item
The vertices of $\GamaA$ are univalent or trivalent. 
Each of the vertices of $\GamaA$ that lies on the boundary $\partial\tcl$ is univalent 
and lies in the interior of a labeled one-cell (i.e.\ of a defect line).
 \item
Each (oriented) edge of $\GamaA$ is labeled by the object $A$. 
 \item
Each trivalent vertex is labeled either by the product $\mu$ or by the coproduct $\Delta$,
and each univalent vertex in the interior of $\tcl$ either by the unit $\eta$ or by the
counit $\eps$. Each univalent vertex on a defect line in $\partial\tcl$ is labeled by
the representation morphism $\rho_X$ or $\ohr_X$ for the left or right $A$-module
structure that is carried by the defect line.
\end{itemize}
~\\[-1.5em]
A \emph{\ffg\ on $\tcl$} is a \fg\ $\GamaA$ on $\tcl$ such that in addition the
following condition is satisfied:
 \Itemize
 \item
$\GamaA$ is \emph{full} in the sense (compare \cite{kmrs}) that adding any further 
edges and vertices in a way compatible with the previous requirements results in
a graph that can be reduced to $\GamaA$ by a Frobenius move.
\end{itemize}
\end{defi}

When drawing a \fg, it is convenient to use the simplified graphical notation
for morphisms involving symmetric Frobenius algebras in which, as explained at the end
of Appendix \ref{app:G}, orientations as well as univalent vertices are omitted.
We refer to the thus obtained simplified version of a (full) \fg\ as a 
\emph{simplified} (full) \fg.

\begin{exa}
The following pictures show an example of a \ffg\ on a two-cell (left) and its 
simplified version (right). 
  \be
  \scalebox{1.0}{\input{SnC0.tikz}} \quad~ \equiv \quad
  \scalebox{1.0}{\input{SnC1.tikz}}
  \label{eq:exa:ffgraph}
  \ee
Note that this \fg\ can be reduced to the one that appears in \eqref{eq:pic:PMAN}
by a suitable Frobenius move. 
\end{exa}

\begin{lem} \label{lem:ffgs}
Any two \ffg s on a two-cell $\tcl \,{\subseteq}\, \S$
are related by a Frobenius move up to isotopy.
\end{lem}

\begin{proof}
The assertion follows by combining the following two statements, each of which is
easy to verify: First, any two realizations of a \emph{simplified} (not necessarily 
full) \fg\ as an ordinary \fg\ are related by a Frobenius move up to isotopy.
Second, any two simplified \ffg s on $\tcl$ can be related by a finite sequence of
the elementary moves that are shown in \eqref{eq:def:a-move}, \eqref{eq:def:b-move}
and \eqref{eq:def:r-move} in Appendix \ref{app:G}.  
\end{proof}

\begin{rem}
(i)
A convenient way to construct a (simplified) \ffg\ $\GamaA$ on a given two-cell $\tcl$
is to remove any number of disks $D_i$ from $\tcl$ and take $\GamaA$ to be the graph that
is obtained from $\tcl {\setminus} \mbox{\large$\cup$}_i D_i$ as a deformation retract.
(This procedure of `punching holes' is similar to the way in which the presence of 
a triangulation of the world sheet in the TFT construction of correlators \cite{fuRs4}
is explained in \cite{kaSau3} and Section 6 of \cite{fusV}.)
 \\[3pt]
(ii) The definition of \ffg\ applies also to the situation that the boundary
$\partial\tcl$ contains a physical boundary labeled by some right $A$-module $M$,
by regarding $M$ as a defect line labeled by $M$ taken as a $\one$-$A$-bimodule.
 \\[3pt]
(iii) If the two-cell $\tcl$ is labeled by the monoidal unit $\one \iN \C$, like e.g.\ any 
of the additional two-cells that result from turning physical boundaries into defect 
lines, then each edge of a full $\one$-graph on $\tcl$ is labeled by $\one$ and each 
vertex by an identity morphism. As a consequence the graph is \emph{transparent} and 
may be replaced by the empty graph, whereby one deals with a situation that has 
already been discussed, from a somewhat different point of view, in \cite{scYa}.
\end{rem}

Next we note that to any boundary circle $b$ of a world sheet $\S$ without physical
boundaries we can naturally associate two string nets on the cylinder over the circle
$b$: First, the string net $[\Gama_{\!\bd}^\circ]$ that directly results from taking the
cylinder over the circle $b$, with boundary value $\bvo \iN \Cylo b$ given by the
boundary datum $\bd$ on $b$, like in the construction of the partial defect network 
$\GamaS$. And second, the string net 
  \be
  p_b^{} := [\Gama_{\!\bd}^\circ \,{\cup}\, \Cup{}_{\tcl} P(\tcl)]
  \ee
on the so obtained cylinder that is obtained by inserting into each of its two-cells 
$\tcl$ the graph $P(\tcl)$ that represents the idempotent for the tensor product over the
Frobenius algebra $A(\tcl)$ which labels the one-cell on $b$ adjacent to $\tcl$.
(In the simplified graphical calculus, $P(\tcl)$ consists of a single $A(\tcl)$-line,
see the picture \eqref{eq:pic:PMAN}.)
The string net $p^{}_b$ is an idempotent in $\Cylo b$ and thus defines an object 
  \be
  \bv = (\bvo,p_b^{})
  \label{bv=bvo,pbd}
  \ee
in the cylinder category $\Cyl b$.
The image of $\bv$ under the equivalence \eqref{Cyl=Z} between $\Cyl b$ and the
Drinfeld center is isomorphic to the field insertion for the circle $b$,
  \be
  \Im(\Phi(p_b^{})) \cong \FF(b) \,.
  \ee
 
\begin{exa} \label{exa:bvo2pbd}
An example already encountered before is the idempotent $p^{}_{X,Y}$ in \eqref{eq:def:pXY},
with image \eqref{eq:ImPhipXY=DXY}.
Another example for an idempotent $p_b$ obtained this way from a boundary circle
$b$ is shown in the following picture:
  \be
  b ~=~~ 
  \scalebox{1.2}{\input{SnC2.tikz}}
  ~~~~\Longrightarrow~~~~~
  p_b^{} ~=~~ 
  \scalebox{1.2}{\input{SnC3.tikz}}
  \label{eq:cylo2cyl}
  \ee
Here the Frobenius algebra labels of the one-cells of $b$ are indicated by different 
colors.
\end{exa}

If $\S$ has physical boundaries, $p_b^{}$ can be defined analogously, after first
replacing $\S$ by the \cws\ $\widetilde\S$.

Combining this observation with Proposition \ref{prop:SN=HomZ} (and recalling the
isomorphism \eqref{eq:cong:Bl}) directly gives

\begin{cor} \label{cor:Bl=SN}
Let $\S$ be a world sheet without physical boundaries. Then there is an isomorphism
  \be
  \Bl_\C(\S) \cong \SN(\surf_\S,\bv)
  \ee
of vector spaces intertwining the action of the mapping class group, where $\surf_\S$
is the surface underlying $\S$ and the boundary value $\bv$ is provided by
the collection of boundary data of $\S$.
\end{cor}

Below we will regard a correlator $\Cor(\S)$ for the world sheet $\S$ as 
a correlator $\Cor(\widetilde\S) \iN \Bl(\widetilde\S)$ for the \cws\ $\widetilde\S$.
In view of Corollary \ref{cor:Bl=SN} we can regard $\Cor(\S)$
as an element of the string-net space $\SN(\surf,\bv)$ with $\surf \eq \surf_\S$ the 
surface underlying $\S$, respectively, in the case of world sheets with physical 
boundaries, with $\surf \eq \surf_{\widetilde\S}$
the surface underlying the \cws\ $\widetilde\S$.
 
\begin{defi} \label{def:CorSN}
(i) Let $\S$ be a world sheet without physical boundaries. Let $\surf_\S$ be the 
surface underlying $\S$ and $\bv$ be the boundary value for $\surf_\S$
that is provided by the collection of boundary data of $\S$.
The \emph{string-net correlator} for $\S$ is the string net 
  \be
  \CorSN(\S) := \big[\,\GamaS \cup \bigcup_{\tcl\subset\S} \GamaA \,\big]
  \in \SN(\surf_\S,\bv)
  \label{eq:CorSN}
  \ee
that is obtained as the union of the partial defect network \eqref{eq:[GamaS]}
with \ffg s on all two-cells $\tcl$ of $\S$.
 \\[3pt]
(ii) Let $\S$ be a world sheet with physical boundaries and $\widetilde\S$ the
\cws\ obtained from $\S$ by turning physical boundaries into defect lines.
Then the string-net correlator is obtained analogously as in (i), but with
$\widetilde\S$ taking over the role of $\S$, i.e.\ $\CorSN(\S) \,{:=}\, 
[\GamaSt \,{\cup}\, \Cup{}_{\!\tcl\subset\widetilde\S\,} \GamaA]$ with $\GamaSt$ the
graph embedded in $\widetilde\S$ that results from the partial defect network $\GamaS$.
\end{defi}

Note that the definition of $\CorSN(\S)$ involves a choice of \ffg\ on each two-cell
of $\S$. However, by Lemma \ref{lem:ffgs} any two such choices are related by a
Frobenius move. Now each Frobenius move corresponds to an equality of morphisms in
\C, and thus also to an equality of string nets. The correlator $\CorSN(\S)$ 
therefore does not depend on these choices.

\begin{exa}
The string-net correlator obtained for the world sheet $\S$ displayed in 
\eqref{eq:pic:exampleS}, for which $\widetilde\S$ is given by 
\eqref{eq:pic:exampleStilde}, is represented by the string net
  \be
  \scalebox{1.6}{\input{WS3.tikz}}  
  \ee
\end{exa}

One of the main results of this 
       book  
is that the string nets \eqref{eq:CorSN}
indeed deserve to be referred to as correlators:

\begin{thm} \label{thm:correlators}
Assigning to any world sheet $\S$ the string-net correlator $\CorSN(\S)$ provides
a consistent system of correlators in the sense of Definition $\ref{def:cosyCor}$.
\end{thm}

\begin{proof}
(i) Invariance under the mapping class group \Maps:
By definition, any $\gamma\iN\Maps$ maps the partial defect network $\GamaS$ on
$\surf_\S$ to itself. Thus when studying the correlator \eqref{eq:CorSN} we only 
have to deal with the \ffg s $\GamaA$ on all two-cells $\tcl \,{\subseteq}\, \S$.
Let thus $\GamaA$ be such a graph and $\gamma\iN\Maps$. Then the graph $\gamma(\GamaA)$
is clearly again a \fg; the following consideration shows that it is even a \ffg:
Let $(\gamma(\GamaA))^+$ be a \fg\ obtained by adding an edge to $\gamma(\GamaA)$. Then
the graph $\gamma^{-1}((\gamma(\GamaA))^+)$ can be obtained from $\GamaA$ by adding an
edge. Since $\GamaA$ is full, it is related to $\gamma^{-1}((\gamma(\GamaA))^+)$ 
by some Frobenius move. Upon applying $\gamma$, that move transports to a Frobenius 
move that relates $\gamma(\GamaA)$ to $(\gamma(\GamaA))^+$. Hence both $\GamaA$ and
$\gamma(\GamaA)$ are \ffg s; by Lemma \ref{lem:ffgs} they are thus related by a 
Frobenius move. Since, as already pointed out, a Frobenius move corresponds to an 
equality of string nets, this implies that the two graphs 
on $\surf_\S$ that are related by replacing $\GamaA$ by $\gamma(\GamaA)$ represent
one and the same string net. Altogether we then have
  \be
  \gamma(\CorSN(\S)) = \big[\,\gamma(\GamaS \cup \bigcup_{\tcl\subset\S} \GamaA) \,\big]
  = \CorSN(\S)
  \ee
for every $\gamma\iN\Maps$.
 \\[3pt]
(ii) Compatibility with sewing: Without loss of generality we can restrict our attention
to what happens when two two-cells $\tcl'$ and $\tcl''$, both labeled by the same
Frobenius algebra, are sewn to a single new two-cell $\tcl \eq \tcl' {\cup_\alpha} \tcl''$,
whereby the \ffg s $\GamaAp$ and $\GamaApp$ combine to a graph $\GamaAp 
{\cup_\alpha} \GamaApp \,{:=}\, \Gama$.
The graph $\Gama$ is clearly a \fg; we must show that it is even a \ffg. To this end,
we add an edge to $\Gama$, resulting in a new graph $\Gama^+$, and show that 
$\Gama^+$ is related to $\Gama$ by a Frobenius move. If the new edge of $\Gama^+$
lies entirely in either $\tcl'$ or $\tcl''$ (regarded as embedded in $\tcl$), then
with obvious notation we have either $\Gama^+ \eq \GamaAp^+ {\cup_\alpha} \GamaApp$ 
or $\Gama^+ \eq \GamaAp {\cup_\alpha} \GamaApp^+$, so that the statement follows
immediately, just because $\GamaAp$ and $\GamaApp$ are both full. Otherwise, i.e.\ if
the new edge lies partly in $\tcl'$ and partly in $\tcl''$, we can perform a
suitable Frobenius move together with isotopy to transform the graph $\Gama$ in 
such a way that we deal again with the previous situation.
 \\[3pt]
(iii) In part (ii) above it is actually implicitly assumed that neither of the
two two-cells $\tcl'$ and $\tcl''$ has arisen from complementing a world sheet.
However, the reasoning still applies to the case of two-cells that result from
complementation: any
\ffg\ in such a two-cell is transparent, so that compatibility with sewing
reduces to the fact that the processes of complementation and of sewing commute.
\end{proof}

  ~ %

\section{Correlators of Particular Interest} \label{sec:specos} %

The prescription presented in the previous chapter applies to all correlators
of the theory. However, a few specific correlators are of particular interest;
concretely, partition functions on the one hand, and correlators which determine 
operator products, i.e.\ composition morphisms on the field objects, on the other
hand. The present chapter provides detailed information about such correlators.
There are two ways of forming a product of defect fields, either along a
defect line or accompanying the fusion of two defect lines.
The former is analyzed in Section \ref{sec:vertical}, and two variants of the latter,
which are related by a braiding, are exhibited in Section \ref{sec:horpros}.
In Section \ref{sec:bulkalgebras} we specialize from defect fields to bulk fields --
the most basic type of field insertions in the bulk, and in Section 
\ref{sec:torusparfu} we discuss the torus partition 
function. The final sections \ref{sec:boOPE} and \ref{sec:buboOPE} of the chapter 
are devoted to the operator product of boundary fields and the
bulk-boundary operator products, respectively.

\subsection{Vertical Operator Product} \label{sec:vertical} %

Consider a gluing circle $b$ with two zero-cells, to which there are attached an incoming
defect line of type $G_1$ and an outgoing defect line of type $G_2$. Recall from Section
\ref{sec:bulkfields} that the field $\FF(b) \iN \ZC$ assigned to $b$ is called 
a \emph{defect field} and is denoted by $\DD^{G_1,G_2}$, and that $G_i$ are
module functors between between \C-modules \M\ and $\M'$.

Let us restrict our attention to the case that the module categories \M\ and $\M'$ are 
indecomposable. Then there
are connected special symmetric Frobenius algebras $A$ and $A'$ in \C\ such that 
$\M \,{\simeq}\, \Mod A$ and $\M' \,{\simeq}\, \Mod A'$, and $A$-$A'$-bimodules
$X_1$ and $X_2$ such that $G_i$ is isomorphic to the functor $G^{X_i} \eq {-} \otA X_i$
for $i \iN \{1,2\}$ (see Appendix \ref{app:M}). We abbreviate $\DD^{G^{X_1},G^{X_2}}$ 
as $\DD^{X_1,X_2}$. Using \eqref{eq:iNat=end} and \eqref{eq:iHom=vee},
the defect field $\DD^{X_1,X_2}$ can thus be written as
  \be
  \bearll
  \DD^{X_1,X_2} = \iNat(G^{X_1},G^{X_2}) \!\! &\dsty 
  =\, \int_{M\in\Mod A} \iHom (M \otA X_1 , M \otA X_2)
  \Nxl3 &\dsty 
  \cong\, (\! \bigoplus_{m\in\I(\Mod A)}\! m \otA X_2 \otAp X_1^\vee \otA m^\vee ,
  \cB_{\DD^{X_1,X_2}}^{}) ~\in \ZC \,.
  \eear
  \label{eq:DDX1X2}
  \ee
The half-braiding $\cB^{}_{\DD^{X_1,X_2}}$ is obtained from the universal coaction of the 
central 
      comonad    %
of \M\ as described in \eqref{eq:Zm-half-braiding}; graphically,
  \be
  \cB^{}_{\DD^{X_1,X_2};Y} ~:= \bigoplus_{m,n\in\I(\Mod A)}\! d_{m} \!\!
  \scalebox{1.0}{\input{VP0.tikz}}
  \quad \text{for} ~Y\iN\C \,.
  \label{eq:half-braiding-DD}
  \ee
Here a summation convention analogous to \eqref{eq_sum-alpha} is used, but with the
implicit $\alpha$-summation being only over module morphisms,

A crucial feature of a full conformal field theory are the \emph{operator products} among
the various types of fields. In the case of defect fields, there are two ways of forming 
a product, either along a given defect line, or such that a fusion of two defect lines is 
involved. In terms of the description \eqref{eq:deffield=iNat} of defect fields, these
correspond to the vertical and horizontal composition of internal natural 
transformations, respectively \cite{fuSc27}. The vertical and horizontal products 
coincide if and only if $\M \eq \M' \eq \C$ and $G_1 \eq G_2 \eq \Id_\C$ 
is the identity functor for \C\ as a bimodule over itself; the corresponding special
defect field $\DD^{\Id,\Id} \iN \ZC$ is known as the object of \emph{bulk fields}.
The purpose of the present subsection is to exhibit how the vertical product 
of defect fields can be related to the string-net construction of
correlators; in the next subsection we will do the same for the horizontal product.

Thinking of defect fields in terms of two-pronged defect junctions, as already visualized
in the picture \eqref{eq:pic:2.19}, the situation on the world sheet that is relevant
for the vertical composition amounts to a sewing operation, according to
  \be
  \sew \Big(\, \scalebox{0.95}{\input{VP13.tikz}}
  \,\sqcup\, \scalebox{0.95}{\input{VP14.tikz}} \,\Big)
  ~~~=~~~ \scalebox{0.95}{\input{VP12.tikz}}
  \label{eq:pic:2.24}
  \ee

It is worth stressing that a string-net correlator as constructed in Section \ref{sec:SNcor}
is directly assigned to a world sheet $\S$, via the underlying surface $\surf_\S$ and the
boundary value for $\surf_\S$ that comes from the boundary data of $\S$. In particular it
does not require the choice of any auxiliary data such as, say, a fine marking of the
surface $\S$ (the latter is, for example,
needed in the Lego-Teichm\"uller based approach of \cite{fuSc22}).
In contrast, relating algebraic data -- like the compositions of internal natural 
transformations that formalize operator products -- to correlators does require such
auxiliary data. Via the string-net construction of correlators one thus achieves a
more invariant description of operator products than what can be seen based on the
underlying purely algebraic structures alone.

\medskip

The vertical composition 
  \be
  \muver \equiv
  \muver(X_1,X_2,X_3) \in \HomZ( \DD^{X_2,X_3}\oti\DD^{X_1,X_2} , \DD^{X_1,X_3})
  \ee
of internal natural transformations is nothing but a particular instance of the 
canonical product of internal Homs \Cite{Def.\,23}{fuSc25}. Just like for ordinary 
natural transformations it amounts to the composition of components, according to
  \be
  \muver ~:=~ \bigoplus_{m\in\I(\M)}
  \raisebox{-5pt}{\scalebox{1.0}{\input{VP1.tikz}}}
  ~ \in \HomZ( \DD^{X_2,X_3}\oti\DD^{X_1,X_2} , \DD^{X_1,X_3}) \,.
  \label{eq:pic:muver}
  \ee

To relate this morphism to a string-net correlator, take the surface to be a 
\emph{pair of pants} $\pop$ -- a sphere with three boundary circles, two of them regarded
as incoming and one as outgoing -- and fix a marking $w$ without cuts (in the sense of
\Cite{Def.\,2.3}{fuSc22}) on $\pop$. Given these data, we will specify a world sheet
$\S_w^{\rm ver} \,{\equiv}\, \S_w^{\rm ver}(X_1,X_2,X_3)$ with underlying surface $\pop$ 
and a boundary value $\bv_w$ on $\pop$ that comes from boundary data for 
$\S_w^{\rm ver}$, as well as a linear isomorphism
  \be
  \varphi_w^{} \Colon \HomZ( \DD^{X_2,X_3}\oti\DD^{X_1,X_2} , \DD^{X_1,X_3})
  \rarr\cong \SN(\pop,\bv_w) \,,
  \label{eq:varphiver}
  \ee
in such a way that $\varphi_w^{}$ maps the vertical composition \eqref{eq:pic:muver} 
to the correlator $\CorSN(\S_w^{\rm ver})$ that the string-net construction yields for 
the world sheet $\S_w^{\rm ver}$, i.e.
  \be
  \varphi_w^{}(\muver) = \CorSN(\S_w^{\rm ver}) \,.
  \label{eq:phi(muver)=SN}
  \ee

We first present the isomorphism \eqref{eq:varphiver} for a particular marking
$w$; the prescription for any other marking without cuts is then obtained via an action of
the mapping class group, as will be described in more detail below. 
Moreover, for describing the isomorphism \eqref{eq:varphiver} it is convenient to split it up
into the composition of two simpler isomorphisms: In a first step, we consider an isomorphism
  \be
  \varphi_w^{\rm can} \Colon \HomZ( \DD^{X_2,X_3}\oti\DD^{X_1,X_2} , \DD^{X_1,X_3})
  \rarr\cong \SN(\pop,\bvcan_w)
  \label{eq:varphican}
  \ee
to the string-net space for a pair of pants with a different boundary value $\bvcan_w$, namely
one that involves the boundary value $\bvcan_Y$, as defined in \eqref{eq:def:bvcan}, for
$Y \iN \ZC$ being
the object $\DD^{X_1,X_2}$, $\DD^{X_2,X_3}$ and $\DD^{X_1,X_3}$, respectively.
The second step then consists of implementing the isomorphism \eqref{eq:bvcan=bvz} between
$\bvcan_{\DD^{X,Y}}$ and the boundary value $\bvz_{X,Y}$ defined in \eqref{eq:def:bvz}.

For a standard pair of pants $\pop$, which we draw as a disk with two holes,
the marking on $\pop$ of our choice is
  \be
  w ~:=~~ \scalebox{0.7}{\input{VP2.tikz}}
  \label{eq:def:varphiver}
  \ee
The corresponding boundary value $\bvcan_w$ is
  \be
  \bvcan_w = 
  (\bvcan_{\DD^{X_2,X_3}})^\vee \boti (\bvcan_{\DD^{X_1,X_2}})^\vee \boti 
  \bvcan_{\DD^{X_1,X_3}} \,.
  \label{eq:bvcanw}
  \ee
The isomorphism $\varphi_w^{\rm can}$ is then defined by
  \be
  \scalebox{1.1}{\input{VP3.tikz}} ~\xmapsto{~\varphi_w^{\rm can}~} \qquad
  \scalebox{1.0}{\input{VP4.tikz}}
  \label{eq:def:varphican}
  \ee
Hereby the vertical composition $\muver$ is mapped to the string net 
  \be
  \bearll
  \varphi_w^{\rm can}(\muver) 
  & = \!\bigoplus_{\scriptstyle i,j,k\in\I(\C) \atop \scriptstyle m,n,p,q\in\I
             (\Mod A)
  } \!\!\!\! \frac{\dd_i\, \dd_j\, \dd_k\, \dd_m\, \dd_n\, \dd_p}{\DC^6}~~
  \scalebox{1.0}{\input{VP8.tikz}} 
  \\[-0.2em]
  & = \sum_{m,n,q\in\I
             (\Mod A)
  } \! \frac{\dd_m\, \dd_n}{\DC^4}~~
  \scalebox{1.0}{\input{VP9.tikz}} 
  \\[-0.5em]~
  \eear
  \label{eq:varphican(muver)}
  \ee
in $\SN(\pop,\bvcan_w)$. Here the second equality holds as a consequence of the
identity \eqref{eq:proofshrink2} and Corollary \ref{cor:shrink}.
(Also, in the first picture -- and likewise in several other pictures later on -- we 
omit, for lack of space, the summation labels of the
module mor\-phisms; the appropriate pairings are instead indicated by matching colors.)
The second step in the construction of $\varphi_w^{}$ then amounts to setting
  \be
  \varphi_w^{}(-) := e^{}_{X_1,X_3} \circ \varphi_w^{\rm can}(-) \circ
  ( r^{}_{X_2,X_3} \,{\boti}\, r^{}_{X_1,X_2} ) \,,
  \label{eq:e.-.rr}
  \ee
with the morphisms $e_{-,-}$ and $r_{-,-}$ as defined in \eqref{eq:eXYrXY}.
Accordingly, the boundary value $\bv_w$ is given by
  \be
  \bv_w = (\bvz_{X_2,X_3})^\vee \boti (\bvz_{X_1,X_2})^\vee \boti \bvz_{X_1,X_3} \,.
  \label{eq:bvw}
  \ee
Via \eqref{eq:e.-.rr}, the string net \eqref{eq:varphican(muver)} gets mapped to
  \be
  \varphi_w^{}(\muver) ~=~~
  \scalebox{1.0}{\input{VP10.tikz}} ~~\in \SN(\pop,\bv_w)
  \label{eq:varphiver(muver)}
  \ee
 
We can now read off that we have indeed achieved to express $\varphi_w^{}(\muver)$ as a
string-net correlator $\CorSN(\S_w^{\rm ver})$, namely as the one for the world sheet
  \be
  \S_w^{\rm ver} ~:=~~
  \scalebox{0.9}{\input{VP11.tikz}}
  \label{eq:Sver}
  \ee
that is, a pair of pants with three defect lines labeled by
            $X_1$, $X_2$ and $X_3$, 
respectively, which pairwise connect the boundary circles.
\begin{rem}	   
By a diffeomorphism, the world sheet \eqref{eq:Sver} can be redrawn as follows:
  \be
  \S_w^{\rm ver} ~=~~
  \scalebox{0.9}{\input{VP15.tikz}}
  \ee

This way of presenting $\S_w^{\rm ver}$ slightly obscures its relevance for the vertical
operator product. On the other hand, it clarifies the relation to other approaches: It
makes it obvious that one deals with an `operator product along a defect line', and it
is precisely what commonly is called the \emph{fundamental world sheet} for three defect 
fields on the sphere, see e.g.\ Section 4.5 of \cite{fuRs10}.
\end{rem}	   

 \medskip

We finally comment on the dependence of the construction on the marking of the surface.
A different choice $w'$ of a marking without cuts on the pair of pants leads to a string
net which, in general, differs from \eqref{eq:varphiver(muver)}. The two world sheets
are related by the unique element $\gamma_{w.w'}$ of the mapping class 
group of $\surf$ that corresponds (see \Cite{Sect.\,3.1}{fuSc22}) to a move of markings 
mapping $w$ to $w'$. In particular, two markings $w$ and $w'$ without cut and with the
same end points on the boundary circles give the same world sheet, and thus the same 
correlator, if and only if they are isotopic (in this case
$\gamma_{w.w'}$ belongs to both $\Map(\S_w^{\rm ver})$ and $\Map(\S_{w'}^{\rm ver})$,
and the two groups actually coincide).

As a particular case, any two markings without cuts give the same correlator when each
of the three defect fields is actually a bulk field, i.e.\ when 
$X_1 \eq X_2 \eq X_3 \eq A \eq A'$.


\subsection{Horizontal Operator Products} \label{sec:horpros} %

We now analyze the horizontal composition of defect fields, in analogy with the 
study of the vertical composition in Section \ref{sec:vertical}.
Horizontal composition happens in combination with the \emph{fusion of defect lines},
in which two parallel segments of defect lines get replaced by a single one. That
such a fusion process is possible is an integrated ingredient of the description of
defect lines in other approaches to CFT. In the present setting, fusion
is algebraically realized as the tensor product over the 
relevant Frobenius algebra: the fusion of two defect lines that are labeled by
an $A$-$A'$-bimodule $X$ and an $A'$-$A''$-bimodule $X'$, respectively,
yields a defect line labeled by the $A$-$A''$-bimodule $X {\otimes_{\!A'_{}}} X'$.

In the same way as done in \eqref{eq:pic:2.24} for the vertical operator product, the 
horizontal composition can be expressed as a specific sewing operation, according to
  \be 
  \sew \Big(\, \scalebox{0.95}{\input{HP1.tikz}}
  \,\sqcup\, \scalebox{0.95}{\input{HP2.tikz}} \,\Big)
  ~~~=~~~ \scalebox{0.95}{\input{HP0.tikz}}
  \label{eq:pic:2.25}
  \ee

When attempting to translate this picture into an expression for composition
  \be
  \DD^{X_2,X_4} \oti \DD^{X_1,X_3} \rarr~ \DD^{X_1\otimes_{\!A'}X_2,X_3\otimes_{\!A'}X_4}
  \ee
of internal natural transformations, some care is needed. To explain the subtleties
involved, it is convenient to recall first how the horizonal composition of 
\emph{ordinary} natural transformations is described in components: The horizonal product
of two natural transformations $\dd^{G_1,G_3}$ between functors $G_1,G_3\colon \M \To \M'$
and $\dd^{G_2,G_4}$ between functors $G_2,G_4\colon \M' \To \M''$ -- a natural
transformation from $G_2 \cir G_1$ to $G_4 \cir G_3$, which are functors from \M\ to
$\M''$ -- amounts to a suitable composition of their components, i.e.\ of
morphisms $\dd^{G_1,G_3}_M \iN \Hom_{\M'}(G_1(M),G_3(M))$ and $\dd^{G_2,G_4}_{M'}
\iN \Hom_{\M''}(G_2(M'),G_4(M'))$, respectively. In more detail, the composition
can be expressed both as
  \be
  G_2 \cir G_1(M) \rarr{G_2^{}(\dd^{G_1,G_3}_M)} G_2 \cir G_3(M) 
  \rarr{\dd^{G_2,G_4}_{G_3(M)}} G_4 \cir G_3(M) \,,
  \label{eq:horNat1}
  \ee
and as 
  \be
  G_2 \cir G_1(M) \rarr{\dd^{G_2,G_4}_{G_1(M)}} G_4^{} \cir G_1(M) 
  \rarr{G_4(\dd^{G_1,G_3}_M)} G_4 \cir G_3(M) \,.
  \label{eq:horNat2}
  \ee
Equality of the two composite morphisms \eqref{eq:horNat1} and \eqref{eq:horNat2}
for all $M \iN \M$ holds by naturality of $\dd^{G_2,G_4}$.
Now in the case of \emph{internal} natural transformations
$\DD^{X_1,X_3}$ and $\DD^{X_2,X_4}$, the components are internal Homs, and accordingly 
it is appropriate to interpret the dinatural structure morphisms 
  \be
  \jj^{}_M \equiv \jj^{F,F'}_M \Colon \iNat(F,F') \rarr~ \iHom_{\N}(F(M),F'(M))
  \label{eq:iNat-jj}
  \ee
of the end \eqref{eq:iNat=end} as the `projection to components'. The analogues
of the morphisms \eqref{eq:horNat1} and \eqref{eq:horNat2} are then the composites
  \be
  \bearl \iNat(G_2,G_4) \otimes \iNat(G_1,G_3)
  \Nxl2 \hsp 3
  \rarr{ \jj^{G_2,G_4}_{G_3(M)} \otimes \jj^{G_1,G_3}_{M} }
  \iHom(G_2 \cir G_3(M),G_4\cir G_3(M)) \otimes \iHom(G_1(M),G_3(M))
  \Nxl2 \hsp 3
  \rarr{ \id \otimes \iG_2 } 
  \iHom(G_2 \cir G_3(M),G_4\cir G_3(M)) \otimes \iHom(G_2\cir G_1(M),G_2\cir G_3(M))
  \Nxl3 \hsp 3
  \rarr{~\imu~} \iHom(G_2\cir G_1(M),G_4\cir G_3(M))
  \eear
  \label{eq:horiNat1}
  \ee
and
  \be
  \bearl \iNat(G_1,G_3) \otimes \iNat(G_2,G_4) 
  \Nxl2 \hsp 3
  \rarr{ \jj^{G_1,G_1}_{M} \otimes \jj^{G_2,G_4}_{G_1(M)} }
  \iHom(G_1(M),G_3(M)) \otimes \iHom(G_2 \cir G_1(M),G_4\cir G_1(M))
  \Nxl2 \hsp 3
  \rarr{ \iG_4  \otimes \id }
  \iHom(G_4\cir G_1(M),G_4\cir G_3(M))\otimes \iHom(G_2 \cir G_1(M),G_4\cir G_1(M))
  \Nxl3 \hsp 3
  \rarr{~\imu~} \iHom(G_2\cir G_1(M),G_4\cir G_3(M))
  \eear
  \label{eq:horiNat2}
  \ee
for $M \iN \M$, respectively.
Here $\imu$ is the standard multiplication \eqref{eq:imu} on internal Homs,
while the morphisms $\iG_2$ and $\iG_4$ are defined as in \eqref{eq:def:iG}.

It is straightforward to see that both of the families
\eqref{eq:horiNat1} and \eqref{eq:horiNat2} are dinatural in $m \iN \M$. 
Owing to the universal property of $\iNat$ as an end, they thus factorize to morphisms
  \be
  \muhore \Colon \iNat(G_2,G_4) \otimes \iNat(G_1,G_3) \rarr~ \iNat(G_2\cir G_1,G_4\cir G_3)
  \phantom{\,,}
  \label{eq:muhore}
  \ee
and
  \be
  \muhorz \Colon \iNat(G_1,G_3) \otimes \iNat(G_2,G_4) \rarr~ \iNat(G_2\cir G_1,G_4\cir G_3) \,,
  \label{eq:muhorz}
  \ee
respectively. It can further be shown that \eqref{eq:muhore} and \eqref{eq:muhorz},
which are defined as morphisms in \C, are in fact even morphisms in \ZC.
When expressing the module functors $G_i$ through bimodules $X_i$ as in \eqref{eq:GY=},
we can write the morphisms \eqref{eq:muhore} and \eqref{eq:muhorz}
in terms of defect fields as
  \be
  \muhore \Colon \DD^{X_2,X_4} \otimes \DD^{X_1,X_3} \rarr~ 
  \DD^{X_1\otimes_{\!A'}X_2,X_3\otimes_{\!A'}X_4} \phantom{\,.}
  \label{eq:muhore2}
  \ee
and
  \be
  \muhorz \Colon \DD^{X_1,X_3} \otimes \DD^{X_2,X_4} \rarr~ 
  \DD^{X_1\otimes_{\!A'}X_2,X_3\otimes_{\!A'}X_4} 
  \label{eq:muhorz2}
  \ee
respectively.

We call the morphisms $\muhore$ and $\muhorz$ the left and right horizontal operator 
product. The reason for this choice of terminology is the description of 
\eqref{eq:muhore2} and \eqref{eq:muhorz2} in terms of (ordinary, non-cyclic)
string diagrams, in which the basis morphisms that are summed over are 
located in the left and right half of the graph, respectively: we have
  \be
  \muhore ~=
  \!\bigoplus_{\scriptstyle m\in\I
             (\Mod A')
  \atop \scriptstyle n\in\I
             (\Mod A)
  }\!\! \dd_m ~~
  \scalebox{1.2}{\input{HP3.tikz}}
  \ee
and
  \be
  \muhorz ~=
  \!\bigoplus_{\scriptstyle m\in\I
             (\Mod A')
  \atop \scriptstyle n\in\I
             (\Mod A)
  }\!\! \dd_m ~~
  \scalebox{1.2}{\input{HP4.tikz}}
  \ee

Let us now express these two compositions through string nets. Recall from Section
\ref{sec:vertical} that this requires the specification of
auxiliary data beyond the structure of a world sheet, concretely the choice of a
marking without cuts. Note that these auxiliary data do not appear in the purely
algebraic treatment of horizontal composition that is given in \cite{fuSc25}.
(Also, in our discussion of the horizontal composition we partly deviate from the
exposition in \cite{fuSc25}.) 

In the same vein as done for the vertical composition $\muver$, we will now determine
the image of the morphism $\muhore \iN \HomZ (\DD^{X_2,X_4} \oti \DD^{X_1,X_3},
\DD^{X_1\otimes_{\!A'}X_2,X_3\otimes_{\!A'}X_4})$ under the composite map
  \be
  \begin{array}{r}
  \varphi_w^{} \Colon \HomZ( \DD^{X_2,X_4} \oti \DD^{X_1,X_3},
  \DD^{X_1\otimes_{\!A'}X_2,X_3\otimes_{\!A'}X_4})
  \rarr{~\varphi_w^{\rm can}} \SN(\pop,\bvcan_w) \phantom{\,,}
  \Nxl2
  \rarr{ e^{}_{\!X_1\otimes_{\!A'}X_2,X_3\otimes_{\!A'}X_4} \circ \,(-)\, \circ\,
  ( r^{}_{\!X_2,X_4} \boxtimes r^{}_{\!X_1,X_3}) } \SN(\pop,\bv^{}_w) \,,
  \eear
  \label{eq:varphiw2}
  \ee
Here $w$ is the marking \eqref{eq:def:varphiver} on the pair of pants $\pop$, the
isomorphism $\varphi_w^{\rm can}$ is defined analogously as in \eqref{eq:def:varphican},
and the boundary values are
  \be
  \bearll
  ~ & \bvcan_w =
  (\bvcan_{\DD^{X_2,X_4}})^\vee \boti (\bvcan_{\DD^{X_1,X_3}})^\vee \boti
  \bvcan_{\DD^{X_1\otimes_{\!A'}X_2,X_3\otimes_{\!A'}X_4}}
  \Nxl3
  {\rm and}\quad &
  \bv_w = (\bvz_{X_2,X_4})^\vee \boti (\bvz_{X_1,X_3})^\vee \boti
  \bvz_{X_1\otimes_{\!A'}X_2,X_3\otimes_{\!A'}X_4}
  \eear
  \ee
analogously as in \eqref{eq:bvcanw} in \eqref{eq:bvw}.

By direct calculation we obtain
  \be
  \bearll
  \varphi_w^{\rm can}(\muhore) & =~~
  \scalebox{1.0}{\input{HP5.tikz}}
  \\[-1.9em] &\dsty
  =\! \sum_{{\scriptstyle i,j,k\in\I(\C) \atop \scriptstyle m,p\in\I(\Mod A')}
  \atop \scriptstyle n,q,r\in\I(\Mod A'')} \!\!\!
  \frac{ \dd_{i}\, \dd_{j}\, \dd_{k}\, \dd_{m}\, \dd_{n}\, \dd_{p}\, \dd_{q} } {\DC^6}~~
  \scalebox{1.0}{\input{HP6.tikz}}
  \\~\\[0.1em]
  &\dsty
  =\! \sum_{\scriptstyle m,p\in\I(\Mod A') \atop \scriptstyle n,q,r\in\I(\Mod A'')} \!\!\!
  \frac{ \dd_{m}\, \dd_{n}\, \dd_{p}\, \dd_{q} }{\DC^6}~~
  \scalebox{1.0}{\input{HP7.tikz}}
  \\[-0.9em] &\dsty \hspace*{12em}
  =\! \sum_{\scriptstyle m\in\I(\Mod A') \atop \scriptstyle n,r\in\I(\Mod A'')} \!\!
  \frac{ \dd_{m}\, \dd_{n} }{\DC^4}~~
  \scalebox{1.0}{\input{HP8.tikz}}
  \\[-1.8em] ~
  \eear
  \ee
~\\[-0.4em]
This implies that $\varphi_w^{}(\muhore)$ gives the world sheet for the left 
horizontal composition as follows:
  \be
  \begin{array}{r}
  \varphi_w^{}(\muhore) ~=~ e^{}_{\!X_1\otimes_{\!A'}X_2,X_3\otimes_{\!A'}X_4} \circ
  \varphi_w^{\rm can}(\muhore) \circ ( r^{}_{\!X_2,X_4} \boti r^{}_{\!X_1,X_3}) 
  \hspace*{8.5em}
  \Nxl4
  =~~~ \scalebox{1.0}{\input{HP9.tikz}}
  =~ \CorSN(\Shore)
  \eear
  \label{eq:CorSNShore}
  \ee
with
  \be
  \Shore ~=~~
  \scalebox{0.9}{\input{HP10.tikz}}
  \label{eq:def:Shore}
  \ee
~\\[0.2em]
Likewise, the image of $\muhorz$ under $\varphi_w^{}$ gives the world sheet for the
right horizontal composition:
  \be
  \varphi_w^{}(\muhorz) ~=~ \CorSN(\Shorz) \qquad{\rm with}\qquad
  \Shorz ~=~~
  \scalebox{0.9}{\input{HP11.tikz}}
  \label{eq:varphi.Shorz}
  \ee

\begin{rem} \label{rem:4.2}
(i) Up to isotopy we have
  \be
  \Shore = \beta(\Shorz) \,,
  \label{eq:Shore=beta.Shorz}
  \ee
where $\beta$ is the action of the \emph{braid move} \cite{BAki}, which 
(when expressed in terms of markings) acts as
  \be
  \scalebox{0.75}{\input{VP2.tikz}}
  ~~\xmapsto{~~~}~~
  \scalebox{0.75}{\input{HP12.tikz}}
  \ee
As we will observe in Remark \ref{rem:muhore=cB.muhorz}
below, this implies that the left and right horizontal compositions $\muhore$ and 
$\muhorz$ merely differ by a half-braiding.
 \\[3pt]
(ii) 
By applying suitable diffeomorphisms, both world sheets $\Shore$ and $\Shorz$
may be redrawn as
  \be
  \S_{\rm hor}^{} ~=~~
  \scalebox{0.9}{\input{HP13.tikz}}
  \ee
which is the form familiar from the sewing operation shown in picture \eqref{eq:pic:2.25}.
\end{rem}


\subsection{Bulk Algebras} \label{sec:bulkalgebras}  %

The most basic type of field insertion in the bulk is the one obtained for a gluing 
circle $b \eq b_{\rm bulk}^A$ of the following form: $b \eq b_{\rm bulk}^A$ has two 
zero-cells, and the two adjacent two-cells are both labeled by one and the same
Frobenius algebra $A$, while the two defect lines attached to the zero-cells are both
transparent, i.e.\ labeled by the functor $G^A \eq {-} \,{\ota} A$, which is 
canonically isomorphic to the identity functor $\Id_{\Mod A}$. Thus 
$b \eq b_{\rm bulk}^A$ is a special case of the gluing circle describing defect fields,
which has already been displayed in the picture \eqref{eq:pic:2.19}:
  \be 
  b_{\rm bulk}^A ~=~~
  \scalebox{1.0}{\input{BA0.tikz}}
  \label{eq:pic:2.19-bulk}
  \ee
The resulting field insertion is the bulk field object
  \be
  \DD^{A,A} \equiv \FF(b_{\rm bulk}^A) \,=\, \iNat(G^A,G^A)
  \,\cong\, \iNat(\Id_{\Mod A},\Id_{\Mod A}) \,= \int_{\!M \in \Mod A} \iHom(M,M) \,.
  \ee
The end $\int_{\!M} \iHom(M,M)$ is the \emph{full center} $Z(A) \iN \ZC$ 
\cite{ffrs,davy20} of the algebra $A \iN \C$. In the
semisimple case of our interest, we thus have
  \be
  \DD^{A,A} \cong Z(A)
  = (\! \bigoplus_{m\in\I(\Mod A)}\! m \otA m^\vee , \cB^{}_{Z(A)} ) ~\in\ZC \,.
  \label{eq:bulkfields}
  \ee
The half-braiding is the one already described in 
\eqref{eq:half-braiding-DD}, which now specializes to
  \be
  \cB^{}_{Z(A);Y} ~= \bigoplus_{m,n\in\I(\Mod A)}\! d_{m} ~
  \scalebox{1.1}{\input{BA1.tikz}}
  \qquad \text{for} ~Y\iN\C \,.
  \ee

As already mentioned in Section \ref{sec:bulkfields}, the special case 
\eqref{eq:bulkfields} of the defect fields \eqref{eq:DDX1X2} is known as the
\emph{bulk fields} of a CFT.
For bulk fields, the vertical product $\muver(A,A,A)$ as defined in 
\eqref{eq:pic:muver} and the left and right horizontal products $\muhore(A,A,A)$
and $\muhorz(A,A,A)$ as defined in \eqref{eq:muhore} and in \eqref{eq:muhorz},
respectively, all coincide. The canonical isomorphism
$\DD^{A,A} \,{\cong}\, Z(A)$ transports each of them to the morphism
  \be
  \mu_{Z(A)}^{} = \bigoplus_{m\in\I(\Mod A)} ~
  \scalebox{1.0}{\input{BA2.tikz}}
  \label{eq:mu(Z(A))}
  \ee

By the same strategy that we already pursued in the more general case of defect fields
we can, after fixing a marking $w$ without cuts on the pair of pants $\pop$, map
the product \eqref{eq:mu(Z(A))} to the string-net correlator on a world sheet $\S$:
as a special case of \eqref{eq:phi(muver)=SN} we get
  \be
  \varphi_w^{}(\mu_{Z(A)}
	     )
  = \CorSN(\S) ~\in \SN(\pop,\bv_{Z(A),w}^{})
  \ee
with
  \be
  \S \,\equiv\, \S_w^{} ~:=~~
  \scalebox{0.9}{\input{BA3.tikz}}
  \ee
 ~

Note that the so obtained world sheet $\S$ is nothing but any of the world sheets
$\S_w^{\rm ver}$, $\Shore$ and $\Shorz$ as displayed in \eqref{eq:Sver}, 
\eqref{eq:def:Shore} and \eqref{eq:varphi.Shorz}, respectively, each specialized to
the case that the two two-cells are both labeled by the same Frobenius algebra and the
three defect lines are all transparently labeled. Owing to their transparency the 
defect lines can be omitted from the world sheet, and the boundary value 
$\bv_{Z(A),w}^{}$ -- a choice of points on $S^1 \,{\sqcup}\, S^1 \,{\sqcup}\, S^1$ 
at which the transparent defect lines start or end -- is immaterial. For the same 
reason, the choice of marking $w$ is actually irrelevant as well.
It follows that $\mu_{Z(A)}$ endows the object $Z(A) \iN \ZC$ with the structure of 
a commutative special symmetric Frobenius algebra, and that the correlator $\CorSN(\S)$
is invariant under the mapping class group $\Map(\pop)$ of the three-holed sphere. 
Further, invoking Theorem 3.4 of \cite{koRu2}, it follows that $Z(A)$
is even a \emph{modular} Frobenius algebra in the sense of \Cite{Def.\,4.9}{fuSc22}.

\begin{rem}
The Frobenius algebra $Z(\one)$ that is obtained when $A$ is the monoidal unit $\one$
of \C\ is the \emph{Cardy bulk algebra} considered in \cite{scYa}. It should be noted
that in \cite{scYa} a less natural boundary value 
$\bvcan_{\tilde F} \eq (S_{\tilde F},\pcan_{\tilde F})$, was considered than here,
consisting of a circle with a point labeled by $\tilde F$ together with the idempotent
  \be
  \pcan_{\tilde F} ~:=~~ 
  \scalebox{0.9}{\input{BA4.tikz}}
  \ee
This choice results in more complicated graphs than the ones appearing above, having
additional edges around the boundary circles. However, the boundary values 
$\bv_{Z(A),w}^{}$ and $\bvcan_{\tilde F}$ are in fact isomorphic, and as a consequence
the present description and the one used in \cite{scYa} indeed yield the same correlator
of three bulk fields on the sphere.
\end{rem}


\subsection{Torus Partition Function} \label{sec:torusparfu}  %

Next we consider the correlator $\Cor(\mathrm T)$ for a torus $\mathrm T$ without
field insertions and without non-transparent defect lines. This correlator, commonly
called the \emph{torus partition function}, is of much interest. For instance, in
rational CFT it allows one to directly read off the decomposition of the bulk field 
object $\DD^{A,A} \eq \FF(b_{\rm bulk}^A)$ into simple objects of \ZC. Also, modular
invariance of $\Cor(\mathrm T)$ is an important constraint on the consistency of a
full conformal field theory, to the extent that often it has
been assumed, erroneously, to be even a sufficient condition for consistency.

The string-net form of the torus partition function follows immediately from 
Definition \ref{def:CorSN}: we have
  \be
  \CorSN(\mathrm T) ~=~~
  \scalebox{1.0}{\input{TP0.tikz}}
  \label{eq:CorSN(T)}
  \ee
Using the expression \eqref{eq:bulkfields} for the bulk field object 
$\DD^{A,A} \,{\cong}\, Z(A)$ together with Corollary \ref{cor:shrink} and the identity
\eqref{eq:proofshrink2}, this can be rewritten as
  \be
  \CorSN(\mathrm T) ~=~
  \!\!\! \bigoplus_{\scriptstyle k\in\I(\C) \atop \scriptstyle m\in\I
             (\Mod A)
  } \!\!\! \frac{\dd_k\, \dd_m}{\DC^2}~
  \scalebox{1.0}{\input{TP1.tikz}}
  ~~=~~
  \scalebox{1.0}{\input{TP2.tikz}} ~~
  \ee
Moreover, since \C\ is semisimple the object $Z(A) \iN \ZC \,{\simeq}\; \C\rev \boti \C$
can be decomposed into a direct sum of simple objects of \ZC\ as
  \be
  Z(A) \cong \bigoplus_{i,j\in\I(\C)}\! \mathrm Z(A)_{i,j} \otik \GC(i \boti j) \,.
  \label{eq:ZA=oplus}
  \ee
Here $\GC$ the equivalence \eqref{eq:def:GC} between $\C\rev \boti \C$ 
and \ZC, while $\mathrm Z(A)_{i,j}$ are the vector spaces
  \be
  \bearll
  \mathrm Z(A)_{i,j} \!\!\!&
  := \HomZ(\GC(i \boti j) , \iNat(G^A,G^A))
  \Nxl2 &
  \,\cong \Hom_{\Funre_\C(\Mod A,\Mod A)} \big( \GC(i \boti j) \act G^A,G^A \big)
  \Nxl3 &
  \,\cong \Hom_{A|A} (i \,{\otimes^-} A \,{\otimes^+} j\,, A) \,.
  \eear
  \ee
This in particular reproduces the decomposition rule \eqref{eq:HomAA()} for 
bulk fields. For the torus partition function we thus get
  \be
  \CorSN(\mathrm T) ~=~ \sum_{i,j\in\I(\C)} \mathrm z(A)_{i,j} ~~
  \scalebox{1.0}{\input{TP3.tikz}}
  \label{eq:CorSN(T)v3}
  \ee
with $\mathrm z(A)_{i,j} \,{:=}\, \mathrm{dim}_\complex(\mathrm Z(A)_{i,j})$.

\begin{rem}
The \ffg\ in the string-net correlator \eqref{eq:CorSN(T)}
is the same as the one appearing in the TFT construction of $\Cor(\mathrm T)$.
In that case, the torus in which the graph is embedded is the subset 
$\{0\} \Times \mathrm T$ of the three-manifold $[-1,1] \Times \mathrm T$
(see \Cite{Eq.\,(5.24)}{fuRs4}).
An analogous relationship between the two constructions holds for all other
partition functions, i.e.\ for the correlator of any world sheet that neither has
field insertions nor contains any physical boundaries or non-trivial defect lines.
\end{rem}

The expression \eqref{eq:CorSN(T)} for $\CorSN(\mathrm T)$ shows manifestly
that the torus partition function is invariant under the geometric action of 
the modular group $\Map(\mathrm T) \,{\cong}\, \mathrm{PSL}(2,\mathbb Z)$ of 
the torus. The result \eqref{eq:CorSN(T)v3} allows us to translate
this geometric modular invariance to the algebraic modular invariance of the
$\I(\C){\times}\I(\C)$-matrix $\mathrm z(A) \eq (\mathrm z(A)_{i,j})$. To
see this, we first note that the set $\{ G_{i,j} \,|\, i,j\iN\I(\C) \}$ with
  \be
  G_{i,j} ~:=~ \sum_{k\in\I(\C)} \frac{\dd_k}{\DC^2} ~~
  \scalebox{0.75}{\input{MT3.tikz}}
  \ee
is a basis of the string-net space $\SN(\mathrm T)$ \Cite{Prop.\,4.8}{runk15}.
It is therefore sufficient to show that the image 
  \be
  S(G_{i,j}) ~=~ \sum_{k\in\I(\C)} \frac{\dd_k}{D^2} ~~
  \scalebox{0.75}{\input{MT4.tikz}}
  \ee
of $G_{i,j}$ under the modular $S$-trans\-formation satisfies
  \be
  S(G_{i,j})
  = \sum_{i',j'\in\I(\C)}\! \frac{s^{}_{i^\vee,i'}\,s^{}_{j,j'}}{\DC^2}\, G_{i',j'} 
  \qquad{\rm with}\qquad s_{i,j}^{} ~:=~
  \scalebox{0.85}{\input{HL.tikz}}
  \label{eq:s(Gij)}
  \ee
The validity of \eqref{eq:s(Gij)} is established by the following chain of equalities:
  \begin{eqnarray}
  && \hspace*{-1.5em}
  \sum_{i',j'\in\I(\C)}\! \frac{s^{}_{i^\vee,i'}\,s^{}_{j,j'}}{\DC^2}\, G_{i',j'}
  ~= \sum_{i',j',k\in\I(\C)}\!\! \frac{s^{}_{i^\vee,i'}\,s^{}_{j,j'}\,\dd_k}{\DC^4}~~
  \scalebox{0.75}{\input{MT5.tikz}}
  \nonumber
  \\&&~\nonumber\\[0.7em] && \hspace*{2.6em}
  = \sum_{i',j',k\in\I(\C)}\!\! \frac{\dd_{i'}\,\dd_{j'}\,\dd_k}{\DC^4}~~
  \scalebox{0.75}{\input{MT6.tikz}}
  ~~= \sum_{i',j',k,l\in\I(\C)}\!\! \frac{\dd_{i'}\,\dd_{j'}\,\dd_k\,\dd_l}{\DC^4}~~
  \scalebox{0.75}{\input{MT7.tikz}} ~~~
  \nonumber
             \end{eqnarray}
             \begin{eqnarray}
  &&~\nonumber\\[0.8em] && \hspace*{2.6em}
  = \sum_{i',j',l\in\I(\C)}\!\! \frac{\dd_{i'}\,\dd_{j'}\,\dd_l}{\DC^4}~~
  \scalebox{0.75}{\input{MT8.tikz}}
  ~~= \sum_{i',j',l,m\in\I(\C)}\!\! \frac{\dd_{i'}\,\dd_{j'}\,\dd_l\,\dd_m}{\DC^4}~~
  \scalebox{0.75}{\input{MT9.tikz}}
  \nonumber
  \\&&~\nonumber\\[0.8em] && \hspace*{2.6em}
  = \sum_{i',l,m\in\I(\C)}\!\! \frac{\dd_{i'}\,\dd_l\,\dd_m}{\DC^4}~~
  \scalebox{0.75}{\input{MT10.tikz}}
  ~~= \sum_{l,m\in\I(\C)}\!\! \frac{\dd_l\,\dd_m}{\DC^2}\, \delta_{l,0}~~
  \scalebox{0.75}{\input{MT11.tikz}}
  \nonumber
  \\&&~\nonumber\\[0.8em] && \hspace*{14em}
  = \sum_{m\in\I(\C)} \frac{\dd_m}{\DC^2}~~
  \scalebox{0.75}{\input{MT12.tikz}}
  ~~=~ S(G_{i,j}) \,.
  \label{eq:notsohard}
  \end{eqnarray}
Recalling the expansion \eqref{eq:CorSN(T)v3}, it follows that
  \be
  \bearll
            \CorSN(\mathrm T)
  \!\!& \dsty = \sum_{i,j\in\I(\C)} \mathrm z(A)_{i,j} \, S(G^{}_{i,j})
  \Nxl2 & \dsty
  = \frac1{\DC^2} \sum_{i,i',j,j'\in\I(\C)} s^{}_{i',i^\vee}\, \mathrm z(A)_{i,j}
  \, s^{}_{j,j'} \, G^{}_{i',j'}
  = \sum_{i',j'\in\I(\C)} (S^{-1}\,\mathrm z(A)\, S)^{}_{i',j'} \, G^{}_{i',j'} 
  \eear
  \ee
with $S_{i,j} \,{:=}\, s_{i,j}/\DC$. Hence indeed we have
  \be
            S(\CorSN(\mathrm T)) = \CorSN(\mathrm T)
  ~\Longleftrightarrow~ [S,\mathrm z(A)] = 0 \,.
  \ee

Similarly, invariance of the correlator 
            $\CorSN(\mathrm T)$ 
under the modular
$T$-transformation is equivalent to $[T,\mathrm z(A)] \eq 0$ which, in turn, 
is equivalent to the statement that the object $Z(A) \iN \ZC$ has trivial twist.

\begin{rem}
A calculation similar to the one in \eqref{eq:notsohard} has been presented in 
\cite{hardi5}. The considerations in \cite{hardi5} are within the framework of
the tube category of the modular fusion category \C\ and of its category of
representations, which correspond to the cylinder category $\Cyloe$ and its
idempotent completion $\Cyle$, respectively.
\end{rem}


\subsection{Boundary Operator Product} \label{sec:boOPE}  %

Among the correlators involving only field insertions on the boundary of the 
world sheet, the most basic one is the correlator for three boundary insertions
on a disk, which describes the operator product of boundary insertions. The
world sheet for this correlator is a disk $\mathrm D$ with a single two-cell, 
labeled by a Frobenius algebra $A$, and with three physical boundary segments
labeled by $A$-modules $M_1$, $M_2$ and $M_3$. We denote this world sheet by
$\mathrm D_{M_1,M_2,M_3}$, and the corresponding \cws\ 
by $\widetilde{\mathrm D}_{M_1,M_2,M_3}$. In pictures,
  \be
  \mathrm D_{M_1,M_2,M_3} ~=~
  \raisebox{-9pt}{\scalebox{1.15}{\input{BP0.tikz}}}
  ~~~~~{\rm and}~~~~~ \widetilde{\mathrm D}_{M_1,M_2,M_3} ~=~
  \scalebox{1.05}{\input{BP1.tikz}}
  \label{eq:pic:bdyOPE1}
  \ee

Analogously as we did in the case of defect fields, we want to obtain the world 
sheet $\widetilde{\mathrm D}$ from a suitable string net on the disk $\mathrm D$
that encodes the algebraic 
information about the boundary operator product. We choose conventions such
that two of the boundary fields, say $\BB^{M_1,M_2}$ and $\BB^{M_2,M_3}$ are 
incoming, while the third is outgoing and is thus given by $\BB^{M_1,M_3}$.
Recall from \eqref{BB-MN} that these fields are internal Homs,
$\BB^{M_i,M_j} \eq \iHom_{\Mod A}(M_i,M_j)$. The natural candidate for 
their operator product is thus the canonical composition
  \be
  \imu(M_1,M_2,M_3) \in \HomC(\iHom(M_2,M_3) \oti \iHom(M_1,M_2) , \iHom(M_1,M_3))
  \ee
of internal Homs (see \eqref{eq:imu}). Accordingly we consider the string net
  \be
  \Gama_{\mathrm D;M_1,M_2,M_3} ~:=~~
  \scalebox{1.1}{\input{BP2.tikz}}
  ~~ \in~ \SN(\mathrm D,\bv_{M_1,M_2,M_3}) \,,
  \label{eq:pic:bdyOPE2}
  \ee
where the boundary value $\bv_{M_1,M_2,M_3}$ is a circle with three points 
labeled by $(\BB^{M_1,M_2})^\vee$, $(\BB^{M_2,M_3})^\vee$ and $\BB^{M_1,M_3}$.
Note that
  \be
  \bearll
  \SN(\mathrm D,\bv_{M_1,M_2,M_3}) \!\!& \cong \HomZ(\one,\LL(\BB^{M_1,M_3}_{}
  {\otimes}\, (\BB^{M_1,M_2})^\vee {\otimes}\, (\BB^{M_2,M_3})^\vee))
  \Nxl3 &
  \cong \HomC(\one,\BB^{M_1,M_3}_{} {\otimes}\, (\BB^{M_1,M_2})^\vee {\otimes}\,
  (\BB^{M_2,M_3})^\vee)
  \Nxl3 &
  \cong \HomC(\BB^{M_2,M_3} \,{\otimes}\, \BB^{M_1,M_2}, \BB^{M_1,M_3}) \,.
  \eear
  \ee

Now 
we can invoke \eqref{eq:iHom=vee} to
write $\BB^{M_i,M_j} \eq M_j^{} \otA M_i^\vee$, and we then
identify the boundary value $\bv_{M_1,M_2,M_3}$ with an isomorphic one 
so as to rewrite the string net \eqref{eq:pic:bdyOPE2} as
  \be
  \Gama_{\mathrm D;M_1,M_2,M_3} ~=~~
  \scalebox{1.0}{\input{BP3.tikz}}
  \ee
Note that the surface obtained this way is the one underlying the \cws\
$\widetilde{\mathrm D}_{M_1,M_2,M_3}$. We can thus read off that we indeed have
  \be
  \Gama_{\mathrm D;M_1,M_2,M_3} = \CorSN(\mathrm D_{M_1,M_2,M_3}) \,.
  \ee


\subsection{Bulk-boundary Operator Product} \label{sec:buboOPE} %

Among the correlators involving field insertions both in the bulk and on the 
boundary, the most basic one is the correlator for one bulk and one boundary 
insertion on a disk. This correlator encodes a connection between bulk and
boundary insertions, which is called the bulk-boundary operator product. Since
in our approach defect fields can be treated in much the same way as bulk 
fields, we consider here the more general situation of one defect and one boundary
insertion on a disk. The corresponding world sheet, which we denote by
$\mathrm D_{X,Y;M}$, is a disk $\mathrm D$ having
two two-cells labeled by Frobenius algebra $A$ and $B$, respectively,
two defect lines labeled by $A$-$B$-bimodules $X$ and $Y$, and a physical
boundary segment labeled by a right $A$-module $M$; it looks as follows:
  \be
  \mathrm D_{X,Y;M} ~=~~
  \scalebox{1.0}{\input{BB0.tikz}}
  \ee

The crucial algebraic datum describing the connection between bulk and boundary
in this situation is the component
  \be
  \hspace*{-1.5em} \bearl
  \dsty \jj_M^{} ~= \bigoplus_{m\in\I(\Mod A)}\!\! \dd_m ~
  \scalebox{1.2}{\input{BB1.tikz}}
  \\~\\[-0.5em] \hspace*{9.2em}
  \in \HomC(\mbox{$\int_{M'\in \Mod A}$}\! M' \,{\ota}\, Y \,{\otimes_{\!B}}\, X^\vee \,{\ota}\,
  M'{}^\vee\,, M \,{\ota}\, Y \,{\otimes_{\!B}}\, X^\vee {\ota}\, M^\vee)
  \eear
  \label{eq:jjM}
  \ee
at the object $ M \,{\ota} Y \,{\otimes_{\!B}}\, X^\vee \,{\ota}\, M^\vee
\eq \BB^{M\ota X,M\ota Y}_{} $ of the structure morphism $\jj$ of the end 
  \be
  \int_{M'\in \Mod A}\! M' \,{\ota}\, Y \,{\otimes_{\!B}}\, X^\vee {\ota}\, M'{}^\vee
  = U(\DD^{X,Y}) \,.
  \ee
Accordingly we consider the string net
  \be
  \Gama_{\mathrm D;X,Y;M} ~:=~~
  \scalebox{1.1}{\input{BB2.tikz}}
  ~~ \in \SN(S^1{\times}I,{(\bvcan_{\DD^{X,Y}})}^{\!\vee} {\boxtimes} \bvm_{X,Y;M}) \,.
  \label{eq:pic:bubdyOPE}
  \ee
Here the boundary value $\bvcan_{\DD^{X,Y}}$ is the one defined in 
\eqref{eq:def:bvcan}, while $\bvm_{X,Y;M}$ in the first place consists of 
four points on $S^1$ labeled by $M$, $Y$, $X^\vee$ and $M^\vee$, respectively,
together with an idempotent consisting of three relative 
tensor product idempotents combined, but is in an obvious manner isomorphic 
to the boundary value that consists of a single point on $S^1$
labeled by $\BB^{M\ota X,M\ota Y}_{}$ together with the
identity cylinder. Note that
  \be
  \bearll
  \SN(S^1{\times}I,{(\bvcan_{\DD^{X,Y}})}^{\!\vee} {\boxtimes} \bvm_{X,Y;M})
  \!\!& \cong \HomZ(\DD^{X,Y},\LL(\BB^{M\ota X,M\ota Y}_{}))
  \Nxl3 &
  \cong \HomC(U(\DD^{X,Y}),\BB^{M\ota X,M\ota Y}_{}) 
  \eear
  \ee
(where we use that, \C\ being semisimple, $\LL$ is also right adjoint to the
forgetful functor $U$).
Inserting the explicit form \eqref{eq:jjM} of $\jj_M^{}$ and the identities
\eqref{eq:proofshrink1} and \eqref{eq:proofshrink2}, we can write
  \be
  \begin{array}{r} \dsty
  \Gama_{\mathrm D;X,Y;M} ~= \!\!\sum_{i\in\I(\C) \atop m,n\in\I(\Mod A)}\!\!\!\!
  \frac{\dd_i\,\dd_m\,\dd_n}{\DC^2} ~~
  \scalebox{1.2}{\input{BB3.tikz}}
  \hspace*{9em}
  \\[-1.4em] \dsty
  = \sum_{m\in\I(\Mod A)} \frac{\dd_m}{\DC^2} ~~
  \scalebox{1.2}{\input{BB4.tikz}}
  \eear
  \ee
Now recall from Example \ref{exa:pXY} the string net $r_{X,Y}$ defined in
\eqref{eq:eXYrXY}, which furnishes an isomorphism
$\bvz_{X,Y} \Rarr\cong \bvcan_{\DD^{X,Y}}$. Precomposing with $r_{X,Y}$ yields
the string net
  \be
  \Gama_{\mathrm D;X,Y;M} \circ r_{X,Y} ~=~~
  \scalebox{1.1}{\input{BB5.tikz}}
  \ee
in $\SN(S^1{\times}I,{(\bvz_{X,Y})}^{\!\vee} {\boxtimes} \bvm_{X,Y;M})$.
By comparison with Definition \ref{def:CorSN}, we thus conclude that
  \be
  \Gama_{\mathrm D;X,Y;M} \circ r_{X,Y} = \CorSN(\mathrm D_{X,Y;M}) \,.
  \ee

In short, the component \eqref{eq:jjM} of the structure morphism of the end
$U(\DD^{X,Y})$ indeed reproduces the string-net correlator for the world sheet
$\mathrm D_{X,Y;M}$.
In particular, specializing to the case that $B \eq A$ and that 
$X \eq A \eq Y$ as bimodules, we obtain the world sheet describing the
bulk-boundary operator product.

  ~ %

\section{Internal Eckmann-Hilton Relation} \label{sec:brEHr} %
The exchange law for ordinary natural transformations turns out to have an
analogue for internal natural transformations. To obtain this analogue, which
we call the \emph{internal Eckmann-Hilton relation}, we introduce three braided
colored operads. The string-net correlator defined in Section \ref{sec:SNcor}
provides a morphism from the braided colored operad of world sheets to the one
of string nets. We compose this morphism with a morphism from the string-net operad to
a certain endomorphism operad. We can then derive the internal Eckmann-Hilton relation
by invoking compatibility of the so obtained composite morphism with operadic 
composition. The result illustrates the utility of string nets for understanding
algebra in braided tensor categories.

\subsection{An Internalized Eckmann-Hilton Argument}  %

Recall the diagram \eqref{eq:Poincare} which gives a Poincar\'e-dual view of 
(two-pronged) defect fields. In this section we consider multiple operator products
of defect fields for which the description in terms of \eqref{eq:Poincare} looks
as follows (for better readability we suppress the labels of the defect fields):
  \be
  \bearl ~\\[-1.8em]
  \begin{tikzpicture}
  \node[left] at (0,0) {$\M$};
  \node[right] at (1.96,0) {$\M'$};
  \node[right] at (4.89,0) {$\M''$};
  \draw[very thick,color=\colorDefect,->] (-0.2,0.23) .. controls (0.7,1.35) and (1.3,1.35) .. (2.2,0.21);
  \node[color=\colorDefect] at (1.02,1.39) {$X_5$};
  \draw[very thick,color=\colorDefect,->] (-0.08,0) -- (2.09,0);
  \node[color=\colorDefect] at (1.46,0.24) {$X_3$};
  \draw[very thick,color=\colorDefect,->] (-0.2,-0.26) .. controls (0.7,-1.35) and (1.3,-1.35) .. (2.18,-0.28);
  \node[color=\colorDefect] at (1.02,-1.44) {$X_1$};
  \draw[thick,double,->] (1,-0.96) --
       node[sloped,xshift=-3pt,yshift=6pt] {$\scriptstyle \DD$} (1,-0.14) ;
  \draw[thick,double,->] (1,0.11) --
       node[sloped,xshift=-3pt,yshift=6pt] {$\scriptstyle \DD$} (1,0.95) ;
  \begin{scope}[shift={(2.93,0)}]
  \draw[very thick,color=\colorDefect,->] (-0.2,0.23) .. controls (0.7,1.35) and (1.3,1.35) .. (2.2,0.21);
  \node[color=\colorDefect] at (1.02,1.39) {$X_6$};
  \draw[very thick,color=\colorDefect,->] (-0.16,0) -- (2.09,0);
  \node[color=\colorDefect] at (1.47,0.24) {$X_4$};
  \draw[very thick,color=\colorDefect,->] (-0.29,-0.26) .. controls (0.7,-1.35) and (1.3,-1.35) .. (2.18,-0.28);
  \node[color=\colorDefect] at (1.02,-1.44) {$X_2$};
  \draw[thick,double,->] (1,-0.96) --
       node[sloped,xshift=-3pt,yshift=6pt] {$\scriptstyle \DD$} (1,-0.14) ;
  \draw[thick,double,->] (1,0.11) --
       node[sloped,xshift=-3pt,yshift=6pt] {$\scriptstyle \DD$} (1,0.95) ;
  \end{scope}
  \end{tikzpicture}
  \\[-1.8em]~ \eear
  \ee
with $\M \eq \Mod A$, $\M' \eq \Mod A'$ and $\M'' \eq \Mod A''$ for simple
special symmetric Frobenius algebras $A,A',A'' \iN \C$.

Now recall that for ordinary natural transformations (and, more generally, for 
2-morphisms in a bicategory), the horizontal and vertical compositions satisfy the 
exchange law
  \be
  \mu_{\rm ver}^{} \circ (\mu_{\rm hor}^{} \oti \mu_{\rm hor}^{})
  = \mu_{\rm hor}^{} \circ (\mu_{\rm ver}^{} \oti \mu_{\rm ver}^{}) \,.
  \label{eq:EH}
  \ee
It should be
appreciated that this is a statement about elements of sets, and it implicitly relies
on the fact that the category $\Set$ of sets is a \emph{symmetric} monoidal category.

In contrast, internal natural transformations are objects of the category \ZC\
which is braided but (generically) not symmetric. However, we can show the 
following generalization to the 
    internalized    %
setting:

\begin{thm} \label{thm:EH}
Let \C\ be a modular fusion category and $A$, $A'$ and $A''$ be
simple special symmetric Frobenius algebras in \C. Then the diagram 
  \be
  \begin{tikzcd}[row sep=3.7em,column sep=-1.7em]
  ~ &
  \DD^{X_4,X_6}{\otimes}\,\DD^{X_2,X_4}{\otimes}\,\DD^{X_3,X_5}{\otimes}\,\DD^{X_1,X_3}~~~~
  \ar[xshift=-2.4em]{dl}[swap,xshift=10pt,yshift=-3pt]
  {\id{\otimes}\cB^{}_{\DD^{X_2,X_4};\DD^{X_3,X_5}}{\otimes}\id}
  \ar[xshift=1.3em]{dr}[xshift=-2pt]{\muver{\otimes}\muver} & ~
  \\
  \DD^{X_4,X_6} {\otimes}\, \DD^{X_3,X_5} {\otimes}\, \DD^{X_2,X_4} {\otimes}\, \DD^{X_1,X_3}
  \ar{d}[swap]{\muhore{\otimes}\muhore}
  & ~ & \DD^{X_2,X_6} {\otimes}\, \DD^{X_1,X_5} \ar{d}{\muhore}
  \\
  \DD^{X_3\otimes_{\!A'}\!X_4,X_5\otimes_{\!A'}\!X_6} {\otimes}\,
  \DD^{X_1\otimes_{\!A'}\!X_2,X_3\otimes_{\!A'}\!X_4} \ar{rr}{\muver}
  & ~ &
  \DD^{X_1\otimes_{\!A'}\!X_2,X_5\otimes_{\!A'}\!X_6}
  \end{tikzcd}
  \label{eq:EHdiagram}
  \ee
involving the horizontal and vertical compositions of internal natural transformations
and the half-braiding $\cB$ of 
            $\,\DD^{X_2,X_4}$
commutes for all $X_1,X_3,X_5 \iN A \MoD A'$ and all $X_2,X_4,X_6 \iN A' \MoD A''$.
\end{thm}

Since the equality \eqref{eq:EH} provides the basis for the Eckmann-Hilton 
argument, we refer to the commutativity of \eqref{eq:EHdiagram} as the
\emph{internal Eckmann-Hilton relation}.

\begin{rem}
(i) Theorem \ref{thm:EH} renders Proposition 14 in \cite{fuSc25} precise. In the 
latter, neither the relevant (left or right) horizontal composition nor the relevant 
(half-)braiding were specified.
 \\[2pt]
(ii) There is, of course, also a variant in which the right rather than left 
horizontal composition appears: the diagram
  \be
  \begin{tikzcd}[row sep=3.7em,column sep=-1.7em]
  ~ &
  \DD^{X_3,X_5}{\otimes}\,\DD^{X_1,X_3}{\otimes}\,\DD^{X_4,X_6}{\otimes}\,\DD^{X_2,X_4}~~~~
  \ar[xshift=-2.4em]{dl}[swap,xshift=10pt,yshift=-3pt]
  {\id{\otimes}\cB^{}_{\DD^{X_1,X_3};\DD^{X_4,X_6}}{\otimes}\id}
  \ar[xshift=1.3em]{dr}[xshift=-2pt]{\muver{\otimes}\muver} & ~
  \\
  \DD^{X_3,X_5} {\otimes}\, \DD^{X_4,X_6} {\otimes}\, \DD^{X_1,X_3} {\otimes}\, \DD^{X_2,X_4}
  \ar{d}[swap]{\muhorz{\otimes}\muhorz}
  & ~ & \DD^{X_1,X_5} {\otimes}\, \DD^{X_2,X_6} \ar{d}{\muhorz}
  \\
  \DD^{X_3\otimes_{\!A'}\!X_4,X_5\otimes_{\!A'}\!X_6} {\otimes}\,
  \DD^{X_1\otimes_{\!A'}\!X_2,X_3\otimes_{\!A'}\!X_4} \ar{rr}{\muver}
  & ~ &
  \DD^{X_1\otimes_{\!A'}\!X_2,X_5\otimes_{\!A'}\!X_6}
  \end{tikzcd}
  \ee
commutes as well. This follows directly by combining Theorem \ref{thm:EH} with the 
identity \eqref{eq:muhore=muhorz.cB} which will be established below.
\end{rem}


\subsection{Three Braided Colored Operads}

Preparing and giving the proof of Theorem \ref{thm:EH} will occupy most of the rest 
of this section.
A main ingredient are certain braided colored operads in $\Set$ \cite{YAuD}. With 
the help of string nets we will establish morphisms between these operads. The
operads in question are:
 \Itemize
 \item
The braided colored \emph{operad $\OWS$ of genus-$0$ world sheets} with labels from \C. 
 \\[2pt]
The colors of $\OWS$ are all possible boundary data of world sheets. The set of 
operations of $\OWS$ consists of all genus-0 world sheets, with the boundary data 
on their boundary circles regarded as in- and outputs, and taken up to isotopy; 
e.g.\ the binary products are given by
  \be
  \OWS\Big( \begin{array}c c\\\!\!a~\,b\!\! \eear \Big) ~=~ \mbox{\Huge\{}
  \scalebox{0.85}{\input{EH0.tikz}}
  \mbox{\Huge\}/}_{\dsty\!\!\rm isotopy} \,,
  \ee
i.e.\ by the set of all genus-0 world sheets with boundary datum 
$a^\vee{\times}\, b^\vee{\times}\,c$, up to isotopy.
The operadic composition on $\OWS$ is the sewing of world sheets, and the braid 
group action on $\OWS$ is obtained by identifying the braid group $\mathrm B_n$
as a subgroup of $\Map(\surf^0_{n+1})$, with $\surf^0_{n+1}$ a standard sphere 
with $n{+}1$ holes.
 \item
The braided colored \emph{operad $\OSN$ of string nets}. 
 \\[2pt]
The colors of $\OSN$ are the objects of the cylinder category $\Cyle$. The set
of operations consists of (the sets underlying) the string-net spaces on genus-0
surfaces with appropriate boundary values, e.g.
  \be
  \OSN\Big( \begin{array}c \mathbf C\\\!\!\mathbf A~\,\bv\!\! \eear \Big)
  = \SN\big(\surf^0_3, {\mathbf A}^{\!\!\vee\,} {\boxtimes}\,\bv^{\!\vee\,}
  {\boxtimes}\,\mathbf C \big) \,.
  \ee
The operadic composition on $\OSN$ is the gluing of string nets, and the braid
group action on $\OWS$ is again obtained by identifying $\mathrm B_n$ as a 
subgroup of $\Map(\surf^0_{n+1})$.
 \item
The braided colored \emph{endomorphism operad} $\OHZ$. 
 \\[2pt]
The colors of $\OHZ$ are the objects of the Drinfeld center \ZC.  The set
of operations consists of (the sets underlying) the morphisms in \ZC, e.g.
  \be
  \OHZ\Big( \begin{array}c Z\\\!\!X~\,Y\!\! \eear \Big)
  = \HomZ(X\oti Y,Z) \,.
  \ee
The operadic composition on $\OHZ$ is the composition of morphisms in \ZC.
We define the braid group action to be generated by
  \be
  \hspace*{-1.0em}\mbox{\Huge(} 
  \scalebox{0.65}{\input{VP2.tikz}}
  \xmapsto{~\beta~} 
  \scalebox{0.65}{\input{HP12.tikz}}
  \mbox{\Huge)} ~~\xmapsto{~~~~} ~~
  \mbox{\Huge(}\hspace*{-0.4em}
  \scalebox{0.8}{\input{EH1.tikz}}
  \hspace*{-0.3em}\xmapsto{~\beta~}\hspace*{-0.7em}
  \scalebox{0.8}{\input{EH2.tikz}}
  \hspace*{-1.0em}\mbox{\Huge)}
  \label{eq:braidOHZ}
  \ee
for $X_1,X_2,Y \iN \ZC$, with the overbraiding standing for the half-braiding
$\cB_{X_2;X_1}^{}$.
\end{itemize}

\begin{rem} \label{rem:lin}
For both $\OSN$ and $\OHZ$ there is an obvious linear version, in which the
operations are given by the vector spaces $\SN(-)$ and $\HomZ(-)$, respectively,
instead of by their underlying sets. We could indeed formulate the present
considerations in a linear setting, by linearizing also the world sheet operad
$\OWS$. The latter can be achieved by regarding a world sheet as colored by the
linear bicategory $\Fr$ of Frobenius algebras (which we will introduce in 
                  Chapter  %
\ref{sec:uCorr}) and treating those $\Fr$-colored world sheets
in the same way as \C-colored graphs are treated in the string-net construction.
\end{rem}


\subsection{Proof of the Relation}

Now we define morphisms
  \be
  \CorSN \Colon \OWS \rarr~ \OSN
  \label{eq:UWS-OSN}
  \ee
and
  \be
  \phiH \Colon \OSN \rarr~ \OHZ
  \label{eq:OSN-OHZ}
  \ee
of braided colored operads in $\Set$ by the following prescriptions:
 \Itemize
 \item
$\CorSN$ acts on colors by sending a boundary datum to the corresponding object in
$\Cyle$ (compare Example \ref{exa:bvo2pbd}). $\CorSN$ acts on operations by sending
a world sheet $\S$ to its string-net correlator $\CorSN(\S)$. Compatibility with the
operadic compositions and braid group equivariance of these prescriptions are evident.
 \item
The definition of $\phiH$ is slightly more involved. It depends on
two types of auxiliary data: for each genus-0 surface $\surf$ a marking 
$w \eq w(\surf)$ without cuts, and for each object $\bv \iN \Cyle$ an isomorphism 
$\psi \eq \psi(\bv)$ from $\bv$ to its ``canonical form'' $\bvcan_{\Phi(B)}$, i.e.\ the one that
is analogous to \eqref{eq:def:bvcan}, with $\Phi \colon \Cyle \Rarr\simeq \ZC$ 
the equivalence \eqref{Cyl=Z}.  With a fixed choice for these data, $\phiH$ acts on colors
as $\bv \,{\xmapsto{~~}}\, \Phi(\bv)$. Its action on operations is
by the inverses of the isomorphisms $\HomZ(...,...) \Rarr\cong \SN(\surf,\bv_w)$
that are analogous to \eqref{eq:e.-.rr} and \eqref{eq:varphiw2}. To give an example,
when choosing again the marking \eqref{eq:def:varphiver} on the pair of pants we have
  \be
  \bearll
  \phiH \Colon
  \scalebox{0.8}{\input{VP10.tikz}}
  & \xmapsto{~~\Psi~~} \sum_{m,n,q\in\I(\M)}\!\! \frac{\dd_m\,\dd_n}{\DC^4}~~
  \scalebox{0.8}{\input{VP9.tikz}}
  \\[-0.5em]
  & ~~~~=~~ 
  \scalebox{0.8}{\input{EH3.tikz}}
  \\~\\[-0.3em]
  & \xmapsto{~~\varphi_w^{-1}~}~ \muver ~\in \OHZ\Big( \begin{array}c
  \DD^{X_1,X_3} \\\!\! \DD^{X_2,X_3}~\, \DD^{X_1,X_2}\!\!\! \eear \Big) \,.
  \eear
  \ee
where $\Psi$ stands for pre- and post-composition with appropriate isomorphisms 
$\psi$ and $\varphi_w$ is given by \eqref{eq:e.-.rr}.
Similarly, with $\varphi_w$ as in \eqref{eq:varphiw2} we have
  \be
  \scalebox{0.8}{\input{HP9.tikz}}  
  ~ \xmapsto{~~\phiH~~} ~~ \muhore ~\in \OHZ\Big( \begin{array}c
  \!\! \DD^{X_1\otimes_{\!A'}\!X_2,X_3\otimes_{\!A'}\!X_4 } \!\!
  \\  \DD^{X_2,X_4}~~\, \DD^{X_1,X_3} \eear \Big) \,.
  \ee

Operadic composition results in a marking with cuts with an internal edge $e$ which
connects the outgoing circle of the inserted surface to the new root of the marking.
Compatibility with the operadic compositions is achieved by complementing the
prescription for $\phiH$ given above by the requirement to replace this marking by
the marking without cuts that is obtained by contracting the internal edge $e$.
(This is analogous to the F-move on markings, compare Figure 5.7 of 
	  \cite{BAki}.)
\end{itemize}

\begin{rem}
A different choice of the isomorphisms between the boundary values and their canonical
forms results in composing every element in a given morphism space with one and the
same isomorphism or its inverse. Similarly, a different choice of the markings results
in composing them with braiding and twist isomorphisms. The choices we make are
particularly convenient for revealing the
	internal  %
Eckmann-Hilton relation, but other choices will lead to that same relation as well.
\end{rem}

Let us verify that our prescription \eqref{eq:braidOHZ} leads to the correct
braid group action on $\OHZ$: Under the move $\beta$ that is shown on the 
left hand side of \eqref{eq:braidOHZ} we have
  \be
  \hspace*{-1.2em}
  \scalebox{0.8}{\input{EH4.tikz}}
  ~~ \xmapsto{~~~~} ~~ 
  \scalebox{0.8}{\input{EH5.tikz}}
  ~~ = ~~ 
  \scalebox{0.8}{\input{EH6.tikz}}
  \label{eq:checkcB}
  \ee
where the equality holds by the cloaking relation \Cite{Lemma\,3.7}{scYa}
  \be
  \scalebox{0.8}{\input{CL1.tikz}}
  ~~~ =
  \scalebox{0.8}{\input{CL2.tikz}}
  \ee
The overbraiding on the right hand side of \eqref{eq:checkcB} is the
half-braiding $\cB_{X_2;X_1}^{}$. Thus we indeed obtain the braid group action 
on $\OHZ$ as defined in \eqref{eq:braidOHZ}.
(A priori, in \eqref{eq:braidOHZ} one might have considered instead the
inverse half-braiding $\cB_{X_1;X_2}^{-1}$; the present calculation
shows that the choice made in \eqref{eq:braidOHZ} is the correct one.)

 \medskip

\noindent
\emph{Proof of Theorem {\rm \ref{thm:EH}}.}
Consider the composite
  \be
  \CorSNt := \phiH \circ \CorSN \Colon \OWS \rarr~ \OHZ
  \ee
of the morphisms \eqref{eq:UWS-OSN} and \eqref{eq:OSN-OHZ}. We have
  \be
  \S^w_{\rm ver} \stackrel{\eqref{eq:Sver}}= ~~
  \scalebox{0.8}{\input{EH7.tikz}}
  ~~ \xmapsto{~\CorSNt~} ~ \muver ~~
  \ee
and
  \be
  \Shore \stackrel{\eqref{eq:def:Shore}}= ~~
  \scalebox{0.8}{\input{EH8.tikz}}
  ~~ \xmapsto{~\CorSNt~} ~ \muhore \,.
  \ee
Invoking the compatibility of $\CorSNt$ with operadic composition, we further get
  \be
  \S_{\rm h;v,v} ~:=~~
  \scalebox{0.9}{\input{EH9.tikz}}
  ~~~ \xmapsto{~\CorSNt~} ~~
  \scalebox{1.1}{\input{EH10.tikz}}
  \label{eq:CorSNt4hvv}
  \ee
as well as
  \be
  \S_{\rm v;h,h} ~:=~~
  \scalebox{0.9}{\input{EH11.tikz}}
  ~~~ \xmapsto{~\CorSNt~} ~
  \scalebox{1.1}{\input{EH12.tikz}}
  \label{eq:CorSNt4vhh}
  \ee
where we use the short-hand $X_{i,j} \,{:=}\, X_i \,{\otimes_{\!A'}} X_j$.
Moreover, applying the braid group element $\beta_{2,3}$ to 
\eqref{eq:CorSNt4vhh} gives
  \be
  \beta_{2,3}(\S_{\rm v;h,h}) ~=~~ 
  \scalebox{0.9}{\input{EH9.tikz}}
  ~~~ \xmapsto{~\CorSNt~} ~
  \scalebox{1.1}{\input{EH13.tikz}}
  \label{eq:cBCorSNt4vhh}
  \ee
Notice that
  \be 
  \beta_{2,3}(\S_{\rm v;h,h}) = \S_{\rm h;v,v} \,.
  \ee
By comparison of \eqref{eq:CorSNt4hvv} and \eqref{eq:cBCorSNt4vhh} we arrive at the
desired equality that states the commutativity of the diagram \eqref{eq:EHdiagram}.
~ \hfill $\Box$

\begin{rem} \label{rem:muhore=cB.muhorz}
By similar arguments one sees that the equality $\Shore \eq \beta(\Shorz)$ obtained
in \eqref{eq:Shore=beta.Shorz} gets mapped under $\CorSNt$ to
  \be
  \muhore = \beta(\muhorz) = \muhorz \circ \cB^{}_{\DD^{X_2,X_4};\DD^{X_1,X_3}} ,
  \label{eq:muhore=muhorz.cB}
  \ee
showing that the left and right horizontal compositions are related by a half-braiding.
\end{rem}

\begin{rem}
Consider the special case that $A \eq A' \eq A''$ and that all defect fields 
involved are actually bulk fields $\DD^{A,A}$, for which vertical and horizontal
compositions coincide (see \eqref{eq:mu(Z(A))}). Then commutativity of the
diagram \eqref{eq:EHdiagram} amounts to the statement that the bulk algebra
$\DD^{A,A}$ in \ZC\ is braided commutative (compare Corollary 15 in \cite{fuSc25}).
\end{rem}

  ~ %

\section{Outlook: Universal Correlators} \label{sec:uCorr}  %

It can happen that two world sheets with different defect networks and thus
different \emph{defect patterns} --
that is, assignments of labels to the cells of the world sheet -- turn out to have the
same correlator. Existing approaches to correlators, including the string-net
construction considered in the previous chapters, treat such world sheets as
different even though they cannot be distinguished by their correlators.
Here we propose a variant of the string-net construction which
captures the additional relations by which world sheets with identical correlators
get identified. This novel description involves a new bicategorical string-net space 
and a linear map (introduced in the diagram \eqref{eq:def:UCor}) from that space to
ordinary string nets, to which we refer as the \emph{universal correlator}.
The proposal crucially uses the fact that the collections of
defect labels form a bicategory. This approach allows us in particular 
to sharpen the concept of a mapping class group action.

\subsection{String Nets Colored by a Pivotal Bicategory}  %

According to Definition \ref{def:S} the labels involved in the specification of
a world sheet $\S$ are simple special symmetric Frobenius algebras in the monoidal
category \C\ (for the two-cells of $\S$), bimodules between these (for the one-cells),
and bimodule morphisms (for the zero-cells). We call this assignment a \emph{defect
pattern}. A pertinent feature of a defect pattern is that these labels form the
objects, 1-morphisms, and 2-morphisms, respectively, of a linear pivotal bicategory 
(as defined e.g.\ in \Cite{Eq,\,(2.12)}{caRun3} or \Cite{Def.\,3.6}{camSc}). We denote
this bicategory by $\Fr$. So far we have ignored this feature of our construction; now
we put it into use. (A further motivation for considering the bicategory $\Fr$ is
the desire to linearize the operad $\OWS$ of world sheets, see Remark \ref{rem:lin}.)

As any pivotal bicategory, $\Fr$ comes with a graphical calculus on disks,
which relates disks having the same boundary value but different defect patterns
in the interior. 
As in the case of the ordinary string-net construction,
we promote this feature to local relations, so that we can identify
different defect patterns locally on disks embedded in a world sheet $\S$.

This observation motivates us to consider a variant of the string-net construction for 
which the coloring is not by a pivotal fusion category \C, but rather by some pivotal 
bicategory $\Bic$. Taking the quotient of the
free 
vector space spanned by all defect 
patterns with given boundary values by the local relations, we obtain a vector space
$\SNob(\surf,\bvb)$, which we call again the bare string-net space. A boundary value 
for this string-net construction is now an object $\bvb$ in the appropriate cylinder
category $\Cylob{\partial\surf}$ instead of $\Cylo{\partial\surf}$.
Let us illustrate the structure of $\Cylob{\partial\surf}$, restricting for brevity
to the case that $\partial\surf \eq S^1$: As an example, consider
  \be
  \scalebox{1.4}{\input{UC0.tikz}}
  \ee
where $a$, $b$ and $c$ are objects in $\Bic$ and $X \iN \HomB(a,b)$, $Y \iN \HomB(c,b)$
and $Z \iN \HomB(c,a)$. Here the symbol $\vee$ indicates the pivotal duality,
compare \Cite{Def.\,3.6}{camSc}. Morphisms in $\Cylob{S^1}$ are equivalence classes
of graphs on $\cyl$ modulo the local relations resulting from the graphical calculus
for $\Bic$. The picture 
  \be
  \scalebox{1.0}{\input{UC1.tikz}}
  \ee
shows an example of a morphism; here $\alpha$ and $\beta$ are 2-morphisms in $\Bic$. 
The mapping class group of $\surf$ acts on the bicategorical string-net spaces.


\subsection{Universal Correlators from String Nets Colored by $\Fr$}

In the present context -- the construction of CFT correlators for given chiral data,
captured by a pivotal fusion category \C\ --  the pivotal bicategory of our interest is 
the bicategory $\Fr$ of defect labels. 
In this case a world sheet $\S$ with underlying surface $\surf \eq \surf_\S$
and boundary value $\bvf \iN \Cylof{\partial\surf}$ provides us with a defect pattern
which determines a vector in the bare string-net space. We describe it as a linear map
  \be
  \hspace*{2.2em} \begin{array}{rl}
  \delta_\S \Colon
          \ko
  & \!\!\! \rarr~\, \SNof(\surf,\bvf) \,,
  \Nxl2
  1 & \!\!\! \xmapsto{~~~}\, [\S] \,,
  \eear
  \label{eq:def:deltaS}
  \ee
where $[\S]$ is the equivalence class
of the world sheet regarded as a $\Fr$-colored 2-graph on $\surf$.
On the other hand, we can view the vector in the string-net space based on $\calc$
that describes the string-net correlator $\CorSN(\S)$ as a linear map
  \be
  \begin{array}{rl}
  \CorSNb(\S) \Colon 
          \ko
  & \!\!\! \rarr~\, \SN(\surf,\FFt(\bvf)) \,,
  \Nxl2
  1 & \!\!\! \xmapsto{~~~}\, \CorSN(\S) \,.
  \eear
  \ee
Here 
  \be
  \FFt \Colon \Cylof{\partial\surf} \rarr~\, \Cyl{\partial\surf}
  \ee
is the obvious map between the two objects of the two categories; as
an illustration, we have
     (displaying on the right hand side the idempotents that define the objects)
  \be
  \scalebox{1.2}{\input{UC2.tikz}}
  ~~~\xmapsto{~~~~}~~~
  \scalebox{1.2}{\input{UC3.tikz}}
  \ee

With the help of $\FFt$ the field map, which was introduced in \eqref{eq:theFF} 
just as a map from the set of geometric boundary circles of the world sheet
to the class of objects in the Drinfeld center, can be promoted to a 
functor.  For brevity, we describe this functor only for the case that 
$\partial\surf \,{\cong}\, (S^1)^{\sqcup n}$ with all $n$
boundary components incoming. Then the field map becomes a functor
  \be
  \FF \Colon \Cylof{\partial\surf}\opp \rarr~ \ZC^{\boxtimes n} ,
  \ee
namely the composite
  \be
  \FF \Colon \Cylof{\partial\surf}\opp \rarr{\FFt\opp_{}} \Cyl{\partial\surf}\opp
  \rarr{{(\phi^{-1}_{})}_*} \Cyle^{\boxtimes n} 
  \rarr{\Phi^{\boxtimes n}} \ZC^{\boxtimes n} ,	  
  \ee
with $\phi\colon (S^1)^{\sqcup n} \To \partial\surf$ an orientation reversing
diffeomorphism describing the parametrization of the $n$ ingoing boundaries
and $\Phi$ the functor \eqref{Cyl=Z}.

 \medskip

Using the bare string nets based on the pivotal bicategory $\Fr$, we can enrich our
description of correlators with further information and sharpen the concept
of a mapping class group action. First, note that world sheets related by local
relations from $\Fr$ have the same boundary data and thus have the same space
$\Bl_\C$ of conformal blocks. Their correlators thus take value in the same vector
space. We claim that a much stronger statement holds: such world sheets that give
the same vector in the bicategorical string-net spaces even have the same correlator
$\CorSN$. Concretely, 
consider a world sheet $\S$ with underlying surface $\surf$, and let $\delta_\S$ be 
the linear map \eqref{eq:def:deltaS} which maps $1 \iN 
          \ko
$ 
to $[\S] \iN \SNof(\surf,\bvf)$ with $\bvf \iN \Cylof{\partial\surf}$. We claim that 
there exists a unique linear map $\UCorSN(\surf,\bvf)$ such that the triangle
  \be
  \begin{tikzcd}[row sep=2.9em]
          \ko
  \ar{r}{\delta_\S}
  \ar[xshift=-5pt,yshift=-2pt]{dr}[swap,xshift=5pt]{\CorSNb(\S)}
  & \SNof(\surf,\bvf) \ar[dashed]{d}{\UCorSN(\surf,\bvf)}
  \\[9pt]
  ~ & ~~~~~~~~~~\SN(\surf,\FFt(\bvf))
  \end{tikzcd}
  \label{eq:def:UCor}
  \ee
commutes for every world sheet $\S$ with underlying surface $\surf$ and boundary
data corresponding to $\bvf$. 
Thus, given any surface $\surf$ and boundary value $\bvf$, besides the pertinent
space of conformal blocks realized as the string-net space $\SN(\surf,\FFt(\bvf))$,
there is naturally associated another vector space $\SNof(\surf,\bvf)$, and this space
comes with a natural map to $\SN(\surf,\FFt(\bvf))$ through 
which the correlator for any world sheet with underlying surface $\surf$ factorizes.
The vector space $\SNof(\surf,\bvf)$ is by definition spanned by string nets. We 
therefore define $\UCorSN(\surf,\bvf)$ as the map that sends 
the vector in the bare string-net space that is determined by the defect pattern
of a world sheet $\S$ to the string-net correlator of $\S$.

In this way every correlator can be obtained by pre-composition with the map 
$\delta_\S$ from \eqref{eq:def:deltaS} which captures the information in the
defect pattern, i.e.\ the additional structure that a world sheet carries beyond its
underlying surface. Since pre-composition is nothing but pull-back along a function, the
situation may be regarded as an analogue -- linearized and at one categorical level 
lower -- of the description of a principal $G$-bundle via the universal bundle on the 
classifying space $BG$ of a group $G$. Accordingly we refer to the linear map
$\UCorSN(\surf,\bvf)$ as the \emph{universal correlator} for $(\surf,\bvf)$.

It is a major insight that it is the equivalence class $[\S]$ in the vector space 
$\SNof(\surf,\bvf)$, rather than $\S$ itself, that is amenable to the `observation'
of any aspects of correlators. It is thus appropriate to refer to the class 
$[\S] \iN \SNof(\surf,\bvf)$ as the `observable world sheet' that corresponds to a
given `classical' world sheet $\S$, or also as the corresponding 
\emph{quantum world sheet}. Note that for being able to introduce this concept,
allowing for the presence of defects is essential.

The crucial observation is now that the 
natural map $\UCorSN(\surf,\bvf)$ is actually well defined and that, moreover 
it intertwines the action of the mapping class group $\Map(\surf)$ on
$\SNof(\surf,\bvf)$ and $\SN(\surf,\FFt(\bvf))$.
The proof of these statements requires technical details which have not appeared in
our previous discussions. Some of them are in fact beyond the scope of the present
book, so that we refrain from presenting the proof here. Instead we content
ourselves to mention a few pertinent ingredients:
 \Itemize
 \item
There is a 2-functor $\UC\colon \Fr \Rarr~ \BC$ from the bicategory $\Fr$ of defect labels
to the delooping bicategory $\BC$ of the monoidal category \C, i.e.\ to \C\
viewed as a bicategory with a single object $\ast\,$.
$\UC$ maps every object of $\Fr$ to the single object of $\BC$
and every bimodule and bimodule morphism to the underlying object and morphism in \C,
respectively; schematically,
  \be
  \begin{tikzcd}
  A 
  \ar[bend right=60,xshift=3pt,start anchor={[yshift=0.8ex]},end anchor={[xshift=1.3pt,
      yshift=-0.6ex]}]{d}[swap,xshift=1pt]{X}[name=LL,xshift=-0.2ex]{}
  \ar[bend left=60,xshift=-2pt,start anchor={[yshift=0.8ex]},end anchor={[yshift=-0.6ex]}]
      {d}[xshift=-1pt]{Y}[name=RR,xshift=-0.65ex]{}
  \\
  B \ar[Rightarrow,to path=(LL) -- (RR) \tikztonodes]{}[yshift=1.2pt]{\alpha}
  \end{tikzcd}
  \quad\raisebox{-2pt}{$\xmapsto{~~\UC~~}$}\quad
  \begin{tikzcd}
  \ast
  \ar[bend right=60,xshift=3pt,start anchor={[yshift=0.8ex]},end anchor={[xshift=1.3pt,
      yshift=-0.6ex]}]{d}[swap,xshift=1pt]{\dot X}[name=LL,xshift=-0.2ex]{}
  \ar[bend left=60,xshift=-2pt,start anchor={[yshift=0.8ex]},end anchor={[yshift=-0.6ex]}]
      {d}[xshift=-1pt]{\dot Y}[name=RR,xshift=-0.65ex]{}
  \\
  \ast \ar[Rightarrow,to path=(LL) -- (RR) \tikztonodes]{}[yshift=1.2pt]{\dot\alpha}
  \end{tikzcd}
  \ee
 \item
The functor $\UC$ comes with canonical lax and oplax structures, with the help of
which it produces the correlator $\CorSN(\S)$ except for the full Frobenius
graphs on the two-cells of $\S$. The lax and oplax
structures are compatible in a way akin to the algebra and coalgebra structures of a
Frobenius algebra for which the product is left-inverse to the coproduct, and they fit
with the dualities in $\Fr$ and $\BC$. As a consequence, the change of color afforded
by $\UC$ with its lax and oplax structures is compatible with contraction along
edges and loops.
 \item
That the prescription is also compatible with the merging of vertices is ensured
by the presence of the full Frobenius graphs in the actual correlator $\CorSN(\S)$.
\end{itemize}

\medskip

We end this chapter with a few immediate observations about universal correlators:
 \Itemize
 \item
In the context of universal correlators there is a more comprehensive notion of the
mapping class group of a world sheet that supersedes the Definition \ref{def:Maps} of 
$\Maps$, namely as the stabilizer of the vector 
$[\S] \eq \delta_\S(1)$ 
under the $\Map(\surf)$-action on $\SNof(\surf,\bvf)$:
  \be
  \Mapr := \mathrm{Stab}_{\Map(\surf)}([\S]) \,.
  \label{eq:def:Mapr}
  \ee
As an illustration, 
let $X$ 
        be 
an invertible bimodule and consider the following world sheets $\S_1$ and
$\S_2$ which are identified by the local relations in $\SNof$: 
  \be
  \S_1 \,:=~\,
  \scalebox{1.05}{\input{UC4.tikz}}
            ~~~\sim~~~ \frac{\dim(B)}{\dim(X)} ~  %
  \scalebox{1.05}{\input{UC5.tikz}}
  ~\,=:\, \S_2 \,.
  \label{eq:elf}
  \ee
In this case, the mapping class group $\Map(\S_1)$ 
      of the world sheet $\S_1$   %
does not include Dehn twists along the boundary circles,
whereas    %
$\Map(\S_2)$, and hence also $\widehat{\Map(\S_1)} \,{\equiv}\, \widehat{\Map(\S_2)}$ does.

That the two world sheets in \eqref{eq:elf} have the same correlator is a consequence
of the invertibility of the bimodule $X$. This type of equality of correlators is also 
at the basis of the fact that invertible bimodules describe symmetries of
conformal field theories
\cite{ffrs3,ffrs5}.

    \item
More generally, the fact that the linear map $\UCorSN(\surf,\bvf)$ from 
$\SNof(\surf,\bvf)$ to $\SN(\surf,\FFt(\bvf))$
is well defined and intertwines the action of $\Map(\surf)$ implies that the
correlator $\CorSN(\S)$ is invariant under the action of the group $\Mapr$ as
defined in \eqref{eq:def:Mapr}.
 \item
That the collection $\{ \UCorSN(\surf,\bvf) \}_{(\surf,\bvf)}$ also intertwines
with the sewing of surfaces is equivalent to the statement that the string-net
construction of correlators is compatible with all sewing constraints.
\end{itemize}

\newpage

\appendix
\section{Appendix}

The purpose of this appendix is to provide the reader with 
pertinent notions from category theory. Readers who have not been much exposed
to such notions are advised to read through the appendix thoroughly,
while those with a working knowledge of braided fusion categories may have 
to resort to it for looking up notations and conventions. The concepts covered
include: fusion categories; pivotal, spherical and braided structures on them;
algebras, coalgebras and Frobenius algebras internal to fusion categories;
modules and bimodules over such algebras as well as module categories; the internal
Hom of an object in a module category over a fusion category \C, which is in 
particular a source of algebras in \C; the Drinfeld center \ZC\ of \C, and the
central monad and comonad, which allow one to construct adjoints of the forgetful 
functor from \ZC\ to \C\ explicitly.
In addition, Section \ref{app:G} surveys the graphical calculus for spherical
fusion categories, which is used freely throughout the paper, including the 
simplified graphical notation for morphisms that involve symmetric Frobenius algebras.

\subsection{Spherical Fusion Categories} \label{app:C} %

Let \ko\ be an algebraically closed field of characteristic zero.
A \emph{fusion category} over \ko\ is a finitely semisimple \ko-linear rigid 
monoidal category with simple monoidal unit. 

To fix notations, we briefly comment on the qualifications in this definition.
A \emph{monoidal} category is a category \C\ equipped with a tensor product functor,
which we denote by $\otimes\colon \C\Times\C \To \C$, and with a monoidal unit object,
which we write as $\one \iN \C$, and with associativity and unit constraints satisfying
pentagon and triangle identities. Without loss of generality we take the tensor product
to be strict, i.e.\ $\otimes$ is associative on the nose and 
$\one \oti C \eq C \eq C \oti \one$ for all $C\iN\C$.
\ko-\emph{linearity} of a monoidal category means that the morphism sets $\HomC(-,-)$ are 
\ko-vector spaces and that composition as well as the tensor product of morphisms is bilinear.

Being \emph{finitely semisimple} means that the isomorphism classes of simple objects of \C\ 
form a finite set and that every object is a finite direct sum of simple objects.
We select a set $\I(\C)$ of representatives for the isomorphism classes of simple objects,
in such a way that $\one \iN \I(\C)$.
In a \ko-linear finitely semisimple category the morphism spaces are finite-dimensional.
These finiteness properties are crucial. Keeping them while removing the semisimplicity
requirement generalizes fusion categories to the class of
\emph{finite tensor categories}, which in many respects behave much like fusion categories.
A fusion category is the same as a semisimple finite tensor category.
Since $\one$ is simple, we have $\HomC(\one,\one) \eq \ko\,\id_\one$, and we can 
canonically identify $\HomC(\one,\one)$ with \ko.
 
Finally, in a \emph{rigid} monoidal category every object $C \iN \C$ has a 
\emph{right dual object} $C^\vee$ accompanied by an \emph{evaluation} morphism
$\evr C \iN \HomC(C^\vee \oti C,\one)$ and by a \emph{coevaluation} morphism 
$\coevr C \iN \HomC(\one,C \oti C^\vee)$ that satisfy the snake identities 
  \be
  (\id_C \oti \evr C) \cir (\coevr C \oti \id_C) = \id_C \qquad\text{and}\qquad
  (\evr C \oti \id_{C^\vee}) \cir (\id_{C^\vee} \oti \coevr C) = \id_{C^\vee} \,,
  \label{eq:snake}
  \ee
as well as a \emph{left} dual object $\Vee C$ accompanied by corresponding evaluation and
coevaluation morphisms. The right and left dualities extend to functors $(-)^\vee$ and
$\Vee(-)$ from \C\ to its opposite category with opposite tensor product.

A \emph{pivotal structure} on a rigid monoidal category is a monoidal natural isomorphism
$\pi$ from the identity functor to the double (right, say) dual functor $(-)^{\vee\vee}$.
Again without loss of generality we take a pivotal structure (if it exists) to be
strict, i.e.\ take $\pi$ to be the identity natural transformation.
A category admitting a pivotal structure is called a pivotal category.
In a pivotal category \C\ one can associate to any endomorphism $f \iN \HomC(C,C)$
its right and left \emph{traces}, which in the strict case are given by
  \be
  \tr_\text{r}(f) = \evr {c^\vee} \circ (f \oti \id_{c^\vee}) \circ \coevr c
  \qquad \text{and} \qquad
  \tr_\text{l}(f) = \evr c  \circ (\id_{c^\vee} \oti f) \circ \coevr {c^\vee} \,,
  \ee
respectively. A pivotal structure is called \emph{spherical} if the right and left
traces of any endomorphism coincide.

A \emph{spherical fusion category} over \ko\ is a fusion category over \ko\ admitting,
and endowed with, a spherical pivotal structure.
In a spherical category one can identify left and right duals, and we tacitly do so.

We denote the (right and left) \emph{dimension} of an object $C$ of a spherical fusion 
category by $\dim(C) \eq \tr(\id_C)$; for simple objects $i \iN \I(\C)$ we abbreviate
  \be
  \dim(i) =: \dd_i \, \id_i\,.
  \ee
The \emph{global dimension} of a spherical fusion category \C\ is the number
  \be
  \DC^2 := \sum_{i \in \I(\C)} \dd_i^2 
  \label{eq:globaldim}
  \ee
(no choice of square root implied); this number is non-zero \Cite{Thm.\,2.3}{etno}.

One quantity in which the global dimension shows up is the
\emph{canonical color} (or \emph{Kirby color}, or \emph{surgery color}). This is by
definition the morphism
  \be
  \cancolor := \sum_{i\in\I(\C)} \frac{\dd_i}{\DC^2}\; \id_i
  ~\in \End_\C\big( \bigoplus\nolimits_{i\in\I(\C)} i \big) \,.
  \label{eq:cancolor}
  \ee


\subsection{Frobenius Algebras} \label{app:A} %

A unital associative algebra object internal to a monoidal category \C\ -- or an \emph{algebra} 
in \C, for short -- is a triple $A \,{\equiv}\, (A,\mu,\eta)$ consisting of an object $A \iN \C$,
a multiplication morphism $\mu \iN \HomC(A \oti A,A)$ that satisfies the associativity property
  \be
  \mu \circ (\id_A \oti \mu) = \mu \circ (\mu \oti \id_A) \,,
  \ee
and a unit morphism $\eta \iN \HomC(\one,A)$ obeying the unit properties
  \be
  \mu \circ (\eta \oti \id_A) = \id_A = \mu \circ (\id_A \oti \eta) \,.
  \ee
An algebra $A$ in \C\ is called \emph{haploid}, or \emph{connected},
iff $\dim_\ko \HomC(\one,A) \eq 1$.

Dually, a \emph{coalgebra} in \C\ is a triple $(C,\Delta,\eps)$ consisting of an object
$C \iN \C$ and comultiplication and counit morphisms $\Delta \iN \HomC(C,C \oti C)$ and
$\eps \iN \HomC(C,\one)$ obeying the coassociativity property
$(\id_C \oti \Delta) \cir \Delta \eq (\Delta \oti \id_C) \cir \Delta$
and $(\eps \oti \id_C) \cir \Delta \eq \id_C \eq (\id_C \oti \eps) \cir \Delta$.

A \emph{Frobenius algebra} in \C\ is a quintuple $F \,{\equiv}\, (F,\mu,\eta,\Delta,\eps)$ 
such that $(F,\mu,\eta)$ is an algebra in \C, $(F,\Delta,\eps)$ is a coalgebra in \C,
and such that in addition the compatibility conditions
  \be
  (\mu \oti \id_F) \circ (\id_F \oti \Delta) = \Delta \circ \mu
  = (\id_F \oti \mu)  \circ (\Delta \oti \id_F)
  \label{eq:Frobprop}
  \ee
between multiplication and comultiplication are satisfied. (The conditions
\eqref{eq:Frobprop} are not independent; when combined with the algebra and coalgebra
structures, each one implies the other.)

A Frobenius algebra $F$ in a spherical fusion category \C\ is called \emph{special} iff 
$\dim(F) \,{\ne}\, 0$ and the two equalities
  \be
  \mu \circ \Delta = \id_F \qquad\text{and}\qquad \eps \circ \eta = \dim(F) \, \id_\one
  \ee
are satisfied, and $F$ is called \emph{symmetric} iff the equality
  \be
  \big( (\eps\cir\mu) \oti \id_{F^\vee_{}} \big) \circ (\id_F \oti \coevr F)
  = \big( \id_{F^\vee_{}}  \oti (\eps\cir\mu) \big) \circ (\coevl F \oti \id_F)
  \label{eq:symmFrob}
  \ee
of morphisms from $F$ to $F^\vee \eq \Vee F$ holds. 

Owing to the Frobenius relation \eqref{eq:Frobprop}, the morphisms \eqref{eq:symmFrob} are 
actually \emph{iso}morphisms. Indeed, any Frobenius algebra in a rigid monoidal category is 
self-dual, i.e.\ is isomorphic to its right and to its left dual, and moreover the 
Frobenius-Schur indicator of $F$ is $+1$. Specifically, a symmetric Frobenius algebra $F$ 
in a spherical fusion category is self-dual in a particularly strong form: we can 
(and do) \emph{identify} $F$ with its (right and left) dual object
and can choose both the right and left evaluation and coevaluation morphisms of $F$ 
in such a way that they are expressible in terms of the structure morphisms of the
algebra-coalgebra $F$, according to
  \be
  \evr F = \eps \circ \mu = \evl F \qquad\text{and}\qquad 
  \coevr F = \Delta \circ \eta = \coevl F \,.
  \label{eq:ev-coev-F}
  \ee


\subsection{Modules and Bimodules} \label{app:M} %

A \emph{left module} $M \,{\equiv}\, (M,\rho)$ over an algebra $(A,\mu,\eta)$
in a monoidal category \C\ is a pair consisting of an object $M \iN \C$ and a
\emph{representation morphism} $\rho \iN \HomC(A \oti M,M)$ satisfying
  \be
  \rho \circ (\id_A \oti \rho) = \rho \circ (\mu \oti \id_M) \qquad\text{and}\qquad
  \rho \circ (\eta \oti \id_M) = \id_M \,.
  \ee
Similarly, a \emph{right module} $N \,{\equiv}\, (N,\ohr)$ over $A$ is a pair consisting 
of an object $N$ and a morphism $\ohr \iN \HomC(N \oti A,N)$ such that
$\ohr \cir (\ohr \oti \id_A) \eq \ohr \cir (\id_N \oti \mu)$ and
$\ohr \cir (\id_N \oti\eta) \eq \id_N$.
Further, given two algebras $(A,\mu,\eta)$ and $(A',\mu',\eta')$ in \C, an
$A$-$A'$-bimodule in \C\ is a triple $B \eq (B,\rho,\ohr)$ such that
        $(B,\rho)$ 
is a left $A$-module, $(B,\ohr)$ is a right $A'$-module and such that the left 
$A$-action and right $A'$-action commute, i.e.
  \be
  \rho \circ (\id_A \oti \ohr) = \ohr \circ (\rho \oti \id_{A'}) \,.
  \ee
Every algebra $A$ has a natural structure of an $A$-bimodule. An algebra is called
\emph{simple} iff it is simple as a bimodule over itself; every connected algebra is
simple.

Dual to the notion of a left module over an algebra is the one of a left \emph{comodule}
$W \,{\equiv}\, (W,\delta)$ over a coalgebra $C \,{\equiv}\, (C,\Delta,\eps)$: an object
$W \iN \C$ together with a morphism $\delta \iN \HomC(W,C\oti W)$ such that
  \be
  (\id_C \oti \delta) \circ \delta = (\Delta \oti \id_W) \circ \delta 
  \qquad\text{and}\qquad (\eps \oti \id_W) \circ \delta = \id_W \,.
  \ee
Right comodules are defined analogously.  
If $A$ is a special Frobenius algebra, then any right $A$-module $M \eq (M,\ohr)$ 
has a canonical structure of a right $A$-comodule, with the right $A$-coaction given by
  \be
  \reflectbox{$\delta$} := (\ohr \oti \id_A) \circ (\id_M \oti (\Delta\cir\eta) ) \,,
  \label{eq:def:deltaM}
  \ee
and analogously for left modules.

Given algebras $A$, $A'$ and $A''$ and an $A$-$A'$-bimodule $(B,\rho,\ohr)$ as well as an
$A'$-$A''$-bimodule $(B',\rho',\ohr')$, the \emph{tensor product over $A'$} of $B$ and $B'$ 
is the $A$-$A''$-bimodule $B \otAp B'$ that is defined by the universal property of 
coequalizing the two morphisms
  \be
  B \oti A' \oti B' \xrightrightarrows [~\id_B \otimes \rho'~] {~\ohrs \otimes \id_{B'}~}
  B \oti B' \,.
  \ee
In case $A'$ is a special Frobenius algebra, the tensor product bimodule $B \otAp B'$
is isomorphic to the image of the split idempotent
  \be
  P_{B\otap B'} := (\ohr \oti \rho') \circ (\id_B \oti (\Delta' \cir \eta') \oti \id_{B'}) \,.
  \label{eq:def:PMAN}
  \ee

For any object $X\iN\C$ the morphism space $\HomC(X,M\otA N)$ is isomorphic to
the subspace $\HomC^{(A)}(X,M\oti N) \,{\subseteq}\, \HomC(X,M\oti N)$ that consists of
those morphisms which are invariant under post-composition with the idempotent $P_{M \ota N}$
defined according to \eqref{eq:def:PMAN}. We use this isomorphism to tacitly identify 
  \be
  \HomC(X,M\otA N) = \HomC^{(A)}(X,M\oti N) \,.
  \label{eq:HomCota=HomAoti}
  \ee
Analogously we identify the space $\HomC(M\otA N,X)$ with the subspace $\HomC^{(A)}(M\oti N,X)$
of $\HomC(M\oti N,X)$ of morphisms invariant under pre-composition with $P_{M \ota N}$, 
and similarly for morphism spaces involving objects of the form 
$M \otimes_{\!A_1}\! B_1 \otimes_{\!A_2} \dots\, \otimes_{\!A_n}\! B_n 
\otimes_{\!A_{n+1}}\! N$ with $M$ a right $A_1$-module, $N$ a left $A_{n+1}$-module,
and, for $i \iN \{1,2,...\,,n\}$, $B_i$ an $A_i$-$A_{i+1}$-bimodule.

\medskip

Given an algebra $A$ in a monoidal category \C, the right $A$-modules constitute
the objects of the \emph{category of
right $A$-modules}, which we denote by
$\Mod A$. The morphisms in $\Mod A$ are those morphisms in \C\ which intertwine the
$A$-action, i.e.
  \be
  \Hom_{\Mod A}((M,\ohr),(M',\ohrp))
  = \{ f \iN \HomC(M,M') \mid \ohrp \cir (f \oti \id_A) \eq f \cir \ohr\, \} \,.
  \ee
The categories $A$-mod of left $A$-modules and $A$-mod-$A'$ of $A$-$A'$-bimodules
are defined analogously. We abbreviate
  \be
  \Hom_{\Mod A} \equiv \Hom_A \qquad\text{and}\qquad
  \Hom_{A\text{-mod-}A'} \equiv \Hom_{A|A'} \,.
  \ee

If \C\ is rigid and $X \eq (\dot X,\rho)$ is a left module over an algebra $A$
in \C, then the object $\dot X^\vee$ that is right dual to the underlying object 
$\dot X$ admits a natural structure of a right $A$-module, while the left dual
object $\Vee\dot X$ admits a natural structure of a right ${}^{\vee\vee\!}A$-module.
Similarly, for a right module $Y \eq (\dot Y,\ohr)$ over an algebra $B$ in \C,
$\Vee\dot Y$ and $\dot Y^\vee$ admit natural structures of a left $B$-module 
and a left $B^{\vee\vee}$-module, respectively.
In particular, if \C\ is strictly pivotal and $A$ and $B$ are symmetric Frobenius
algebras in \C, then for $X$ an $A$-$B$-bimodule, the dual $\dot X^\vee \eq \Vee\dot X$ 
admits a natural structure of a $B$-$A$-bimodule. We call this $B$-$A$-bimodule the
bimodule \emph{dual} to $X$ and denote it by $X^\vee$. The assignment 
$X \,{\mapsto}\, X^\vee$ is functorial. Strict pivotality also allows us to identify
  \be
  X^{\vee\vee} = X
  \ee
as $A$-$B$-bimodules.


\subsection{Module Categories} \label{app:Mcat} %

A (\emph{left}) \emph{module category} \M\ over a monoidal category \C\ is an abelian 
category equipped with an exact functor $\act\colon \C \Times \M \To \M$, called the
\emph{action} of \C\ on \M, and with natural
families of isomorphisms $(C\oti C') \act M \To C \act (C' \act M)$ and
$\one \act M \To M$ for $C,C' \iN \C$ and $M \iN \M$ that satisfy pentagon and 
triangle relations analogous to those for the associator and unitors of a non-strict
monoidal category. By invoking the equivalence $\calm \,{\simeq}\, \Mod A_\calm$ 
(see below) and strictifying the underlying monoidal category,
for our purposes we can (and do) take also these isomorphisms to be identities.
Right module categories as well as \C-$\mathcal D$-\emph{bimodule categories} over a 
pair of monoidal categories \C\ and $\mathcal D$ are defined analogously. 

A \emph{module functor} between two module categories \M\ and \N\ over \C\ is a
functor $F\colon \M \To \N$ together with a natural family of isomorphisms
$F(X \act M) \To X \act F(M)$ for $X \iN \C$ and $M \iN M$ that satisfy a pentagon
and a triangle identity that express compatibility with the module constraints of
\M\ and \N. Since we take the latter to be identities, we can take the module
structure morphism of a module functor to be the identity as well.

The \emph{direct sum} $\M_1 \,{\oplus}\, \M_2$ of two module categories $\M_1$ and $\M_2$ 
over \C\ is the direct sum of $\M_1$ and $\M_2$ as abelian categories
with \C-action (and analogously mixed associator and unitor) given by the sum 
of the \C-actions for $\M_1$ and $\M_2$. A module category is called
\emph{indecomposable} iff it is not equivalent to a non-trivial direct sum of module 
categories.

For a (left) module category \M\ over a monoidal category \C\ and any pair of objects
$M,M'\iN\calm$, the \emph{internal Hom} $\iHom_\calm(M,M')$ is an object in \C\
together with a natural family of isomorphisms 
  \be
  \HomC(C,\iHomM(M,M')) \xrightarrow{~\cong~} \HomM(C\act M,M')
  \ee
for $C\iN\calc$. In other words, for every $M'\in\M$ the functor $\iHom_\M(-,M')$ is
right adjoint to the action functor $\Act$. For a generic module category internal Homs
need not exist. But they do exist under suitable finiteness and exactness conditions, 
which are in particular met for finite module categories over finite tensor categories,
and in particular for all semisimple module categories over fusion categories. For any 
triple $M,M',M''$ of objects in \C\ there is an associative multiplication morphism
  \be
  \imu \equiv \imu{}_{M,M',M''}^{}\Colon
  \iHom(M',M'')\otimes \iHom(M,M') \xrightarrow{\phantom{~~}} \iHom(M,M'')
  \label{eq:imu}
  \ee
in \C. In particular, for any $M \iN \M$ the internal End $\iHom(M,M)$ comes with a
natural structure of an algebra in \C.

Let further $G$ be a right exact \C-module functor between \C-module categories
\M\ and \N. Then for any $M,M'\iN\M$, by taking the composite
  \be
  \begin{array}r
  \Hom_\C(C,\iHom_\M(M,M')) \rarr\cong \Hom_\M(C\act M,M') \hsp{13.6}
  \Nxl2
  \rarr{~G~} \Hom_\M(G(C\act M),G(M')) \rarr\cong \Hom_\M(C\act G(M),G(M')) \hsp{2.3}
  \Nxl2
  \rarr\cong \Hom_\C(C,\iHom_\N(G(M),G(M'))
  \eear
  \ee
for $C \iN \C$ one defines a natural transformation from $\C^{\rm opp}$ to the category
of vector spaces and thus a morphism
  \be
  \iG \Colon \iHom_\M(M,M') \rarr~ \iHom_\N(G(M),G(M'))
  \label{eq:def:iG}
  \ee
of internal Hom objects in \C.

\medskip

A \emph{pivotal} module category over a pivotal finite tensor category \C\ is a module 
category \M\ over \C\ such that for any $M,M' \iN \M$ there are functorial isomorphisms 
  \be
  \iHom(M,M')^\vee \xrightarrow{~\cong~} \iHom(M',M)
  \ee
between internal Homs and dual internal Homs which are compatible with the pivotal 
structure of \C.

In case that the finite tensor category \C\ is semisimple and thus a fusion category -- 
the situation studied in the main text -- an indecomposable pivotal module category over
\C\ is the same
as a \emph{semisimple} module category. We prefer to use the qualification `pivotal'
also in the semisimple case because various results for pivotal (and thus semisimple)
module categories over fusion categories extend to pivotal module categories over
finite tensor categories, see \cite{schaum,shimi20,fuSc25}.
For instance, if the module category \M\ is pivotal, then for any $M\iN\calm$ the algebra
$\iHom(M,M)$ in \C\ has a natural structure of a \emph{symmetric Frobenius} algebra.
If \C\ is semisimple and thus a fusion category, then the Frobenius algebra
$\iHom(M,M)$ is in addition special \Cite{Thm.\,6.6}{schaum}.

\medskip

By setting $C \act (M,\ohr) \,{:=}\, (C \oti M,\id_C \oti \ohr)$, the category $\Mod A$ of 
\emph{right} modules over an algebra $A$ in \C\ becomes a \emph{left} module category 
over \C. Conversely, under suitable finiteness conditions which are in particular 
fulfilled in the case of fusion 
categories, for any left module category \M\ over \C\ there exists an algebra $A_\calm$ in
\C\ such that \M\ and $\Mod A_\calm$ are equivalent as module categories over \C. This 
algebra is not unique; algebras $A$ and $A'$ such that $\Mod A$ and $\Mod A'$ are 
equivalent as module categories are called \emph{Morita equivalent}.
In fact, any such algebra is an internal Hom $\iHom(M,M)$ for some $M \iN \M$. 

If \M\ is indecomposable, then the Morita class contains a representative that is
connected. Thus in particular for every indecomposable pivotal module category over 
a pivotal fusion category \C\ there is a connected special symmetric Frobenius algebra
$A$ in \C\ such that $\Mod A$ is equivalent to \M\ as a module category.
Moreover, for any two such algebras $A$ and $A'$, every right exact \C-module functor
from $\M \,{\simeq}\, \Mod A$ to $\M' \,{\simeq}\, \Mod A'$ is isomorphic to the functor
  \be
  G^X := {-} \otA X \Colon \M \To \M'
  \label{eq:def:GB}
  \ee
of forming the tensor product, over $A$, with some $A$-$A'$-bimodule $X$.
Finally, in terms of the Frobenius algebra $A$ the internal Hom is given by
  \be
  \iHom_{\Mod A}(M,M') = M' \otA M^\vee \,.
  \label{eq:iHom=vee}
  \ee


\subsection{Drinfeld Center} \label{app:Z} %

A \emph{braiding} $\cb$ on a monoidal category \C\ is a natural isomorphism from the
tensor product functor to the opposite tensor product, i.e.\ a natural family of
isomorphisms $\cb_{C,C'}\colon C\oti C' \To C' \oti C$ that satisfy two 
\emph{hexagon identities}, which express compatibility with the tensor product and
for strict \C\ reduce to
  \be
  \begin{array}{rl}
  & \cb_{C,C'\otimes C''} = (\id_{C'} \oti \cb_{C,C''}) \circ (\cb_{C,C'} \oti \id_{C''})
  \Nxl2
  \text{and} \quad &
  \cb_{C\otimes C',C''} = (\cb_{C,C''} \oti \id_{C'}) \circ (\id_{C} \oti \cb_{C',C''}) 
  \eear
  \ee
for $C,C',C'' \iN \C$.
A \emph{half-braiding} $\cB$ on an object $C \iN \C$ is a natural family of isomorphisms
  \be
  \cB_X \Colon C \oti X \rarr{} X \oti C
  \ee
for $X \iN \C$ obeying the single hexagon identity
$\cB_{X\otimes X'} \eq (\id_X \oti \cB_{X'}) \cir (\cB_X \oti \id_{X'})$ for $X,X'\iN \C$.

The \emph{Drinfeld center} \ZC\ of a monoidal category \C\ is the category whose objects
are pairs $(C,\cB)$ consisting of an object $C \iN \C$ and a half-braiding on $C$,
and whose morphisms are those morphisms in \C\ which intertwine the half-braiding, i.e.\
$\HomZ((C,\cB),(C',\cB'))$ is the subset of those morphisms $f \iN \HomC(C,C')$ for which
$(f \oti \id_X) \cir \cB_X^{} \eq \cB'_X \cir (\id_X \oti f)$ for all $X\iN\calc$.

The half-braidings of its objects endow the Drinfeld center of any monoidal category
with a braiding. The Drinfeld center of a fusion category is again a fusion category,
and the Drinfeld center of a finite tensor category is again a finite tensor category.
The Drinfeld center of a spherical fusion category is naturally a modular fusion category.
If \M\ and \N\ are finite module categories over a finite tensor category \C, then
the finite category $\Funre_\C(\M,\N)$ of right exact module functors is a finite
module category over the Drinfeld center \ZC.
Indeed, given a functor $F\iN \Funre_\C(\M,\N)$ and an object 
$X \eq (\dot X,\cB) \iN \ZC$, the right exact functor $X\act F$ defined by
$(X\act F)(M) \,{:=}\, \dot X\act (F(M))$ is endowed with the structure of a
\C-module functor via the isomorphisms
  \be
  \bearll
  (X\act F)(C\act M) = \dot X\act F(C\act M) \!\!\!& \rarr{~\cong~} (\dot X \oti C)\act F(M)
  \Nxl2 &
  \xrightarrow[\cong]{\,\cB_C^{}\act F(M)}
  (C \oti \dot X)\act F(M) = C\act \big( (X\act F)(M) \big) 
  \eear
  \label{eq:RexZCmodule}
  \ee
for $C \iN \C$.

\medskip

Given a braiding $\cb$ on a monoidal category \C, the family $\cb\rev$ consisting of the
isomorphisms $(\cb\rev)_{C,C'}^{} \eq \cb_{C',C}^{-1}$ is a braiding on \C\ as well.
We denote the braided category with underlying monoidal category \C\ and braiding
$\cb\rev$ by $\C\rev$. A pivotal braided fusion category (and, more generally, a pivotal
braided finite tensor category) \C\ is called \emph{modular} iff a (distinguished)
braided monoidal functor 
  \be 
  \GC \Colon \C\rev \boti \C \xrightarrow~ \ZC
  \label{eq:def:GC}
  \ee 
from the Deligne product the categories 
$\C\rev$ and \C\ to the Drinfeld center of \C\ is an equivalence.


\subsection{The Central Monad and Comonad} \label{app:ZZ} %

A \emph{monad} on a category \C\ is an endofunctor $T\colon \C \To \C$ endowed with
the structure of an algebra in the monoidal category of endofunctors of \C. Analogously,
a \emph{comonad} on \C\ is an endofunctor endowed with a coalgebra structure. A (left)
\emph{module} over a monad $T$ on \C\ is a pair $(C,\rho)$ consisting of an object $C\iN \C$
and a representation morphism $\rho\colon T(C) \To C$ that is compatible in the standard 
way with the algebra structure of $T$.

On any finite tensor category there are canonically a monad $\Zm\,{\equiv}\,\Zm_\C$
and a comonad $\Zc\,{\equiv}\,\Zc_\C$, called the \emph{central monad} and \emph{central 
comonad}, respectively. On objects their underlying endofunctors are given by
  \be
  \Zm\Colon C \,\xmapsto{~~} \int^{X\in\calc}\! X^\vee \oti C\oti X
  \qquad \text{and} \qquad
  \Zc\Colon C \,\xmapsto{~~} \int_{X\in\calc} X \oti C \oti X^\vee ,
  \label{eq:def:Zm,Zc}
  \ee
respectively. Here $\int^X$ is a \emph{coend} and $\int_X$ an \emph{end},
i.e.\ $\Zm$ and $\Zc$ come equipped with families of structure morphisms
  \be
  \imath^\Zm_{C;Y} \Colon Y^\vee \oti C\oti Y \rarr{} \Zm(C)
  \qquad \text{and} \qquad
  \imath^\Zcs_{C;Y} \Colon \Zc(C) \rarr{} Y \oti C \oti Y^\vee
  \label{eq:imathZ}
  \ee
for $C,Y\iN\C$ that are dinatural in $Y$ and obey a universal property with respect to
dinaturality. 

Being defined through a universal property, the objects $\Zm(C)$ and
$\Zc(C)$ are (if they exist) unique up to unique isomorphism.
If \C\ is semisimple and thus a fusion category -- the case considered in the main text --
as well as pivotal, so that we can identify left and right duals, they can both be 
expressed as a finite direct sum 
  \be
  \Zm(C) = \Zc(C) = \bigoplus_{i \iN \I(\C)} i^\vee \oti C \oti i
  \ee
over the isomorphism classes of simple objects.

The forgetful functor $U$ from the Drinfeld center \ZC\ to \C\ that omits the half-braiding
is exact and thus,
owing to the finiteness of \C\ and \ZC,
has both a left adjoint $U\la$ and a right adjoint $U\ra$. The 
endofunctors $\Zm$ and $\Zc$ on \C\ are related to these by
  \be
  \Zm = U \circ U\la \qquad\text{and}\qquad \Zc = U \circ U\ra .
  \label{eq:Z-vs-UI}
  \ee
Moreover, the adjunctions $U\la \,{\vdash}\, U$ and $U \,{\dashv}\, U\ra$ 
are monadic and comonadic respectively, implying that
there are canonical equivalences between \ZC\ and the  categories of $\Zm$-modules and 
of $\Zc$-comodules. Indeed, the family of morphisms
  \be
  \partial^C_X := (\id_X \oti \imath^\Zm_{C;X}) \circ (\coevr X \oti \id_C \oti \id_X) \Colon
  C \oti X \rarr{} X \oti \Zm(C)
  \ee
for $C,X \iN \C$ generalizes the natural right coaction 
$\partial^\one_X\colon X \,{\to}\, X \oti \Zm(\one)$ of the coalgebra $\Zm(\one)$ on
$X \iN \C$; this family is called the \emph{universal coaction} associated with $\Zm$
\Cite{Ch.\,9.1.2}{TUvi}. Given any $\Zm$-module $(C,\rho)$, the family
  \be
  \cB^{(C,\rho)}_X := (\id_X \oti \rho) \circ \partial^C_X \Colon C \oti X \rarr{} X \oti C
  \label{eq:Zm-half-braiding}
  \ee
of morphisms in \C\
defines a half-braiding on $C$. This allows one to define a functor $\Zm\moD \To \ZC$ by
setting $(C,\rho) \,{\mapsto}\, (C,\cB^{(C,\rho)})$. One can show (see e.g.\ 
\Cite{Thm.\,9.3}{TUvi}) that 
this functor is a monoidal isomorphism. The case of $\Zc$-comodules can be treated dually.


\subsection{Coends for Generating Subcategories} %

\begin{defi} {}\Cite{Def.\,5.1.6}{KEly} \label{def:5.1.6} \\
Let $\cala$ be a finite \ko-linear abelian category. A subset $\mathrm U$ of objects of 
$\cala$ is said to \emph{p-generate} $\cala$ iff for any object $X \iN \cala$ there exists
an epimorphism $h_X\colon \mbox{\Large$\oplus$}_{i \in I\,} U_i \,{\twoheadrightarrow}\, X$
with a finite set $I$ and $U_i \iN \mathrm U$ for $i \iN I$, and for any $U \iN \mathrm U$
and any morphism $f\colon U \To X$ there exists a morphism 
$g\colon U \To \mbox{\Large$\oplus$}_{i \in I\,} U_i$ such that $f \eq h_X \cir g$.
We denote by $\mathcal U$ the full subcategory of $\cala$ having $\mathrm U$
as its set of objects.
\end{defi}

\begin{lem} {}{\rm \Cite{Thm.\,5.1.7}{KEly}} \label{lem:5.1.7} \\
Let $\cala$
        and 
$\calb$ be finite \ko-linear abelian categories and 
$\mathcal U \,{\subset} \cala$ be a full subcategory that p-generates $\cala$. 
Let $G\colon \cala\Times\cala\opp \To \calb$ be an exact functor and 
$G'\colon \mathcal U\Times\mathcal U\opp \To \calb$ its restriction to $\mathcal U$. 
Then the coends of $G$ and $G'$ exist, and the canonical morphism
  \be
  \int^{X\in\mathcal U}\! G'(X,\overline X) \rarr{~} \int^{Y\in\cala}\! G(Y,\overline Y)
  \ee
is an isomorphism.
\end{lem}

\begin{lem} \label{lem:Kar-pgen}
Let $\cald$ be a \ko-linear additive category. Then $\cald$ p-generates its
Karoubian envelope $\mathrm{Kar}(\cald)$.
\end{lem}

\begin{proof}
Given any object $X \eq (\dot X,p_X) \iN \mathrm{Kar}(\cald)$, with $\dot X \iN \cald$, we 
can set $h_X \,{:=}\, p_X \colon (\dot X,\id_{\dot X}) \,{\twoheadrightarrow}\, X$. Further,
for any object $Y \eq (\dot Y,\id_{\dot Y}) \iN \cald \,{\subseteq}\, \mathrm{Kar}(\cald)$
and any morphism $f \colon Y \To X$, we have $h_X \cir f \,{\equiv}\, p_X \cir f \eq f$.
Thus the conditions in Definition \ref{def:5.1.6} are satisfied by just taking 
$h_X \eq p_X$ and $g \eq f$.
\end{proof}


\subsection{Graphical Calculus for Spherical Fusion Categories} \label{app:G} %

A convenient tool for manipulating the string nets of our interest is the 
graphical calculus for morphisms in monoidal categories (see e.g.\ \Cite{Ch.\,2}{TUvi}),
including the treatment of algebras and their modules and bimodules (see e.g.\ 
\Cite{App.\,A}{fjfrs}). This calculus is tailored to the case that the monoidal category
is strict, albeit as long as one deals with equalities between morphisms it still makes
full sense for the general case, since all associativity and unit constraints involved
are easily reconstructed.
 
In the context of \emph{spherical fusion} categories, the graphical
calculus can be simplified in several ways: 
First, taking the pivotal structure to be strict we can unambiguously identify
  \be 
  \scalebox{1.0}{\input{GC0.tikz}} \,~~=~~
  \scalebox{1.0}{\input{GC1.tikz}}
  \ee
and can identify the left (co)evaluation morphism for $X$ with the right (co)evaluation 
morphism for $X^\vee \eq \Vee X$.

It also follows that
for any pair $X$ and $Y$ of objects we have a distinguished linear isomorphism 
$\zeta_{X,Y}\colon \HomC(\one,X \oti Y) \Rarr\cong \HomC(\one,Y \oti X)$ given by
  \be 
  \zeta_{X,Y}\Colon ~~
  \scalebox{1.0}{\input{GC2.tikz}} ~~~\longmapsto~~
  \scalebox{1.0}{\input{GC3.tikz}}
  \ee
Pivotality implies that $\zeta_{Y,X} \cir \zeta_{X,Y} \eq \id$. As an important 
consequence, we can think of a morphism
$\varphi \iN \HomC(\one,X_1 \,{\otimes} \cdots {\otimes}\, X_n)$ as depending only 
on the cyclic order of the tensor factors $X_1,X_2,...\,,X_n$. 
This can be made precise by associating to the cyclically ordered set 
$(X_1,X_2,...\,,X_n)$ not the individual vector space 
$\HomC(\one,X_1 \,{\otimes} \cdots {\otimes}\, X_n)$, but rather the limit of the diagram
  \be
  \bearl
  \cdots \Rarr{\zeta^{}_{X_{i+1}\otimes\cdots\otimes X_{i-1},X_{i}}}
  \HomC(\one,X_i \oti X_{i+1} \,{\otimes} \cdots {\otimes}\, X_{i-1})
  \Nxl2
  \phantom{\cdots} \Rarr{\zeta^{}_{X_i\otimes\cdots\otimes X_{i-2},X_{i-1}}}
  \HomC(\one,X_{i-1} \oti X_{i} \,{\otimes} \cdots {\otimes}\, X_{i-2})
  \rarr{~~~} \cdots \,,
  \eear
  \label{eq:lim-cyclic}
  \ee
but for ease of notation we represent this limit by any of the vector spaces appearing 
in the diagram.
To stress this cyclic aspect of morphisms, it is convenient to introduce the graphical 
notation \cite{kirI24}
  \be 
  \varphi ~~=~~ 
  \scalebox{1.2}{\input{GC4.tikz}}
  \label{eq:roundcoupon}
  \ee
in which the edges attached to a round coupon are to be thought of as being only cyclically
ordered.
Further, for suitable pairs of morphisms of this type there are partial composition maps
  \be
  \begin{array}{rr}
  \circ_X : ~~ & \HomC(\one,Y \oti X^\vee) \otimes_\ko \HomC(\one,X \oti Y')
  \rarr~ \HomC(\one,Y \oti Y') \,,
  \Nxl2
  & \varphi \otimes \varphi' \xmapsto{~~~} 
  (\id_Y \oti \evr X \oti \id_{Y'}) \circ (\varphi \oti \varphi') \,.
  \eear
  \ee
In particular for each $X\iN\C$ this provides a non-degenerate pairing
  \be
  \HomC(\one,X^\vee) \otimes_\ko \HomC(\one,X) \rarr~ \ko \,.
  \ee
Given a basis $\{\varphi^{\alpha}\}_{\alpha\in\Alpha}$ of $\HomC(\one,X)$ we denote the
basis of $\HomC(\one,X^\vee)$ that is dual to $\{\varphi^{\alpha}\}$ with respect to this
pairing by $\{\varphi_{\alpha}\}_{\alpha\in\Alpha}$. Following \cite{kirI24} we 
suppress the symbols for summing over a pair of such bases, i.e.\ use the
short-hand notation
  \be
  \sum_{\alpha\in\Alpha} \varphi^{\alpha} \oti \varphi_{\alpha} ~~=:~~
  \scalebox{1.2}{\input{GC5.tikz}}
  \label{eq_sum-alpha}
  \ee

In addition to the convention \eqref{eq:roundcoupon}, it is frequently convenient to
regard a morphism in $\HomC(\one,X)$ as an element of an isomorphic morphism space that
is related to $\HomC(\one,X)$ by rigidity, and also (in particular when product, coproduct
or representation morphisms are involved) to suppress the round coupon in the notation.
As an illustration, the completeness relation for $\HomC(\one,X)$ that follows from
the semisimplicity of \C, which we draw as
  \be
  \bigoplus_{i\in\I(\C)} \dd_i
  \scalebox{1.1}{\input{GC6.tikz}} ~~=~~~
  \scalebox{1.1}{\input{GC7.tikz}}
  \label{eq:dominance}
  \ee
is in fact a family of equalities in various different morphism spaces, rather than a
single equality: we may e.g.\ regard it as an equality in $\End_\C(X_1^{} \,{\otimes}
\cdots {\otimes}\, X_n^{})$, in which case the right hand side stands for the identity
morphism $\id_{X_1^{}\otimes\cdots\otimes X_n^{}}$, while when regarding it as an equality
in $\HomC(\one,X_1^{} \,{\otimes} \cdots {\otimes}\, X_n^{} \oti X_n^\vee {\otimes} \cdots
{\otimes}\, X_1^\vee)$, the right hand side stands instead for the coevaluation
$\coevr{X_1^{}\otimes\cdots\otimes X_n^{}}$.

As an example in which a coupon is omitted, we often think of the components of the 
half-braiding of an object $Y \eq (U(Y),\cB) \iN \ZC$ as morphisms in 
$\HomC(U(Y)\oti X,X\oti U(Y))$ and accordingly draw them as
  \be 
  \scalebox{1.1}{\input{GC8.tikz}} ~~=~~~
  \scalebox{1.1}{\input{GC9.tikz}} ~~=~~~
  \scalebox{1.1}{\input{GC10.tikz}}
  \label{pic:cB}
  \ee
Here in the last picture we slightly abuse notation by using the label $Y \iN \ZC$
in the description of a morphism in \C, which has the advantage that the label at the 
over-crossing cannot be mixed up with a braiding in \C\ and can thus be omitted. (For 
clarity, in addition we draw here strands labeled by objects of \C\
       and
strands labeled by objects of \ZC\ in two different shadings
(in the color version, in blue and green, respectively).
  
Similarly, instead of thinking of the product $\mu$ of an algebra $A$ as an element 
of the morphism space $\HomC(A\oti A,A)$, we may also regard it as an element of
$\HomC(\one,A\oti A^\vee,A^\vee)$, which amounts to the identification
  \be
  \scalebox{1.1}{\input{GC11.tikz}} ~~=~~
  \scalebox{1.1}{\input{GC12.tikz}}
  \ee
The defining relations for an algebra-coalgebra 
       in a spherical fusion category
to be Frobenius, special, and symmetric, respectively, are then drawn as
  \be
  \scalebox{1.1}{\input{GC13.tikz}} ~~=~~
  \scalebox{1.1}{\input{GC14.tikz}} ~~=~~
  \scalebox{1.1}{\input{GC15.tikz}}
  \ee
and as
  \be
  \scalebox{1.1}{\input{GC16.tikz}} ~~=~~
  \scalebox{1.1}{\input{GC17.tikz}} \qquad \text{and} \qquad
  \scalebox{1.1}{\input{GC18.tikz}} ~~=~~
  \scalebox{1.1}{\input{GC19.tikz}}
  \ee
respectively.

Finally, when dealing with a symmetric Frobenius algebra in a spherical fusion
category we 
may further simplify the graphical description of morphisms involving $A$ by using
that its evaluation and coevaluation morphisms can be taken to expressed through its 
structural morphisms as in \eqref{eq:ev-coev-F}. Specifically, we can identify $A$ with
its (left and right) dual and accordingly omit the orientation of an $A$-line, and we
can use \eqref{eq:ev-coev-F} to remove all occurrences of the unit and counit of
$F$. The resulting \emph{simplified graphs} have unoriented $A$-lines and do not involve
univalent coupons for the unit or counit. 
Let us give two examples of such simplified graphs:
First, the $A$-coaction defined in \eqref{eq:def:deltaM} can be simplified according to
  \be
  \reflectbox{$\delta$} ~~=~~
  \scalebox{1.1}{\input{GC20.tikz}} ~=:~~ 
  \scalebox{1.1}{\input{GC21.tikz}}
  \label{eq:pic:atled}
  \ee
Second, the idempotent $P_{B\otap B'}$ defined in \eqref{eq:def:PMAN} simplifies as
  \be
  P_{B\otap B'} ~~=~~  
  \scalebox{1.1}{\input{GC22.tikz}} ~=:~~
  \scalebox{1.1}{\input{GC23.tikz}}
  \label{eq:pic:PMAN}
  \ee
  
As one specific application of the simplified graphical calculus we note that
any \emph{Frobenius move} between two (simplified) \ffg s in the sense
of Definition \ref{defi:ffgraph} can be presented as a sequence of elementary
moves of the following form:
 \Itemize
 \item
The a-move
  \be
  \scalebox{1.1}{\input{GC24.tikz}} ~~\longleftrightarrow~~
  \scalebox{1.1}{\input{GC25.tikz}}
  \label{eq:def:a-move}
  \ee
 \item
The b-move
  \be
  \scalebox{1.1}{\input{GC26.tikz}} ~~\longleftrightarrow~~
  \scalebox{1.1}{\input{GC27.tikz}}
  \hspace{2.0em}
  \label{eq:def:b-move}
  \ee
 \item
The r-move
  \be
  \hspace{0.5em}
  \scalebox{1.1}{\input{GC28.tikz}} ~\longleftrightarrow~~~
  \scalebox{1.1}{\input{GC29.tikz}}
  \label{eq:def:r-move}
  \ee
\end{itemize}
Note that each of these elementary moves between simplified graphs stands for
a whole collection of moves between non-simplified graphs. For instance, the moves
  \be
  \begin{array}{c}
  \scalebox{1.1}{\input{GC16.tikz}} ~~\longleftrightarrow~
  \scalebox{1.1}{\input{GC17.tikz}} ~\qquad \text{and} \qquad~~
  \scalebox{1.1}{\input{GC30.tikz}} \longleftrightarrow~
  \scalebox{1.1}{\input{GC17.tikz}}
  \\~\\[7pt]
  ~~\text{and} \qquad~
  \raisebox{1.8em}{\scalebox{1.1}{\input{GC31.tikz}}} \longleftrightarrow~
  \scalebox{1.1}{\input{GC32.tikz}}
  \eear
  \ee
~\\[-1.9em]
constitute
three possible realizations of the b-move as a move between non-simplified graphs.


\subsection{Three Useful Lemmas} \label{app:appl} %

In this appendix we provide the proofs of observations that are used in
the main text. We freely use the graphical calculus.

\medskip

The following result is an important ingredient in the definition of the field map in
Section \ref{sec:fields}, which allows us to treat boundary fields on the same footing
as bulk and defect fields:

\begin{lem} \label{lem:Lxy=Lyx}
Let \C\ be a finite tensor category and $C$ and $D$ be any two objects in \C. Then
the objects $\LL(C\oti D)$ and $\LL(D\oti C)$ in \ZC\ are isomorphic.
\end{lem}

\begin{proof}
Consider morphisms $f\colon \LL(C\Oti D) \Rarr{} \LL(D\Oti C)$ and
$g\colon \LL(D\Oti C) \Rarr{} \LL(C\Oti D)$ in \ZC\ that are defined, by means
of the dinatural structure morphisms \eqref{eq:imathZ}, by the equalities
  \be
  \scalebox{1.1}{\input{LM7.tikz}}
  ~~=~~~
  \scalebox{1.1}{\input{LM8.tikz}}
  \ee
and
  \be
  \scalebox{1.1}{\input{LM9.tikz}}
  ~~=~~~
  \scalebox{1.1}{\input{LM10.tikz}}
  \ee
respectively. The composite $g \cir f$ satisfies
  \be
  \scalebox{1.1}{\input{LM11.tikz}}
  ~~=~~~
  \scalebox{1.1}{\input{LM12.tikz}}
  ~~~=~~~
  \scalebox{1.1}{\input{LM13.tikz}}
  \ee
where the second equality follows by dinaturality. Together with the snake identity
\eqref{eq:snake} this shows that $g \cir f \eq \id_{\LL(C\oti D)}$. Similarly one sees
that $f\cir g \eq \id_{\LL(D\oti C)}$. Moreover, it is readily checked that all 
morphisms involved are compatible with half-braidings, so that they are actually 
morphisms in \ZC. It thus follows that $f$ is a (non-canonical) isomorphism
$\LL(C\oti D) \Rarr{\cong} \LL(D\oti C)$ in \ZC.
\end{proof}

\begin{lem} \label{lem:char-eps}
For $A$ a simple special symmetric Frobenius algebra in a
	  spherical  %
fusion category \C\ and $M$ a right $A$-module the equality
  \be
  \scalebox{1.1}{\input{LM17.tikz}}
  ~~=~~ \frac{\dim(M)}{\dim(A)}~
  \raisebox{0.9em}{\scalebox{1.1}{\input{LM15.tikz}}}
  \label{eq:dm/dA-eps}
  \ee
holds, where the morphism on the right hand side is the counit $\eps$ of $A$.
\end{lem}

\begin{proof}
Consider the endomorphism
  \be
  f_M ~:= ~~
  \scalebox{1.1}{\input{LM16.tikz}}
  \ee
of $A$, 
obtained by combining \eqref{eq:dm/dA-eps} with the coproduct of $A$.
Since $A$ is symmetric Frobenius, $f_M$
is not just a morphism in \C, but even a morphism of $A$-bimodules (with 
respect to the regular $A$-bimodule structure on $A$). Since $A$ is simple, this
implies that $f_M$ is a multiple $\xi_M\,\id_A$ of the identity morphism. 
Post-composing with the counit then shows that the morphism on the left hand side
of \eqref{eq:dm/dA-eps} equals $\xi_M\,\eps$. Further pre-composing with the
unit gives $\dim(M) \eq \xi_M\, \eps\cir\eta$. Since $A$ is special, this implies
\eqref{eq:dm/dA-eps}.
\end{proof}

\begin{lem} \label{lem:shrink0}
For $A$ a 
	  special    %
symmetric Frobenius algebra in a 
	  spherical  %
fusion category \C\ we have
  \be
  \sum_{m\in\I(\Mod A)}\! \dd_m^2 = \dim(A) \, \DC^2 \,,
  \ee
with $\DC^2$ the global dimension of \C.
\end{lem}

\begin{proof}
We have
  \be
  \hspace*{-0.8em}
  \bearll
  ~\\[-2.5em] 
  \dim(A) \, \DC^2 & \dsty
  \stackrel{\eqref{eq:globaldim}}=~ \sum_{i\in\I(\C)} \dd_i ~
  \scalebox{1.0}{\input{LM0.tikz}}
  ~~~= \!\sum_{i\in\I(\C)\atop m\in\I(\Mod A)}\!\!\! \dd_i\, \dd_m ~
  \scalebox{1.0}{\input{LM1.tikz}}
  \\~\\[-0.4em] & \dsty
  = \sum_{m\in\I(\Mod A)}\! \dd_m ~
  \scalebox{1.0}{\input{LM2.tikz}}
  ~~= \sum_{m\in\I(\Mod A)}\!\! \dd_m^2 \,.  
  \eear
  \label{eq:proofshrink}
  \ee
Here in the second equality we use the identity
  \be
  \bigoplus_{m\in\I(A\MoD B)}\! \dd_m ~~
  \scalebox{1.0}{\input{LM3.tikz}}
  ~~=~~
  \scalebox{1.0}{\input{LM4.tikz}}
  \label{eq:proofshrink1}
  \ee
valid for any right $A$-$B$-bimodule $X$, where the $\alpha$-summation is over a 
basis of $\Hom_{A\MoD B}(m,X)$; this is a variant of the completeness relation
\eqref{eq:dominance}. The third equality of \eqref{eq:proofshrink} follows from
the identity \Cite{Lemma\,4.3}{fjfrs}
  \be
  \bigoplus_{i\in\I(\C)} \dd_i 
  \scalebox{1.0}{\input{LM5.tikz}}
  ~~=~~~
  \scalebox{1.0}{\input{LM6.tikz}}
  \label{eq:proofshrink2}
  \ee
which holds for any right $A$-module $M$ and left $A$-module $N$ (with the 
$\alpha$-summation being over a basis of 
$\HomC^{(A)}(i,M\oti N) \,{\cong}\, \HomC(X,M\otA N)$), 
\end{proof}

Combining Lemma \ref{lem:char-eps} and Lemma \ref{lem:shrink0} we obtain

\begin{cor} \label{cor:shrink}
For $A$ a simple special symmetric Frobenius algebra in a 
	  spherical  %
fusion ca\-te\-go\-ry the equality 
  \be
  \sum_{m\in\I(\Mod A)} \frac{\dd_m}{\DC^2} ~
  \scalebox{1.1}{\input{LM14.tikz}}
  ~~=~~
  \scalebox{1.1}{\input{LM15.tikz}}
  \ee
holds.
\end{cor}

\newpage

\newcommand\wb{\,\linebreak[0]} \def\wB {$\,$\wb}
\newcommand\Arep[2]  {{\em #2}, available at {\tt #1}}
\newcommand\Bi[2]    {\bibitem[#2]{#1}}
\newcommand\inBO[9]  {{\em #9}, in:\ {\em #1}, {#2}\ ({#3}, {#4} {#5}), p.\ {#6--#7} {\tt [#8]}}
\newcommand\J[7]     {{\em #7}, {#1} {#2} ({#3}) {#4--#5} {{\tt [#6]}}}
\newcommand\JJ[4]    {{\em #4}, {#1} ({#2}) {{\tt [#3]}}}
\newcommand\JO[6]    {{\em #6}, {#1} {#2} ({#3}) {#4--#5} }
\newcommand\JP[7]    {{\em #7}, {#1} ({#3}) {{\tt [#6]}}}
\newcommand\BOOK[4]  {{\em #1\/} ({#2}, {#3} {#4})}
\newcommand\PhD[2]   {{\em #2}, Ph.D.\ thesis #1}
\newcommand\Prep[2]  {{\em #2}, preprint {\tt #1}}
\newcommand\uPrep[2] {{\em #2}, unpublished preprint {\tt #1}}
\def\adma  {Adv.\wb Math.}
\def\anma  {Ann.\wb Math.}
\def\anop  {Ann.\wb Phys.}
\def\apcs  {Applied\wB Cate\-go\-rical\wB Struc\-tures}
\def\atmp  {Adv.\wb Theor.\wb Math.\wb Phys.}   
\def\cocm  {Com\-mun.\wb Con\-temp.\wb Math.}
\def\comp  {Com\-mun.\wb Math.\wb Phys.}
\def\coma  {Con\-temp.\wb Math.}
\def\entc  {Electronic Notes in Theoretical Computer Science}
\def\foph  {Fortschr.\wb Phys.}
\def\geat  {Geom.\wB and\wB Topol.}
\def\imrn  {Int.\wb Math.\wb Res.\wb Notices}
\def\jims  {J.\wb Indian\wb Math.\wb Soc.}
\def\jams  {J.\wb Amer.\wb Math.\wb Soc.}
\def\joal  {J.\wB Al\-ge\-bra}
\def\joms  {J.\wb Math.\wb Sci.}
\def\jopa  {J.\wb Phys.\ A}
\def\jktr  {J.\wB Knot\wB Theory\wB and\wB its\wB Ramif.}
\def\mama  {ma\-nu\-scripta\wB mathematica\wb}
\def\mams  {Memoirs\wB Amer.\wb Math.\wb Soc.}
\def\momj  {Mos\-cow\wB Math.\wb J.}
\def\nupb  {Nucl.\wb Phys.\ B}
\def\pams  {Proc.\wb Amer.\wb Math.\wb Soc.}
\def\phrl  {Phys.\wb Rev.\wb Lett.}
\def\pnas  {Proc.\wb Natl.\wb Acad.\wb Sci.\wb USA}
\def\pspm  {Proc.\wb Symp.\wB Pure\wB Math.}
\def\qjmo  {Quart.\wb J.\wb Math.\wB Oxford}
\def\quto  {Quantum Topology}
\def\sema  {Selecta\wB Mathematica}
\def\sigm  {SIGMA}
\def\slnm  {Sprin\-ger\wB Lecture\wB Notes\wB in\wB Mathematics}
\def\taac  {Theo\-ry\wB and\wB Appl.\wb Cat.}
\def\tams  {Trans.\wb Amer.\wb Math.\wb Soc.}
\def\toap  {Topology\wB Applic.}
\def\topo  {Topology}
\def\trgr  {Trans\-form.\wB Groups}

 \end{document}

%% file: sample.tikzstyles
\tikzstyle{nc}=[fill=white, draw=black, shape=circle, minimum size=0.40cm]
\tikzstyle{nr}=[fill=white, draw=black, shape=rectangle, minimum width=0.60cm, minimum height=0.60cm]
\tikzstyle{nrg}=[fill={rgb,255: red,199; green,199; blue,199}, draw=black, shape=rectangle]
\tikzstyle{sr}=[fill=white, draw=black, shape=rectangle, minimum width=0.60cm, minimum height=0.40cm, inner sep=0pt]
\tikzstyle{lsr}=[fill=white, draw=black, shape=rectangle, minimum width=0.60cm, minimum height=0.40cm, inner sep=0pt, rotate=30]
\tikzstyle{rsr}=[fill=white, draw=black, shape=rectangle, minimum width=0.60cm, minimum height=0.40cm, inner sep=0pt, rotate=-30]
\tikzstyle{rrrsr}=[fill=white, draw=black, shape=rectangle, minimum width=0.60cm, minimum height=0.40cm, inner sep=0pt, rotate=-120]
\tikzstyle{ivsr}=[fill=none, draw=none, shape=rectangle, minimum width=0.60cm, minimum height=0.40cm, inner sep=0pt]
\tikzstyle{ivlsr}=[fill=none, draw=none, shape=rectangle, minimum width=0.60cm, minimum height=0.40cm, inner sep=0pt, rotate=30]
\tikzstyle{ivrsr}=[fill=none, draw=none, shape=rectangle, minimum width=0.60cm, minimum height=0.40cm, inner sep=0pt, rotate=-30]
\tikzstyle{ivrrrsr}=[fill=none, draw=none, shape=rectangle, minimum width=0.60cm, minimum height=0.40cm, inner sep=0pt, rotate=-120]
\tikzstyle{mr}=[fill=white, draw=black, shape=rectangle, minimum width=1.00cm, minimum height=0.50cm, inner sep=0pt]
\tikzstyle{product}=[fill={rgb,255: red,255; green,5; blue,80}, draw=black, shape=circle, inner sep=0pt, minimum size=0.15cm]
\tikzstyle{unit}=[fill={rgb,255: red,255; green,166; blue,217}, draw=black, shape=circle, inner sep=0pt, minimum size=0.15cm]
\tikzstyle{coproduct}=[fill={rgb,255: red,2; green,145; blue,255}, draw=black, shape=circle, inner sep=0pt, minimum size=0.15cm]
\tikzstyle{counit}=[fill={rgb,255: red,165; green,219; blue,255}, draw=black, shape=circle, inner sep=0pt, minimum size=0.15cm]
\tikzstyle{antipode}=[fill=white, draw=black, shape=circle, inner sep=0pt, minimum size=0.15cm]
\tikzstyle{region}=[fill={rgb,255: red,191; green,191; blue,191}, draw={rgb,255: red,191; green,191; blue,191}, shape=circle, minimum size=0.80cm]
\tikzstyle{boundarydisc}=[fill={rgb,255: red,128; green,128; blue,128}, draw=black, shape=circle, minimum size=0.80cm]
\tikzstyle{identity}=[fill=black, draw=black, shape=circle, inner sep=0pt, minimum size=0.15cm]
\tikzstyle{S}=[fill=white, draw=black, shape=rectangle, minimum width=0.15cm, minimum height=0.15cm]
\tikzstyle{S-}=[fill={rgb,255: red,199; green,199; blue,199}, draw=black, shape=rectangle, minimum width=0.15cm, minimum height=0.15cm]
\tikzstyle{dot}=[fill=white, draw=black, shape=circle, inner sep=0pt, minimum size=0.07cm]
\tikzstyle{bdot}=[fill=black, draw=black, shape=circle, inner sep=0pt, minimum size=0.07cm]
\tikzstyle{m}=[fill={rgb,255: red,255; green,5; blue,80}, draw={rgb,255: red,255; green,5; blue,80}, shape=circle, inner sep=0pt, minimum size=0.14cm, tikzit shape=circle]
\tikzstyle{mod}=[fill={rgb,255: red,172; green,229; blue,232}, draw={rgb,255: red,14; green,115; blue,177}, shape=circle]
\tikzstyle{smod}=[fill={rgb,255: red,172; green,229; blue,232}, draw={rgb,255: red,14; green,115; blue,177}, shape=circle, inner sep=0pt, minimum size=0.15cm]
\tikzstyle{smod1}=[fill={rgb,255: red,150; green,255; blue,138}, draw={rgb,255: red,5; green,130; blue,3}, shape=circle, inner sep=0pt, minimum size=0.15cm]
\tikzstyle{smod2}=[fill={rgb,255: red,255; green,202; blue,239}, draw={rgb,255: red,126; green,0; blue,124}, shape=circle, inner sep=0pt, minimum size=0.15cm]
\tikzstyle{smod3}=[fill={rgb,255: red,255; green,213; blue,164}, draw={rgb,255: red,121; green,57; blue,2}, shape=circle, inner sep=0pt, minimum size=0.15cm]
\tikzstyle{mid-nc}=[fill=white, draw=black, shape=circle, inner sep=0pt, minimum size=0.25cm]

\tikzstyle{di}=[draw=black, ->]
\tikzstyle{da}=[-, dashed]
\tikzstyle{da-di}=[dashed, ->]
\tikzstyle{th}=[-, thick]
\tikzstyle{th-di}=[->, thick]
\tikzstyle{th-da}=[-, dashed, thick]
\tikzstyle{gr}=[-, draw={rgb,255: red,0; green,186; blue,143}, thick]
\tikzstyle{gr-di}=[draw={rgb,255: red,0; green,186; blue,143}, ->, thick]
\tikzstyle{gr-da}=[-, dashed, draw={rgb,255: red,0; green,186; blue,143}, thick]
\tikzstyle{pu}=[-, draw={rgb,255: red,124; green,95; blue,239}, thick]
\tikzstyle{pu-da}=[-, dashed, draw={rgb,255: red,124; green,95; blue,239}, thick]
\tikzstyle{re}=[-, draw={rgb,255: red,255; green,5; blue,80}, thick]
\tikzstyle{re-di}=[->, draw={rgb,255: red,255; green,5; blue,80}, thick]
\tikzstyle{re-da}=[-, draw={rgb,255: red,255; green,5; blue,80}, dashed, thick]
\tikzstyle{bl}=[-, draw={rgb,255: red,14; green,115; blue,177}, thick]
\tikzstyle{bl-di}=[draw={rgb,255: red,14; green,115; blue,177}, ->, thick]
\tikzstyle{bl-da}=[-, draw={rgb,255: red,14; green,115; blue,177}, dashed, thick]
\tikzstyle{lgr}=[-, draw={rgb,255: red,188; green,220; blue,156}, thick]
\tikzstyle{lgr-da}=[-, draw={rgb,255: red,188; green,220; blue,156}, thick, dashed]
\tikzstyle{lbl}=[-, draw={rgb,255: red,168; green,207; blue,245}, thick]
\tikzstyle{lbl-da}=[-, draw={rgb,255: red,168; green,207; blue,245}, thick, dashed]
\tikzstyle{pi}=[-, draw={rgb,255: red,236; green,181; blue,181}, thick]
\tikzstyle{pi-da}=[-, draw={rgb,255: red,236; green,181; blue,181}, thick, dashed]
\tikzstyle{sh}=[-, draw={rgb,255: red,171; green,171; blue,171}]
\tikzstyle{sh-da}=[-, draw={rgb,255: red,171; green,171; blue,171}, dashed]
\tikzstyle{fi-gr}=[-, fill={rgb,255: red,194; green,228; blue,162}, draw={rgb,255: red,194; green,228; blue,162}]
\tikzstyle{fi-bl}=[-, fill={rgb,255: red,175; green,215; blue,255}, draw={rgb,255: red,175; green,215; blue,255}]
\tikzstyle{fi-pi}=[-, fill={rgb,255: red,246; green,188; blue,188}, draw={rgb,255: red,246; green,188; blue,188}]
\tikzstyle{fi-pu}=[-, fill={rgb,255: red,230; green,203; blue,246}, draw={rgb,255: red,230; green,203; blue,246}]
\tikzstyle{fi-ye}=[-, fill={rgb,255: red,255; green,255; blue,155}, draw={rgb,255: red,255; green,255; blue,155}]
\tikzstyle{fi-dg}=[-, fill={rgb,255: red,152; green,170; blue,139}, draw={rgb,255: red,152; green,170; blue,139}]
\tikzstyle{fi-sh}=[-, fill={rgb,255: red,222; green,222; blue,222}, draw={rgb,255: red,222; green,222; blue,222}]
\tikzstyle{fi-dsh}=[-, fill={rgb,255: red,204; green,204; blue,204}, draw={rgb,255: red,204; green,204; blue,204}]
\tikzstyle{vt}=[-, very thick]
\tikzstyle{vt-di}=[->, very thick]
\tikzstyle{vt-da}=[-, dashed, very thick]
\tikzstyle{vt-gr}=[-, draw={rgb,255: red,0; green,186; blue,143}, very thick]
\tikzstyle{vt-gr-di}=[draw={rgb,255: red,0; green,186; blue,143}, ->, very thick]
\tikzstyle{vt-gr-da}=[-, dashed, draw={rgb,255: red,0; green,186; blue,143}, very thick]
\tikzstyle{vt-pu}=[-, draw={rgb,255: red,124; green,95; blue,239}, very thick]
\tikzstyle{vt-pu-da}=[-, dashed, draw={rgb,255: red,124; green,95; blue,239}, very thick]
\tikzstyle{vt-re}=[-, draw={rgb,255: red,255; green,5; blue,80}, very thick]
\tikzstyle{vt-re-di}=[->, draw={rgb,255: red,255; green,5; blue,80}, very thick]
\tikzstyle{vt-re-da}=[-, draw={rgb,255: red,255; green,5; blue,80}, dashed, very thick]
\tikzstyle{vt-bl}=[-, draw={rgb,255: red,14; green,115; blue,177}, very thick]
\tikzstyle{vt-bl-di}=[draw={rgb,255: red,14; green,115; blue,177}, ->, very thick]
\tikzstyle{vt-bl-da}=[-, draw={rgb,255: red,14; green,115; blue,177}, dashed, very thick]
\tikzstyle{vt-lgr}=[-, draw={rgb,255: red,188; green,220; blue,156}, very thick]
\tikzstyle{vt-lgr-di}=[draw={rgb,255: red,188; green,220; blue,156}, ->, very thick]
\tikzstyle{vt-lgr-da}=[-, draw={rgb,255: red,188; green,220; blue,156}, very thick, dashed]
\tikzstyle{vt-lbl}=[-, draw={rgb,255: red,168; green,207; blue,245}, very thick]
\tikzstyle{vt-lbl-di}=[->, draw={rgb,255: red,168; green,207; blue,245}, very thick]
\tikzstyle{vt-lbl-da}=[-, draw={rgb,255: red,168; green,207; blue,245}, very thick, dashed]
\tikzstyle{vt-lpp}=[-, draw={rgb,255: red,223; green,199; blue,240}, very thick]
\tikzstyle{vt-pi}=[-, draw={rgb,255: red,236; green,181; blue,181}, very thick]
\tikzstyle{vt-pi-da}=[-, draw={rgb,255: red,236; green,181; blue,181}, very thick, dashed]
\tikzstyle{vt-sh}=[-, fill=none, draw={rgb,255: red,171; green,171; blue,171}, very thick]
\tikzstyle{vt-sh-di}=[draw={rgb,255: red,171; green,171; blue,171}, very thick, ->]
\tikzstyle{fi-bx}=[-, fill={rgb,255: red,172; green,229; blue,232}, draw={rgb,255: red,14; green,115; blue,177}, thick]
\tikzstyle{ut-gr}=[-, draw={rgb,255: red,0; green,186; blue,143}, line width=3pt]
\tikzstyle{ut-gr-di}=[draw={rgb,255: red,0; green,186; blue,143}, ->, line width=3pt]
\tikzstyle{ut-ig}=[-, draw={rgb,255: red,255; green,5; blue,80}, line width=3pt]
\tikzstyle{ut-ig-di}=[draw={rgb,255: red,255; green,5; blue,80}, line width=3pt, ->]
\tikzstyle{fi-lbl}=[-, fill={rgb,255: red,213; green,237; blue,255}, draw={rgb,255: red,213; green,237; blue,255}]
\tikzstyle{fi-vlbl}=[-, fill={rgb,255: red,229; green,251; blue,255}, draw={rgb,255: red,229; green,251; blue,255}]
\tikzstyle{op-sh}=[-, draw={rgb,255: red,113; green,113; blue,113}, fill={rgb,255: red,113; green,113; blue,113}, opacity=0.5]

%% file: figures/GC0.tikz
\begin{tikzpicture}
	\begin{pgfonlayer}{nodelayer}
		\node [style=none] (0) at (0, 2) {};
		\node [style=none] (1) at (0, -2) {};
		\node [style=none] (2) at (0, -2.5) {$X$};
	\end{pgfonlayer}
	\begin{pgfonlayer}{edgelayer}
		\draw [style=vt-bl-di] (1.center) to (0.center);
	\end{pgfonlayer}
\end{tikzpicture}

%% file: figures/GC1.tikz
\begin{tikzpicture}
	\begin{pgfonlayer}{nodelayer}
		\node [style=none] (0) at (0, 2) {};
		\node [style=none] (1) at (0, -2) {};
		\node [style=none] (2) at (0, -2.5) {$X^{\vee}$};
	\end{pgfonlayer}
	\begin{pgfonlayer}{edgelayer}
		\draw [style=vt-bl-di] (0.center) to (1.center);
	\end{pgfonlayer}
\end{tikzpicture}

%% file: figures/GC17.tikz
\begin{tikzpicture}
	\begin{pgfonlayer}{nodelayer}
		\node [style=none] (0) at (0, -2) {};
		\node [style=none] (6) at (0, -2.5) {$A$};
		\node [style=none] (12) at (0, 2) {};
	\end{pgfonlayer}
	\begin{pgfonlayer}{edgelayer}
		\draw [style=vt-lgr-di] (0.center) to (12.center);
	\end{pgfonlayer}
\end{tikzpicture}

%% file: figures/GC27.tikz
\begin{tikzpicture}
	\begin{pgfonlayer}{nodelayer}
		\node [style=none] (0) at (0, 2) {};
		\node [style=none] (1) at (0, -2) {};
	\end{pgfonlayer}
	\begin{pgfonlayer}{edgelayer}
		\draw [style=vt-lgr] (0.center) to (1.center);
	\end{pgfonlayer}
\end{tikzpicture}

%% file: figures/LM15.tikz
\begin{tikzpicture}
	\begin{pgfonlayer}{nodelayer}
		\node [style=counit] (2) at (0, 0.5) {};
		\node [style=none] (4) at (0, -1.5) {};
		\node [style=none] (6) at (0, -2) {$A$};
		\node [style=none] (7) at (0, 2) {};
	\end{pgfonlayer}
	\begin{pgfonlayer}{edgelayer}
		\draw [style=vt-lgr-di, in=-90, out=90, looseness=0.75] (4.center) to (2);
	\end{pgfonlayer}
\end{tikzpicture}

%% file: figures/LM4.tikz
\begin{tikzpicture}
	\begin{pgfonlayer}{nodelayer}
		\node [style=none] (15) at (0, 2.5) {};
		\node [style=none] (16) at (0, -2.5) {};
		\node [style=none] (18) at (0, -3) {$X$};
	\end{pgfonlayer}
	\begin{pgfonlayer}{edgelayer}
		\draw [style=vt-bl-di] (16.center) to (15.center);
	\end{pgfonlayer}
\end{tikzpicture}